\documentstyle{article}

\def\aA{\mbox{$\mathcal A$}}

\def\T{\mbox{$\mathcal T$}}
\def\cirk{\,{\raisebox{.3ex}{\tiny $\circ$}}\,}
\def\kon{\wedge}
\def\str{\rightarrow}
\def\rts{\leftarrow}
\def\mj{\mbox{\bf 1}}
\def\koc{{\raisebox{-.2ex}{$\Box$}}}
\def\pl{\!+\!}
\def\od{\!-\!}
\def\prop#1#2{\vspace{2ex} \noindent{\sc #1.} {\it #2} \par \vspace{2ex}}
\def\dkz{\noindent{\sc Proof. }}
\def\qed{\hfill $\dashv$\vspace{2ex}}
\def\ks{\mbox{\footnotesize$\;\xi\;$}}
\def\kst{\raisebox{1pt}{\mbox{\tiny$\xi$}}}

\def\b#1#2{\stackrel{\raisebox{-2pt}{\mbox{\tiny $#1$}}}
{\raisebox{0pt}{$b$}}^{\raisebox{-7pt}{\scriptsize $#2$}}}

\def\d#1#2{\stackrel{\raisebox{-2pt}{\mbox{\tiny $\,#1$}}}
{\raisebox{0pt}{$\delta$}}^{\raisebox{-7pt}{\scriptsize $#2$}}}

\def\s#1#2{\stackrel{\raisebox{-2pt}{\mbox{\tiny $#1$}}}
{\raisebox{0pt}{$\sigma$}}^{\raisebox{-7pt}{\scriptsize $#2$}}}

\def\c#1{\stackrel{\raisebox{-2pt}{\mbox{\tiny $\,#1$}}}
{\raisebox{0pt}{$c$}}}

\def\p#1{\stackrel{\raisebox{-2pt}{\mbox{\tiny $\,#1$}}}
{\raisebox{0pt}{$\psi$}}}

\def\pM#1{\stackrel{\raisebox{-2pt}{\mbox{\tiny $\,\top$}}}
{\raisebox{0pt}{$\psi$}}^{\raisebox{-7pt}{\scriptsize $#1$}}}

\def\pp#1{\stackrel{\raisebox{-2pt}{\mbox{\tiny $\,#1$}}}
{\raisebox{0pt}{$\overline{\psi}$}}}

\def\ppM#1{\stackrel{\raisebox{-2pt}{\mbox{\tiny $\,\bot$}}}
{\raisebox{0pt}{$\overline{\psi}$}}^{\raisebox{-7pt}{\scriptsize
$#1$}}}

\def\w#1#2{\stackrel{\raisebox{-2pt}{\mbox{\tiny $#1$}}}
{\raisebox{0pt}{$w$}}^{\raisebox{-7pt}{\scriptsize $#2$}}}

\def\k#1#2{\stackrel{\raisebox{-2pt}{\mbox{\tiny $#1$}}}
{\raisebox{0pt}{$k$}}^{\raisebox{-7pt}{\scriptsize $#2$}}}

\def\ma{\mbox{$\mathbf{A}$}}
\def\mn{\mbox{$\mathbf{N}$}}
\def\mk{\mbox{$\mathbf{K}$}}
\def\mna{\mbox{$\mathbf{NA}$}}
\def\mka{\mbox{$\mathbf{KA}$}}
\def\mns{\mbox{$\mathbf{NS}$}}
\def\mks{\mbox{$\mathbf{KS}$}}
\def\ml{\mbox{$\mathbf{L}$}}

\def\Ck{\mbox{$\mathbf{C^k}$}}
\def\ACk{\mbox{$\mathbf{AC^k}$}}
\def\lz{|\![}
\def\dz{]\!|}

\def\nad#1#2{
\begin{picture}(14,12)(0,8)
\put(1,12){\line(1,0){9}} \put(3,15){$#1$} \put(3,5){$#2$}
\end{picture}}

\def\naad#1#2{
\begin{picture}(14,12)(4,8)
\put(2,12){\line(1,0){18}} \put(3,15){$\,\,#1$} \put(3,5){$#2$}
\end{picture}}

\def\nadd#1#2{
\begin{picture}(13,12)(3,8)
\put(1,12){\line(1,0){9}} \put(12,8){\scriptsize *}
\put(3,15){$#1$} \put(3,5){$#2$}
\end{picture}}

\def\ms{\mbox{$\mathbf{S}$}}
\def\SCk{\mbox{$\mathbf{SC^k}$}}
\def\Kon{\mbox{$\kon$}}

\begin{document}

\title{Intermutation}
\author{\small {\sc Kosta Do\v sen} and {\sc Zoran Petri\' c}
\\[1ex]
{\small Mathematical Institute, SANU}\\[-.5ex]
{\small Knez Mihailova 36, p.f. 367, 11001 Belgrade,
Serbia}\\[-.5ex]
{\small email: \{kosta, zpetric\}@mi.sanu.ac.rs}}
\date{}
\maketitle

\vspace{-3ex}

\begin{abstract}
\noindent This paper proves coherence results for categories with
a natural transformation called \emph{intermutation} made of
arrows from $(A\kon B)\vee(C\kon D)$ to ${(A\vee C)\kon(B\vee
D)}$, for $\kon$ and $\vee$ being two biendofunctors.
Intermutation occurs in iterated, or \emph{n}-fold, monoidal
categories, which were introduced in connection with \emph{n}-fold
loop spaces, and for which a related, but different, coherence
result was obtained previously by Balteanu, Fiedorowicz, Schw\"
anzl and Vogt. The results of the present paper strengthen up to a
point this previous result, and show that two-fold loop spaces
arise in the manner envisaged by these authors out of categories
of a more general kind, which are not two-fold monoidal in their
sense. In particular, some categories with finite products and
coproducts are such.

Coherence in Mac Lane's ``all diagrams commute'' sense is proved
here first for categories where for $\kon$ and $\vee$ one assumes
only intermutation, and next for categories where one also assumes
natural associativity isomorphisms. Coherence in the sense of
coherence for symmetric monoidal categories is proved when one
assumes moreover natural commutativity isomorphisms for $\kon$ and
$\vee$. A restricted coherence result, involving a proviso of the
kind found in coherence for symmetric monoidal closed categories,
is proved in the presence of two nonisomorphic unit objects. The
coherence conditions for intermutation and for the unit objects
are derived from a unifying principle, which roughly speaking is
about preservation of structures involving one endofunctor by
another endofunctor, up to a natural transformation that is not an
isomorphism. This is related to weakening the notion of monoidal
functor. A similar, but less symmetric, justification for
intermutation was envisaged in connection with iterated monoidal
categories. Unlike the assumptions previously introduced for
two-fold monoidal categories, the assumptions for the unit objects
of the categories of this paper, which are more general, allow an
interpretation in logic.

\end{abstract}

\vspace{.3cm}


\vspace{.5ex}

\noindent {\small {\it Keywords$\,$}: coherence, associativity,
commutativity, monoidal categories, symmetric monoidal categories,
iterated monoidal categories, loop spaces}

\section{\large\bf Introduction} For $\kon$ and $\vee$ being two
biendofunctors, we call \emph{intermutation} the natural
transformation $c^k$ whose components are the arrows
\[
c^k_{A_1,A_1',A_2,A_2'}\!:(A_1\kon A_1')\vee(A_2\kon
A_2')\str(A_1\vee A_2)\kon(A_1'\vee A_2')
\]
(the notation $c^k$ is used in \cite{DP04}). Intermutation has
been investigated in connection with \emph{n}-fold loop spaces in
\cite{BFSV}, where coherence of categories involving
intermutation, called \emph{iterated}, or \emph{n-fold},
\emph{monoidal categories}, is of central concern. (The wider
context of algebraic topology within which the results of
\cite{BFSV} should be placed is described in the introduction of
\cite{FSS}.) Here we take over from \cite{BFSV} the main coherence
conditions concerning intermutation, and strengthen the previous
coherence result, either by considering notions more general than
that of two-fold monoidal category with respect to unit objects,
or by adding symmetry, i.e.\ natural commutativity isomorphisms
for $\kon$ and $\vee$. Instead of having \emph{n-fold} for every
${n\geq 2}$, we deal only with the case when ${n=2}$, and leave
open the question to what an extent our approach could be extended
to ${n>2}$. We suppose however that this can be achieved by
relying on the technique of Section 14 of this paper.

Before \cite{BFSV}, intermutation was taken in \cite{JS93} to be
an isomorphism, which led to an isomorphism between $\kon$ and
$\vee$. In \cite{BFSV} it is not assumed that intermutation is an
isomorphism, and $\kon$ and $\vee$ need not be isomorphic, but the
two corresponding unit objects, $\top$ and $\bot$, are assumed to
be isomorphic (actually, they coincide), which delivers arrows
from ${A\vee B}$ to ${A\kon B}$. In our approach, $\top$ and
$\bot$ are not assumed to be isomorphic; we must only have an
arrow from $\bot$ to $\top$, and we need not have arrows from
${A\vee B}$ to ${A\kon B}$. This generalization of the notion of
iterated monoidal category does not sever the connection with loop
spaces. Two-fold loop spaces arise in the manner envisaged by
\cite{BFSV} out of categories of a more general kind, which are
not two-fold monoidal in the sense of \cite{BFSV} (see Sections
12, 13 and 15 for details). In particular some categories with
finite products and coproducts are such.

We differ also from \cite{BFSV} in being able to lift from the
coherence result in the associative nonsymmetric context, such as
the context of \cite{BFSV}, a proviso appropriate for symmetric
contexts. This proviso can be formulated either with the help of
graphs, or by limiting the number of occurrences of letters, as it
is done in \cite{BFSV}, and as we will do in the symmetric context
of Sections 14 and 16, where we have natural commutativity
isomorphisms for $\kon$ and $\vee$. Our coherence results in the
nonsymmetric context are results that say that all arrows with the
same source and target are equal. This is coherence in Mac Lane's
``all diagrams commute'' sense.

If $\kon$ and $\vee$ are interpreted as meet and join
respectively, and $\str$ is replaced by $\leq$, then intermutation
corresponds to an inequality that holds in any lattice. So if
$\kon$ and $\vee$ are interpreted as conjunction and disjunction
respectively, then intermutation corresponds to an implication
that is a logical law, whose converse is not a logical law. This
conjunction and this disjunction need not be those of classical or
intuitionistic logic, which are tied to distributive lattices.
They may be any \emph{lattice} conjunction and disjunction (see
Section 15), such as the \emph{additive} conjunction and
disjunction of linear logic. (Intermutation does not hold for the
\emph{multiplicative} conjunction and disjunction of linear
logic.) Our approach is in tune with this logical interpretation
also in the presence of the unit objects, which is not the case
for \cite{BFSV}.

The law generalizing intermutation in logic is \emph{quantifier
shift}:
\[
\exists x\forall yA\str\forall y\exists xA.
\]
When $A$ is $xRy$ and ${R\subseteq \{0,1\}^2}$, quantifier shift
amounts to intermutation:
\[
(0R0\kon 0R1)\vee(1R0\kon 1R1)\str (0R0\vee 1R0)\kon(0R1\vee 1R1).
\]
Other generalizations of intermutation in logic are the laws
\[
\begin{array}{l}
(\forall xA\vee\forall xB)\str\forall x(A\vee B),\\
\exists x(A\kon B)\str(\exists xA\kon\exists xB).
\end{array}
\]
Analogously, the logical laws
\[
\begin{array}{l}
\forall x\forall yA\str\forall y\forall xA,\\[.5ex]
\exists x\exists yA\str\exists y\exists xA
\end{array}
\]
generalize respectively the implications underlying the arrows
\[
\begin{array}{l}
\hat{c}^m_{A_1,A_1',A_2,A_2'}\!:(A_1\kon A_1')\kon(A_2\kon
A_2')\str(A_1\kon A_2)\kon(A_1'\kon A_2'),\\[.5ex]
\check{c}^m_{A_1,A_1',A_2,A_2'}\!:(A_1\vee A_1')\vee(A_2\vee
A_2')\str(A_1\vee A_2)\vee(A_1'\vee A_2'),
\end{array}
\]
which we will encounter in Section 14. (In \cite{DP06a}, the
arrows $\hat{c}^m$ are investigated in the same spirit as
intermutation here.)

In a similar sense, the logical laws
\[
\begin{array}{l}
A\kon\exists xB\str\exists x(A\kon B),\\[.5ex]
A\vee \forall xB\str \forall x(A\vee B),
\end{array}
\]
provided $x$ is not free in $A$, together with the converse
implications, generalize respectively distribution of conjunction
over disjunction and distribution of disjunction over conjunction;
the logical law ${\forall xA\str\exists xA}$ generalizes the
implication ${A\kon B\str A\vee B}$ (which is the type of the
\emph{mix} arrows of \cite{DP04}, Chapter~8).

By restricting the quantifiers in quantifier shift, we obtain the
logical laws
\[
\begin{array}{c}
(\exists x\in\emptyset)\forall yA\str\forall
y(\exists x\in\emptyset)A,\\[.5ex]
\exists x(\forall y\in\emptyset)A\str(\forall
y\in\emptyset)\exists xA,\\[.5ex]
(\exists x\in\emptyset)(\forall
y\in\emptyset)A\str(\forall y\in\emptyset)(\exists
x\in\emptyset)A.
\end{array}
\]
These laws correspond respectively to the arrows
\[
\begin{array}{rl}
\hat{w}^{\rts}_{\bot}\!\!\!\!&:\bot\str\bot\kon\bot,\\[.5ex]
\check{w}^{\str}_{\top}\!\!\!\!&:\top\vee\top\str\top,\\[.5ex]
\kappa\!\!\!\!&:\bot\str\top,
\end{array}
\]
which we will encounter for the first time in Section~3, and which
will play an important role in Sections 12 and 16 especially.
These assumptions for $\top$ and $\bot$ flow out of a unifying
logical principle, which delivers also intermutation.

As $\hat{c}^m$ serves to \emph{atomize} the conjunctive indices of
diagonal arrows, so intermutation serves to atomize the
disjunctive indices of diagonal arrows, or the conjunctive indices
of codiagonal arrows (see the proof of the Proposition in Section
15). Diagonal and codiagonal arrows correspond in proof theory to
the structural rules of contraction on the left and on the right
respectively. This atomization was exploited in a
proof-theoretical, not categorial, context in \cite{BT01},
\cite{B03} and \cite{B06}, where intermutation is called
\emph{medial} (an unfortunate denomination, since the different
principle underlying $\hat{c}^m$ is called so in universal
algebra; see \cite{JK83}). Following these papers, a categorial
investigation of intermutation was started in \cite{L06} and
papers cited therein, where one finds various proposals for
axiomatizing structures involving more than what we consider,
without concentrating on coherence. The notation of \cite{L06}
(Section 2.3), analogous to that which may be found in \cite{LS97}
(Session 26), resembles the rectangular notation of Section~8
below. Intermutation plays an important role in the distributive
lattice categories of \cite{DP04} (Chapter 11; see also Sections
9.4 and 13.2).

Here is a summary of our paper. After some preliminary matters in
Section~2, in Sections 3 and 4 we justify the introduction of
intermutation and of the equations for arrows involving it. We do
the same for the arrows $\hat{w}^{\rts}_{\bot}$,
$\check{w}^{\str}_{\top}$ and $\kappa$ mentioned above. This
justification is governed by a unifying principle, like the
principle of \cite{DP06a}, and it is related to the justification
provided by \cite{BFSV}, which is however less symmetrical.
Roughly speaking, this principle is about preservation of
structures involving one endofunctor by another endofunctor, up to
a natural transformation that is not an isomorphism. This has to
do with weakening the notion of monoidal functor (see \cite{EK66},
Sections II.1 and III.1, \cite{JS93}, \cite{ML71}, second edition,
Section XI.2, and \cite{DP04}, Section 2.8).

Sections 5-7 present auxiliary coherence results involving the
unit objects $\top$ and $\bot$. In Sections 8-9 we prove coherence
for intermutation in the absence of additional assumptions
concerning the biendofunctors $\kon$ and $\vee$. In Sections 10-12
we prove our central coherence result for intermutation in the
presence of natural associativity isomorphisms for $\kon$ and
$\vee$. In Section 11, in the absence of the unit objects, we have
a full coherence result in Mac Lane's ``all diagrams commute''
sense, and in Section 12, in the presence of the unit objects, we
have a restricted coherence result. The restriction is of the kind
Kelly and Mac Lane had for their coherence result for symmetric
monoidal closed categories in \cite{KML71}. Section 10 is about
deciding whether there is an arrow with given source and target,
which is a problem some authors take as being a part of the
coherence problem (cf.\ \cite{KML71}, Theorem 2.1, and
\cite{BFSV}, Theorem 3.6.2). In Section 13 we compare the
coherence results of Sections 11-12 with the related, but
different, coherence result of \cite{BFSV}. We show that our
restricted coherence result of Section 12 is sufficient for the
needs of \cite{BFSV} when $n$ in \emph{n-fold} is 2. We leave open
the question whether this result can be extended to ${n>2}$. As we
said above, we suppose that this can be achieved by relying on the
technique of Section 14.

In Sections 14-16 we prove coherence for intermutation in the
presence of symmetry, i.e.\ natural commutativity isomorphisms for
$\kon$ and $\vee$, besides natural associativity isomorphisms. In
Section 14 we have a full coherence result in the absence of the
unit objects, and in Section 16 a restricted coherence result in
the presence of the unit objects, the restriction being analogous
to the restriction of Section 12. We formulate coherence in the
presence of symmetry without mentioning graphs, but by limiting
the number of occurrences of letters. This is however equivalent
to coherence with respect to graphs. In Section 15 we show that
with  natural associativity and commutativity isomorphisms
together with intermutation we have caught an interesting fragment
of categories with finite nonempty products and coproducts. This
fact can serve to obtain loop spaces in the style of \cite{BFSV}
out of categories not envisaged by \cite{BFSV}.

\section{\large\bf Biassociative and biunital categories} This
section is about preliminary matters. In it we fix
terminology and state some basic results on which we rely.

For an arrow ${f\!:A\str B}$ in a category, the \emph{type} of $f$
is ${A\str B}$, which stands for the ordered pair ${(A,B)}$ made
of the source and target of $f$. We call \emph{categorial}
equations the following usual equations assumed for categories:
\begin{tabbing}
\mbox{\hspace{8.5em}}\=$f\cirk \mj_A=\mj_B\cirk f=f,
\quad{\mbox{\rm for}}\; f\!:A\str B,$
\\*[1ex]
\>$h\cirk (g\cirk f)=(h\cirk g)\cirk f.$
\end{tabbing}
We call \emph{bifunctorial} equations for $\!\ks\!$ the equations
\begin{tabbing}
\mbox{\hspace{8.5em}}\=$f\cirk \mj_A=\mj_B\cirk f=f,
\quad{\mbox{\rm for}}\; f\!:A\str B,$\kill

\>$\mj_A\ks\mj_B=\mj_{A\xi B},$
\\*[1ex]
\>$(f_1\cirk f_1')\ks (f_2\cirk f_2')=(f_1\ks f_2)\cirk (f_1'\ks
f_2').$
\end{tabbing}
The \emph{naturality} equation for $c^k$ is
\[
((f\vee h)\kon(g\vee j))\cirk
c^k_{A,B,C,D}=c^k_{A',B',C',D'}\cirk((f\kon g)\vee(h\kon j)).
\]
We have analogous naturality equations for other natural
transformations to be encountered in the text.

For ${n\geq 0}$, an \emph{n-endofunctor} of a category \aA\ is a
functor from the product category $\aA^n$ to \aA. If ${n=0}$, then
$\aA^0$ is the trivial category \T\ with a unique object $\ast$,
and a unique arrow $\mj_\ast$, and 0-endofunctors of \aA\ amount
to special objects of \aA. \emph{Endofunctors} are 1-endofunctors,
and \emph{biendofunctors} are 2-endofunctors.

Next we introduce some classes of categories for which coherence
results are already known. These results are ultimately based on
Mac Lane's monoidal coherence results of \cite{ML63} (see also
\cite{ML71}, Section VII.2).

We say that ${\langle\aA,\kon,\vee\rangle}$ is a
\emph{biassociative} category when \aA\ is a category, $\kon$ and
$\vee$ are biendofunctors of \aA, and there are two natural
isomorphisms $\b{\xi}{\str}$, for ${\!\ks\!\in\{\kon,\vee\}}$,
with the following components in \aA:
\[
\b{\xi}{\str}_{A,B,C}:A\ks(B\ks C)\str(A\ks B)\ks C,
\]
which satisfy Mac Lane's pentagonal equations:
\[
\b{\xi}{\str}_{A\kst B,C,D}\!\cirk\b{\xi}{\str}_{A,B,C\kst
D}=(\b{\xi}{\str}_{A,B,C}\!\ks\mj_D)\cirk\b{\xi}{\str}_{A,B\kst
C,D}\!\!\cirk(\mj_A\ks\!\b{\xi}{\str}_{B,C,D}).
\]
The natural isomorphism inverse to $\b{\xi}{\str}$ is
$\b{\xi}{\rts}$. We call the $\b{\xi}{\str}$ and $\b{\xi}{\rts}$
arrows collectively \emph{b-arrows}.

Let \ma\ be the free biassociative category generated by a set of
objects. (This set may be conceived as a discrete category.) We
take that the objects of \ma\ are the \emph{formulae} of the
propositional language generated by a set of \emph{letters}
(nonempty if the category is to be interesting) with $\kon$ and
$\vee$ as binary connectives. We will later use for letters $p$,
$q$, $r,\ldots$, sometimes with indices. Formally, we have the
inductive definition:
\begin{tabbing}
\hspace{2.5em}\=every letter is a formula;\\*[.5ex] \>if $A$ and
$B$ are formulae, then ${(A\ks B)}$ is a formula, for
${\!\ks\!\in\{\kon,\vee\}}$.
\end{tabbing}
As usual, we take the outermost parentheses of formulae for
granted, and omit them. We do the same for other expressions of
the same kind later on. The formulae $A$ and $B$ are the
\emph{conjuncts} of ${A\kon B}$, and the \emph{disjuncts} of
${A\vee B}$. For the free biunital category below, and other
categories that have the special objects
${\zeta\in\{\top,\bot\}}$, we enlarge the inductive definition of
formula by the clauses: ``$\zeta$ is a formula.''

The \emph{arrow terms} of \ma\ are defined by assuming first that
$\mj_A$, $\b{\xi}{\str}_{A,B,C}$ and $\b{\xi}{\rts}_{A,B,C}$ are
arrow terms for all formulae $A$, $B$ and $C$. These
\emph{primitive} arrow terms are then closed under composition
$\cirk$ and the operations $\!\ks\!$, provided for composition
that the types of the arrow terms composed make them composable.

These arrow terms are then subject to the equations assumed for
biassociative categories; namely, the categorial equations, the
bifunctorial equations for $\kon$ and $\vee$, the naturality and
isomorphism equations for $\b{\xi}{\str}$ and $\b{\xi}{\rts}$, and
Mac Lane's pentagonal equations. This means that to obtain \ma\ we
factor the arrow terms through an equivalence relation engendered
by the equations, which is congruent with respect to $\cirk$ and
$\!\ks\!$, and we take the equivalence classes as arrows. (A
detailed formal definition of such syntactically constructed
categories may be found in \cite{DP04}, Chapter 2.) We proceed
analogously for other freely generated categories we deal with
later in the text.

An arrow term of \ma\ in which $\cirk$ does not occur, and in
which $\b{\xi}{\str}$ occurs exactly once is called a
$\b{\xi}{\str}$-\emph{term}. For example, $\b{\xi}{\str}_{A,B,C}$,
${\mj_D\kon\b{\xi}{\str}_{A,B,C}}$ and
${(\mj_D\kon\b{\xi}{\str}_{A,B,C})\vee\mj_E}$ are all
$\b{\xi}{\str}$-terms. We define analogously
$\b{\xi}{\rts}$-terms, and other $\beta$-terms for other natural
transformations $\beta$, which will be introduced later in the
text. The subterm of a $\beta$-term that is a component of the
natural transformation $\beta$ is called its \emph{head}. For
example, the head of the $\hat{b}^{\str}$-term
${(\mj_D\kon\hat{b}^{\str}_{A,B,C})\vee\mj_E}$ is
$\hat{b}^{\str}_{A,B,C}$.

An arrow term of the form ${f_n\cirk\ldots\cirk f_1}$, where
${n\geq 1}$, such that for every ${i\in\{1,\ldots,n\}}$ we have
that $f_i$ is composition-free is called \emph{factorized}, and
$f_i$ in such a factorized arrow term is called a \emph{factor}. A
factorized arrow term ${f_n\cirk\ldots\cirk f_1\cirk\mj_A}$ is
\emph{developed} when for every ${i\in\{1,\ldots,n\}}$ we have
that $f_i$ is a $\beta$-term for some $\beta$. Then by using the
categorial and bifunctorial equations we prove easily by induction
on the length of $f$ the following lemma for the category \ma.

\prop{Development Lemma}{For every arrow term $f$ there is a
developed arrow term $f'$ such that ${f=f'}$.}

\noindent The same lemma will hold for other freely generated
categories analogous to \ma\, which we will introduce later.

The category \ma\ is a \emph{preorder}, which means that there is
in \ma\ at most one arrow of a given type, i.e.\ with a given
source and target (for a proof see \cite{DP04}, Section 6.1). We
call this fact \emph{Biassociative Coherence}.

We say that ${\langle\aA,\kon,\vee,\top,\bot\rangle}$ is a
\emph{biunital} category when \aA\ is a category, $\kon$ and
$\vee$ are biendofunctors of \aA, and $\top$ and $\bot$ are
special objects of \ma\ such that there are four natural
isomorphisms $\d{\xi}{\str}$ and $\s{\xi}{\str}$, for
${\!\ks\!\in\{\kon,\vee\}}$, with the following components in \aA:
\[
\d{\xi}{\str}_A:A\ks\zeta\str A, \quad\quad\quad
\s{\xi}{\str}_A:\zeta\ks A\str A,
\]
for ${(\!\ks\!,\zeta)\in\{(\kon,\top),(\vee,\bot)\}}$, which
satisfy the equations
\[
\d{\xi}{\str}_\zeta\,=\;\s{\xi}{\str}_\zeta\!\!.
\]
The natural isomorphisms inverse to $\d{\xi}{\str}$ and
$\s{\xi}{\str}$ are $\d{\xi}{\rts}$ and $\s{\xi}{\rts}$
respectively. We call all these arrows collectively
\emph{$\delta$-$\sigma$-arrows}.

A coherence result analogous to Biassociative Coherence can be
proved for biunital categories (in the style of Normal Biunital
Coherence of Section~5 below). We will not dwell on that proof,
which we do not need for the rest of this work.

A \emph{bimonoidal} category is a biunital category
${\langle\aA,\kon,\vee,\top,\bot\rangle}$ such that
${\langle\aA,\kon,\vee\rangle}$ is a biassociative category, and,
moreover, the following equations are satisfied:
\[
\b{\xi}{\str}_{A,\zeta,C}\,=\;\d{\xi}{\rts}_A\ks\s{\xi}{\str}_C
\]
for ${(\!\ks\!,\zeta)\in\{(\kon,\top),(\vee,\bot)\}}$. A coherence
result analogous to Biassociative Coherence can be proved for
bimonoidal categories (see \cite{DP04}, Section 6.1).

\section{\large\bf Upward and downward functors and intermutation}
In this section we justify the introduction of the arrows $c^k$,
$\hat{w}^{\rts}_{\bot}$, $\check{w}^{\str}_{\top}$ and $\kappa$
(see Section~1), and of the $\delta$-$\sigma$-arrows of biunital
categories (see the preceding section), in terms of notions of
functors that preserve the structure induced by an $m$-endofunctor
up to a natural transformation, which need not be an isomorphism.
This justification proceeds out of a unifying principle.

Let $\vec{A}_m$, where ${m\geq 0}$, be an abbreviation for
${A_1,\ldots,A_m}$. If ${m=0}$, then $\vec{A}_m$ is the unique
object $\ast$ of the trivial category \T\ (see the preceding
section).

Let $M$ and $M'$ be $m$-endofunctors of the categories \aA\ and
$\aA'$ respectively, and let $F$ be a functor from $\aA'$ to \aA.
We say that ${(F,\psi)}$ is an \emph{upward} functor from
${(\aA',M')}$ to ${(\aA,M)}$ when $\psi$ is a natural
transformation whose components are the following arrows of \aA:
\[
\psi_{\vec{A}_m}\!\!:M(\overrightarrow{F\!A\,}\!_m)\str
FM'(\vec{A}_m).
\]
We say that ${(F,\overline{\psi})}$ is a \emph{downward} functor
from ${(\aA',M')}$ to ${(\aA,M)}$ when $\overline{\psi}$ is a
natural transformation whose components are the following arrows
of \aA:
\[
\overline{\psi}_{\vec{A}_m}\!\!:FM'(\vec{A}_m)\str
M(\overrightarrow{F\!A\,}\!_m),
\]
with type converse to that of $\psi_{\vec{A}_m}$.

When ${(F,\psi)}$ is an upward functor and ${(F,\overline{\psi})}$
a downward functor from ${(\aA',M')}$ to ${(\aA,M)}$, and,
moreover, $\psi$ is an isomorphism whose inverse is
$\overline{\psi}$, we say that ${(F,\psi,\overline{\psi})}$ is a
\emph{loyal} functor from ${(\aA',M')}$ to ${(\aA,M)}$.

For an $m$-endofunctor $M$ of a category \aA, we define as usual,
in a coordinatewise manner, the $m$-endofunctor $M^n$ of the
product category $\aA^n$; this means that on objects we have
\[
M^n((A^1_1,{\small \ldots},A^n_1),{\small \ldots},(A^1_m,{\small
\ldots},A^n_m))=_{df}(M(A^1_1,{\small \ldots},A^1_m),{\small
\ldots}, M(A^n_1,{\small \ldots},A^n_m)),
\]
and analogously on arrows. Suppose now that $\aA'$ is the product
category $\aA^n$, and that $M'$ is $M^n$.

We say that $M$ \emph{intermutes} with $F$ when there is an upward
functor ${(F,\p{M})}$ from ${(\aA^n,M^n)}$ to ${(\aA,M)}$. From an
upward functor ${(F,\p{M})}$ from ${(\aA^n,M^n)}$ to ${(\aA,M)}$
we obtain a downward functor ${(M,\pp{F})}$ from ${(\aA^m,F^m)}$
to ${(\aA,F)}$ such that
\[
\pp{F}_{(A^1_1,\ldots,A^1_m),\ldots,(A^n_1,\ldots,A^n_m)}=\;
\p{M}_{(A^1_1,\ldots,A^n_1),\ldots,(A^1_m,\ldots,A^n_m)},
\]
and vice versa. So $M$ intermutes with $F$ iff there is a downward
functor ${(M,\pp{F})}$ from ${(\aA^m,F^m)}$ to ${(\aA,F)}$. (Note
that the relation ``intermutes with'' is not symmetric.)

We have that ${(F,\p{M},\pp{M})}$ is a loyal functor iff
${(M,\p{F},\pp{F})}$ is a loyal functor. If ${(F,\p{M},\pp{M})}$
is a loyal functor, then $M$ and $F$ intermute with each other.

For $\kon$ and $\vee$ being biendofunctors of the category \aA, if
$\vee$ intermutes with $\kon$, then we write
\[
c^k_{A_1,A_1',A_2,A_2'}\!:(A_1\kon A_1')\vee(A_2\kon
A_2')\str(A_1\vee A_2)\kon(A_1'\vee A_2')
\]
for ${\p{\vee}_{(A_1,A_1'),(A_2,A_2')}}$, which is equal to
${\pp{\kon}_{(A_1,A_2),(A_1',A_2')}}$. We call the natural
transformation $c^k$ \emph{intermutation}.

For $\!\ks\!$ being a biendofunctor of \aA\ and $\zeta$ a special
object of \aA, if ${(\!\ks\!,\p{\zeta})}$ is an upward functor
from ${(\aA^2,(\zeta,\zeta))}$ to ${(\aA,\zeta)}$, then we write
\[
\w{\xi}{\rts}_\zeta:\zeta\str\zeta\ks\zeta
\]
for $\p{\zeta}_\ast$, which is equal to $\pp{\xi}_{\ast,\ast}$;
here $\zeta$ intermutes with $\!\ks\!$. If
${(\!\ks\!,\pp{\zeta})}$ is a downward functor from
${(\aA^2,(\zeta,\zeta))}$ to ${(\aA,\zeta)}$, then we write
\[
\w{\xi}{\str}_\zeta:\zeta\ks\zeta\str\zeta
\]
for $\pp{\zeta}_\ast$, which is equal to $\p{\xi}_{\ast,\ast}$;
here $\!\ks\!$ intermutes with $\zeta$. If
${(\!\ks\!,\p{\zeta},\pp{\zeta})}$ is a loyal functor from
${(\aA^2,(\zeta,\zeta))}$ to ${(\aA,\zeta)}$, then $\!\ks\!$ and
$\zeta$ intermute with each other.

For $\top$ and $\bot$ being special objects of the category \aA,
if $\bot$ intermutes with $\top$, then we write
\[
\kappa\!:\bot\str\top
\]
for $\p{\vee}_\ast$, which is equal to $\pp{\kon}_\ast$.
Intuitively, we conceive of $\top$ as nullary $\kon$, and of
$\bot$ as nullary $\vee$.

In terms of the biendofunctor $\!\ks\!$ and the special object
$\zeta$ of \aA\ we define the endofunctors $\!\ks\!\zeta$ and
$\zeta\!\ks\!$ of \aA\ by
\[
(\!\ks\!\zeta)a=_{df}a\ks\zeta,\quad\quad\quad\quad
(\zeta\!\ks\!)a=_{df}\zeta\ks a;
\]
here $a$ stands either for an object or for an arrow of \aA, and
for arrows we read $\zeta$ on the right-hand sides as $\mj_\zeta$.
Let $I$ be the identity functor of \aA. If
${(I,\psi,\overline{\psi})}$ is a loyal functor from
${(\aA,\kon\top)}$ to ${(\aA,I)}$, then we write
\[
\hat{\delta}^{\str}_A\!:A\kon\top\str A \quad \mbox{\rm and} \quad
\hat{\delta}^{\rts}_A\!:A\str A\kon\top
\]
for $\psi_A$ and $\overline{\psi}_A$ respectively. A loyal functor
from ${(\aA,\top\kon)}$ to ${(\aA,I)}$ yields analogously the
isomorphisms
\[
\hat{\sigma}^{\str}_A\!:\top\kon A\str A \quad \mbox{\rm and}
\quad \hat{\sigma}^{\rts}_A\!:A\str\top\kon A.
\]
By replacing $\kon$ and $\top$ by $\vee$ and $\bot$ respectively,
we obtain analogously the remaining $\delta$-$\sigma$-arrows. The
arrows $\d{\xi}{\str}_\zeta$ and $\d{\xi}{\rts}_\zeta$, for
${(\!\ks\!,\zeta)\in\{(\kon,\top),(\vee,\bot)\}}$, guarantee that
$\!\ks\!$ and $\zeta$ intermute with each other.

The introduction of the $b$-arrows of biassociative categories
(see the preceding section) is justified in a similar manner in
\cite{DP06a}.

\section{\large\bf Preservation}
In the preceding section we justified the introduction of the
arrows $c^k$, $\hat{w}^{\rts}_{\bot}$, $\check{w}^{\str}_{\top}$
and $\kappa$, and of the $\delta$-$\sigma$-arrows. In this section
we will justify the equations we will assume later for these
arrows. This justification (like that of \cite{DP06a}) proceeds
out of a unifying principle of preservation of a natural
transformation by an $m$-endofunctor, which resembles something
that may be found in \cite{BFSV}, but is more symmetrical.

First, we define inductively the notion of \emph{shape} and of its
arity:
\begin{itemize}
\item[(0)]$\zeta$ is a shape of arity 0;\vspace{-1ex}
\item[(1)]$\koc$ is a shape of arity 1;\vspace{-1ex} \item[(2)] if
$M$ and $N$ are shapes of arities $m$ and $n$ respectively, then
${(M\ks N)}$ is a shape of arity ${m\pl n}$.
\end{itemize}
As we did for formulae, we take the outermost parentheses of
shapes for granted, and omit them. The shapes we have just defined
will be called ${(\!\ks\!,\zeta)}$-\emph{shapes}. A
$\!\ks\!$-\emph{shape} has the clause (0) omitted.

Let $a_i$, where $1\leq i\leq m$, stand either for an object or
for an arrow of \aA. We shall use the following abbreviations,
like the abbreviation $\vec{A}_m$ of the preceding section:
\vspace{-1.5ex}
\begin{tabbing}
\mbox{\hspace{8em}}\=$\vec{a}_m$ \mbox{\hspace{3em}}\= for \quad
${a_1,\ldots,a_m}$,\\[.5ex]
\>$\vec{a}_{\pi(m)}$ \> for \quad
${a_{\pi(1)},\ldots,a_{\pi(m)}}$,\\[.5ex]
\>$\overrightarrow{a_m, a_m'\,}$ \> for \quad ${a_1,
a_1',\ldots,a_m, a_m'}$,\\[.5ex]
\>$\overrightarrow{a_m\kon a_m'\,}$ \> for \quad ${a_1\kon
a_1',\ldots,a_m\kon a_m'}$,
\end{tabbing}
and other analogous abbreviations, made on the same pattern. If
$m=0$ and $a_i$ stands for an arrow, then $\vec{a}_m$ is the
unique arrow $\mj_\ast$ of the trivial category \T.

A shape $M$ of arity $m$ defines an $m$-endofunctor of \aA, such
that $M(\vec{a}_m)$ is obtained by putting $a_i$ for the $i$-th
$\koc$, counting from the left, in the shape $M$. For arrows, we
read $\zeta$ as $\mj_\zeta$. An $m$-endofunctor defined by a shape
$M$ and a permutation $\pi$ of $\{1,\ldots,m\}$ define an
$m$-endofunctor $M^\pi$ such that
\[
M^\pi(\vec{a}_m)=_{df}M(\vec{a}_{\pi(m)}).
\]

For the definitions below we make the following assumptions:
\begin{tabbing}
\hspace{2.5em}\=${(\vee\kon)}$\quad\quad\=$\vee\,$ \=intermutes
with
$\kon$,\\[.5ex]
\>${(\bot\kon)}$\>$\bot$ \>intermutes with
$\kon$,\\[.5ex]
\>${(\vee\top)}$\>$\vee$ \>intermutes with $\top$.
\end{tabbing}
The assumption ${(\vee\kon)}$ delivers the natural transformation
$c^k$ with components in \aA, the assumption ${(\bot\kon)}$
delivers the arrow ${\hat{w}^\rts_\bot\!:\bot\str\bot\vee\bot}$ of
\aA, and the assumption ${(\vee\top)}$ delivers the arrow
${\check{w}^\str_\top\!:\top\vee\top\str\top}$ of \aA\ (see the
preceding section).

We define now by induction on the complexity of the
${(\vee,\bot)}$-shape $M$ of arity $m$ the natural transformation
$\psi^M\!$ whose components are the following arrows of \aA:
\[
\psi^M_{\scriptsize\overrightarrow{A_m,A_m'}}\!\!:M(\overrightarrow{A_m\kon
A_m'\,})\str M(\vec{A}_m)\kon M(\vec{A}_m').
\]
Here is the definition:
\begin{tabbing}
\mbox{\hspace{2em}}\=$\psi^\bot=\hat{w}^\rts_\bot\!:\bot\str\bot\kon\bot$,\\[1.5ex]
\>$\psi^\Box_{A,A'}=\mj_{A\kon A'}$,\\[2ex]
\>$\psi^{M\vee
N}_{\scriptsize\overrightarrow{A_m,A_m'},\overrightarrow{B_n,B_n'\,}}=
c^k_{M(\vec{A}_m),M(\vec{A}_m'),N(\vec{B}_n),N(\vec{B}_n')}\cirk
(\psi^M_{\scriptsize\overrightarrow{A_m,A_m'}}\vee
\psi^N_{\scriptsize\overrightarrow{B_n,B_n'\,}})$.
\end{tabbing}
In the last clause, $M$ is a shape of arity $m$ and $N$ a shape of
arity $n$.

We define next by induction on the complexity of the
${(\kon,\top)}$-shape $M$ of arity $m$ the natural transformation
$\overline{\psi}^M\!$ whose components are the following arrows of
\aA:
\[
\overline{\psi}^M_{\scriptsize\overrightarrow{A_m,A_m'}}\!\!:M(\vec{A}_m)\vee
M(\vec{A}_m') \str M(\overrightarrow{A_m\vee A_m'\,}).
\]
Here is the definition:
\begin{tabbing}
\mbox{\hspace{2em}}\=$\overline{\psi}^\top=\check{w}^\str_\top\!:\top\vee\top\str\top$,\\[1.5ex]
\>$\overline{\psi}^\Box_{A,A'}=\mj_{A\vee A'}$,\\[2ex]
\>$\overline{\psi}^{M\kon
N}_{\scriptsize\overrightarrow{A_m,A_m'},\overrightarrow{B_n,B_n'\,}}=
 (\overline{\psi}^M_{\scriptsize\overrightarrow{A_m,A_m'}}\kon
\overline{\psi}^N_{\scriptsize\overrightarrow{B_n,B_n'\,}})\cirk
c^k_{M(\vec{A}_m),N(\vec{B}_n),M(\vec{A}_m'),N(\vec{B}_n')}$.
\end{tabbing}
Note that
\[
\psi^{\Box\vee\Box}_{A,A',B,B'}=\overline{\psi}^{\Box\kon\Box}_{A,B,A',B'}=c^k_{A,A',B,B'}.
\]
If $\bot$ and $\top$ do not occur in the shape $M$, we do not need
the arrows $\hat{w}^\rts_\bot$ and $\check{w}^\str_\top$ to define
$\psi^M$ and $\overline{\psi}^M$, and we may do without the
assumptions ${(\bot\kon)}$ and ${(\vee\top)}$.

Let now $\alpha$ be a natural transformation from the
$m$-endofunctor of \aA\ defined by the ${(\vee,\bot)}$-shape $M_1$
to the $m$-endofunctor of \aA\ defined by the
${(\vee,\bot)}$-shape $M_2$ and a permutation $\pi$ of
$\{1,\ldots,m\}$. We say that $\alpha$ is \emph{upward preserved}
by $\kon$ when diagrams of the following form commute in \aA:

\begin{center}
\begin{picture}(300,105)
\put(3,80){$M_1(\vec{A}_m)\kon M_1(\vec{A}_m')$}

\put(180,80){$M_2(\vec{A}_{\pi(m)})\kon M_2(\vec{A}_{\pi(m)}')$}

\put(108,93){$\alpha_{\vec{A}_m}\!\kon\alpha_{\vec{A}_m'}$}

\put(92,85){\vector(1,0){80}}

%%%%%%%%%%%

\put(7,48){$\psi^{M_1}_{\scriptsize\overrightarrow{A_m,A_m'}}$}

\put(237,48){$\psi^{M_2}_{\scriptsize\overrightarrow{A_{\pi(m)},A_{\pi(m)}'}}$}

\put(46,30){\vector(0,1){40}}

\put(233,30){\vector(0,1){40}}

%%%%%%%%%%%%%

\put(8,10){$M_1(\overrightarrow{A_m\kon A_m'})$}

\put(184,10){$M_2(\overrightarrow{A_{\pi(m)}\kon A_{\pi(m)}'})$}

\put(110,6){$\alpha_{\scriptsize\overrightarrow{A_m\kon A_m'}}$}

\put(92,15){\vector(1,0){80}}

\end{picture}
\end{center}
i.e., we have in \aA\ the equation
\begin{tabbing}
\mbox{\hspace{2.5em}}$(\psi\alpha)\quad\quad
\psi^{M_2}_{\scriptsize\overrightarrow{A_{\pi(m)},A_{\pi(m)}'}}\cirk
\alpha_{\scriptsize\overrightarrow{A_m\kon
A_m'}}=(\alpha_{\vec{A}_m}\kon\alpha_{\vec{A}_m'})\cirk
\psi^{M_1}_{\scriptsize\overrightarrow{A_m,A_m'}}$.
\end{tabbing}
(This equation is an instance of the equation $(\psi\alpha)$ of
\cite{DP04}, Section 2.8; see also \cite{DP06a}, Section~2.
Something analogous, but different from what we have, occurs in
\cite{BFSV}, Section~1.)

Let now $\beta$ be a natural transformation as $\alpha$ above save
that the shapes $M_1$ and $M_2$ are not ${(\vee,\bot)}$-shapes but
${(\kon,\top)}$-shapes. We say that $\beta$ is \emph{downward
preserved} by $\vee$ when diagrams of the following form commute
in \aA:

\begin{center}
\begin{picture}(300,105)
\put(3,80){$M_1(\vec{A}_m)\vee M_1(\vec{A}_m')$}

\put(180,80){$M_2(\vec{A}_{\pi(m)})\vee M_2(\vec{A}_{\pi(m)}')$}

\put(108,93){$\beta_{\vec{A}_m}\!\vee\beta_{\vec{A}_m'}$}

\put(92,85){\vector(1,0){80}}

%%%%%%%%%%%

\put(7,48){$\overline{\psi}^{M_1}_{\scriptsize\overrightarrow{A_m,A_m'}}$}

\put(237,48){$\overline{\psi}^{M_2}_{\scriptsize\overrightarrow{A_{\pi(m)},A_{\pi(m)}'}}$}

\put(46,70){\vector(0,-1){40}}

\put(233,70){\vector(0,-1){40}}

%%%%%%%%%%%%%

\put(8,10){$M_1(\overrightarrow{A_m\vee A_m'})$}

\put(184,10){$M_2(\overrightarrow{A_{\pi(m)}\vee A_{\pi(m)}'})$}

\put(110,6){$\beta_{\scriptsize\overrightarrow{A_m\vee A_m'}}$}

\put(92,15){\vector(1,0){80}}

\end{picture}
\end{center}
i.e., we have in \aA\ the equation
\begin{tabbing}
\mbox{\hspace{2.5em}}$(\overline{\psi}\beta)\quad\quad
\beta_{\scriptsize\overrightarrow{A_m\vee A_m'}}\cirk
\overline{\psi}^{M_1}_{\scriptsize\overrightarrow{A_m,A_m'}}=
\overline{\psi}^{M_2}_{\scriptsize\overrightarrow{A_{\pi(m)},A_{\pi(m)}'}}\cirk
(\beta_{\vec{A}_m}\vee\beta_{\vec{A}_m'})$.
\end{tabbing}

Let ${\langle\aA,\kon,\vee\rangle}$ be a biassociative category,
and let us make the assumption ${(\vee\kon)}$ (i.e., $\vee$
intermutes with $\kon$) together with the assumptions
\begin{tabbing}
\hspace{2.5em}\=${(\psi\check{b})}$\quad\quad\=$\check{b}^\str$
\=is upward preserved by $\kon$,\\[.5ex]
\>${(\overline{\psi}\hat{b})}$\>$\hat{b}^\rts$ \>is downward
preserved by $\vee$.
\end{tabbing}
(The biassociative intermuting categories, introduced in Section
10 below, satisfy these assumptions.) It is easy to see that
$\check{b}^\str$ is upward preserved by $\kon$ iff
$\check{b}^\rts$ is, and analogously with $\hat{b}$ and downward
preservation by $\vee$. We take ${}^\str$ in ${(\psi\check{b})}$
and ${}^\rts$ in ${(\overline{\psi}\hat{b})}$ to make duality
apparent.

The assumption ${(\psi\check{b})}$ amounts to the equation
\begin{tabbing}
\hspace{.5em}$(\psi\check{b})\quad c^k_{A_1\vee A_2,A'_1\vee
A'_2,A_3,A'_3}\cirk(c^k_{A_1,A_1',A_2,A_2'}\vee\mj_{A_3\kon
A_3'})\cirk \check{b}^\str_{A_1\kon A_1',A_2\kon A_2',A_3\kon
A_3'}=$\\*[1ex]
\`$(\check{b}^\str_{A_1,A_2,A_3}\kon\check{b}^\str_{A_1',A_2',A_3'})\cirk
c^k_{A_1,A_1',A_2\vee A_3,A_2'\vee A_3'}\cirk(\mj_{A_1\kon
A_1'}\vee c^k_{A_2,A_2',A_3,A_3'})$,
\end{tabbing}
and the assumption ${(\overline{\psi}\hat{b})}$ amounts to the
equation
\begin{tabbing}
\hspace{.5em}$(\overline{\psi}\hat{b})\quad \hat{b}^\rts_{A_1\vee
A_1',A_2\vee A_2',A_3\vee
A_3'}\cirk(c^k_{A_1,A_2,A_1',A_2'}\kon\mj_{A_3\vee A_3'})\cirk
c^k_{A_1\kon A_2,A_3,A_1'\kon A_2',A_3'}=$\\*[1ex]
\`$(\mj_{A_1\vee A_1'}\kon c^k_{A_2,A_3,A_2',A_3'})\cirk
c^k_{A_1,A_2\kon A_3,A_1',A_2'\kon
A_3'}\cirk(\hat{b}^\rts_{A_1,A_2,A_3}\!\vee\hat{b}^\rts_{A_1',A_2',A_3'})$.
\end{tabbing}
We call collectively these two equations ${(c^k b)}$. They stem
from \cite{BFSV} (Section~1), where, however, the $b$-arrows are
identity arrows.

Let ${\langle\aA,\kon,\vee,\top,\bot\rangle}$ be a biunital
category, and let us make the assumptions ${(\vee\kon)}$,
${(\bot\kon)}$ and ${(\vee\top)}$ together with the assumptions
\begin{tabbing}
\hspace{2.5em}\=${(\psi\check{\delta}\check{\sigma})}$\quad\quad\=$\check{\delta}^\str$
and $\check{\sigma}^\str$ \=are upward preserved by
$\kon$,\\*[.5ex]
\>${(\overline{\psi}\hat{\delta}\hat{\sigma})}$\>$\hat{\delta}^\rts$
and $\hat{\sigma}^\rts$ \>are downward preserved by $\vee$.
\end{tabbing}
The assumption ${(\psi\check{\delta}\check{\sigma})}$ amounts to
the equations
\begin{tabbing}
\hspace{2.5em}\=${(\psi\check{\delta}\check{\sigma})}$\quad\quad\=$\check{\delta}^\str$
and $\check{\sigma}^\str$ are upward preserved by $\kon$\kill

\>$(\psi\check{\delta})$\>$\check{\delta}^\str_{A\kon
A'}$\=$=(\check{\delta}^\str_A\kon\check{\delta}^\str_{A'})\,$\=$\cirk
c^k_{A,A',\bot,\bot}\cirk(\mj_{A\kon A'}\vee
\hat{w}^\rts_\bot)$,\\*[1ex]
\>$(\psi\check{\sigma})$\>$\check{\sigma}^\str_{A\kon
A'}$\>$=(\check{\sigma}^\str_A\kon\check{\sigma}^\str_{A'})$\>$\cirk
c^k_{\bot,\bot,A,A'}\cirk(\hat{w}^\rts_\bot\vee\mj_{A\kon A'})$,
\end{tabbing}
and the assumption ${(\overline{\psi}\hat{\delta}\hat{\sigma})}$
amounts to the equations
\begin{tabbing}
\hspace{2.5em}\=${(\psi\check{\delta}\check{\sigma})}$\quad\quad\=$\check{\delta}^\str$
and $\check{\sigma}^\str$ are upward preserved by $\kon$\kill

\>$(\overline{\psi}\hat{\delta})$\>$\hat{\delta}^\rts_{A\vee
A'}$\=$=(\mj_{A\vee A'}\kon\check{w}^\str_\top)$\=$\cirk
c^k_{A,\top,A',\top}\cirk(\hat{\delta}^\rts_A\vee\hat{\delta}^\rts_{A'})$,\\*[1ex]
\>$(\overline{\psi}\hat{\sigma})$\>$\hat{\sigma}^\rts_{A\vee
A'}$\>$=(\check{w}^\str_\top\kon\mj_{A\vee A'})$\>$\cirk
c^k_{\top,A,\top,A'}\cirk(\hat{\sigma}^\rts_A\vee\hat{\sigma}^\rts_{A'})$.
\end{tabbing}
We call collectively these four equations ${(c^k\delta\sigma)}$.

For the definitions below we make the assumptions ${(\bot\kon)}$
and ${(\vee\top)}$ together with the assumption
\begin{tabbing}
\hspace{2.5em}\=${(\psi\check{\delta}\check{\sigma})}$\quad\quad\=$\check{\delta}^\str$
and $\check{\sigma}^\str$ are upward preserved by $\kon$\kill

\>${(\bot\top)}$\>$\bot$ intermutes with $\top$,
\end{tabbing}
which delivers the arrow $\kappa\!:\bot\str\top$ of \aA. We define
by induction on the complexity of the ${(\vee,\bot)}$-shape $M$
the arrow ${\pM{M}\!:M(\top,\ldots,\top)\str\top}$ of \aA:
\[
\begin{array}{l}
\pM{\bot}=\kappa,\\[.5ex] \pM{\Box}=\mj_\top,\\[.5ex]
\pM{M\vee N}=\check{w}^\str_\top\cirk(\pM{M}\vee\pM{N}).
\end{array}
\]
Next we define by induction on the complexity of the
${(\kon,\top)}$-shape $M$ the arrow ${\ppM{M}\!:\bot\str
M(\bot,\ldots,\bot)}$ of \aA:
\[
\begin{array}{l}
\ppM{\top}=\kappa,\\[.5ex] \ppM{\Box}=\mj_\bot,\\[.5ex]
\ppM{M\kon N}=(\ppM{M}\kon\ppM{N})\cirk\hat{w}^\rts_\bot.
\end{array}
\]
When $\top$ and $\bot$ are conceived as nullary $\kon$ and $\vee$
respectively, the definitions of $\pM{M}$ and $\ppM{M}$ are
analogous to those of $\psi^M$ and $\overline{\psi}^M$
respectively.

Let now $\alpha$ and $\beta$ be natural transformations as for the
equations ${(\psi\alpha)}$ and ${(\overline{\psi}\beta)}$. We say
that $\alpha$ is \emph{upward preserved} by $\top$ when in \aA\ we
have the equation
\begin{tabbing}
\hspace{2.5em}\=${(\psi\check{\delta}\check{\sigma})}$\quad\quad\=
$\pM{M_2}\!\!\cirk\,\alpha_{\top,\ldots,\top}\,$\=\kill
\>$(\pM{}\!\alpha)$\>
$\pM{M_2}\!\!\cirk\,\alpha_{\top,\ldots,\top}\,$\>$=\mj_\top\cirk\pM{M_1}$.
\end{tabbing}
We say that $\beta$ is \emph{downward preserved} by $\bot$ when in
\aA\ we have the equation
\begin{tabbing}
\hspace{2.5em}\=${(\psi\check{\delta}\check{\sigma})}$\quad\quad\=
$\pM{M_2}\!\!\cirk\,\alpha_{\top,\ldots,\top}\,$\=\kill

\>$(\ppM{}\!\beta)$\>$
\beta_{\bot,\ldots,\bot}\cirk\ppM{M_1}$\>$=\;\ppM{M_2}\!\!\cirk\mj_\bot$.
\end{tabbing}
We left $\mj_\top$ and $\mj_\bot$ in these two equations to make
apparent the analogy with the equations ${(\psi\alpha)}$ and
${(\overline{\psi}\beta)}$.

Let ${\langle\aA,\kon,\vee,\top,\bot\rangle}$ be a bimonoidal
category, and let us make the assumptions ${(\bot\kon)}$,
${(\vee\top)}$ and ${(\bot\top)}$ together with the assumptions
\begin{tabbing}
\hspace{2.5em}\=${(\psi\check{\delta}\check{\sigma})}$\quad\quad\=
$\pM{M_2}\!\!\cirk\,\alpha_{\top,\ldots,\top}\,$\kill

\>${(\pM{}\!\check{b})}$\>$\check{b}^\str$ \=is upward preserved
by $\top$,\\*
\>${(\pM{}\!\check{\delta}\check{\sigma})}$\>$\check{\delta}^\str$
\=and $\check{\sigma}^\str$ are upward preserved by $\top$,\\[1ex]
\>${(\ppM{}\!\hat{b})}$\>$\hat{b}^\rts$ \>is downward preserved by
$\bot$,\\*
\>${(\ppM{}\!\hat{\delta}\hat{\sigma})}$\>$\hat{\delta}^\rts$
\>and $\hat{\sigma}^\rts$ are downward preserved by $\bot$.
\end{tabbing}

The assumptions ${(\pM{}\!\check{b})}$ and ${(\ppM{}\!\hat{b})}$
amount respectively to the equations
\begin{tabbing}
\hspace{2.5em}\=${(\psi\check{\delta}\check{\sigma})}$\quad\quad\=
$\pM{M_2}\!\!\cirk\,\alpha_{\top,\ldots,\top}\,$\kill

\>${(\pM{}\!\check{b})}$\>$\check{w}^\str_\top\cirk
(\check{w}^\str_\top\vee\mj_\top)\cirk\check{b}^\str_{\top,\top,\top}$\=$=
\check{w}^\str_\top\cirk(\mj_\top\vee\check{w}^\str_\top)$,\\*[.5ex]
\>${(\ppM{}\!\hat{b})}$\>$\hat{b}^\rts_{\bot,\bot,\bot}\!\cirk(\hat{w}^\rts_\bot\kon\mj_\bot)
\cirk\hat{w}^\rts_\bot$\>$=(\mj_\bot\kon\hat{w}^\rts_\bot)\cirk\hat{w}^\rts_\bot$,
\end{tabbing}
which are analogous to the equations ${(c^k b)}$. We call
collectively these two equations ${(wb)}$.

The assumptions ${(\pM{}\!\check{\delta}\check{\sigma})}$ and
${(\ppM{}\!\hat{\delta}\hat{\sigma})}$ amount respectively to the
equations
\begin{tabbing}
\hspace{2.5em}\=${(\psi\check{\delta}\check{\sigma})}$\quad\quad\=
$\pM{M_2}\!\!\cirk\,\alpha_{\top,\ldots,\top}\,$\kill

\>${(\pM{}\!\check{\delta})}$\>$\check{\delta}^\str_\top\,$\=$=\check{w}^\str_\top\cirk
(\mj_\top\vee\kappa)$,\quad\quad\quad\quad\=${(\ppM{}\!\hat{\delta})}$\quad\quad\=
$\hat{\delta}^\rts_\bot\,$\=$=
(\mj_\bot\kon\kappa)\cirk\hat{w}^\rts_\bot$,\\*[.5ex]
\>${(\pM{}\!\check{\sigma})}$\>$\check{\sigma}^\str_\top$\>$=\check{w}^\str_\top\cirk
(\kappa\vee\mj_\top)$,
\>${(\ppM{}\!\hat{\sigma})}$\>$\hat{\sigma}^\rts_\bot$\>$=
(\kappa\kon\mj_\bot)\cirk\hat{w}^\rts_\bot$,
\end{tabbing}
which are analogous to the equations ${(c^k\delta\sigma)}$. We
call collectively these four equations ${(\kappa\delta\sigma)}$.

To define the bimonoidal intermuting categories in Section 12 we
need also the following equations, which we call collectively
${(c^k\kappa)}$:

\begin{tabbing}
\hspace{2.5em}\=${(\psi\check{\delta}\check{\sigma})}$\quad\quad\=
$\pM{M_2}\!\!\cirk\,\alpha_{\top,\ldots,\top}\,$\kill

\>\>$\hat{\delta}^\str_\top\cirk(\check{\delta}^\str_\top\kon\check{\sigma}^\str_\top)\cirk
 c^k_{\top,\bot,\bot,\top}$\=$=\kappa
\cirk\check{\delta}^\str_\bot\cirk(\hat{\sigma}^\str_\bot\vee\hat{\delta}^\str_\bot)$,\\*[1ex]
\>\>$\hat{\delta}^\str_\top\cirk(\check{\sigma}^\str_\top\kon\check{\delta}^\str_\top)\cirk
c^k_{\bot,\top,\top,\bot}$\>$=\kappa
\cirk\check{\delta}^\str_\bot\cirk(\hat{\delta}^\str_\bot\vee\hat{\sigma}^\str_\bot)$.
\end{tabbing}
From these two equations we obtain either two alternative
definitions of $\kappa$ in terms of $c^k_{\top,\bot,\bot,\top}$ or
$c^k_{\bot,\top,\top,\bot}$, or the definitions of these two
arrows in terms of $\kappa$. These equations bear an analogy to
the equations ${(wb)}$ (see the equation ${(i\alpha)}$ in
Section~7). As in the equation ${(\pM{}\!\check{b})}$ the arrow
$\check{b}^\str_{\top,\top,\top}$ is ``shifted'' to $\mj_\top$ by
simplifying arrows made of $\check{w}^\str_\top$, so in the
equations ${(c^k\kappa)}$ the arrows $c^k_{\top,\bot,\bot,\top}$
and $c^k_{\bot,\top,\top,\bot}$ are shifted to $\kappa$ by
analogous simplifying arrows.

A justification of the assumptions for biassociative categories in
the spirit of this section and of the preceding one may be found
in \cite{DP06a}. This involves in particular a justification of
Mac Lane's pentagonal equation. In Sections 14 and 16 we will
mention further justifications in the same spirit of the
assumptions made for categories with symmetry, i.e.\ natural
commutativity isomorphisms.

\section{\large\bf Normal biunital categories} In this and in the
next two sections we present auxiliary coherence results involving
the unit objects $\top$ and $\bot$. These results will serve for
the coherence results of Sections 12 and 16.

A \emph{normal biunital} category is a biunital category
${\langle\aA,\kon,\vee,\top,\bot\rangle}$ (see Section~2) such
that in \aA\ we have the isomorphisms
${\hat{w}^\rts_\bot\!:\bot\str\bot\kon\bot}$ and
$\check{w}^\str_\top\!:\top\vee\top\str\top$. According to the
terminology of Section~3, in a normal biunital category for every
$\!\ks\!$ in ${\{\kon,\vee\}}$ and every $\zeta$ in
${\{\top,\bot\}}$ we have that $\!\ks\!$ and $\zeta$ intermute
with each other. The inverses of $\hat{w}^\rts_\bot$ and
$\check{w}^\str_\top$ are $\hat{w}^\str_\bot$ and
$\check{w}^\rts_\top$ respectively. We call these four arrows
collectively \emph{w-arrows}. The normal biunital category freely
generated by a set of objects is called $\mn_{\top,\bot}$.

For every formula $A$ we define the formula $\nu(A)$, which is the
\emph{normal form} of $A$, in the following way:
\begin{tabbing}
\hspace{2.5em}$\nu(p)=p$, \hspace{.5em}for $p$ a
letter,\hspace{3.5em}$\nu(\zeta)=\zeta$, \hspace{.5em}for
$\zeta\in\{\top,\bot\}$,
\end{tabbing}
for $\!\ks\!\in\{\kon,\vee\}$ and $i\in\{1,2\}$,
\[
\nu(A_1\ks A_2)= \left\{
\begin{array}{l}
\zeta \hspace{2.9em} {\mbox{\rm if }} \nu(A_1)=\nu(A_2)=\zeta,
\\[.5ex]
\nu(A_i) \hspace{1em} {\mbox{\rm if }}
\nu(A_{3-i})=\zeta\neq\nu(A_i)
{\mbox{\rm{ and }}} (\!\ks\!,\zeta)\in\{(\kon,\top),(\vee,\bot)\},\\[.5ex]
\nu(A_1)\ks\nu(A_2) \quad {\mbox{\rm otherwise.}}
\end{array}
\right .
\]
If no letter occurs in $A$, then $\nu(A)$ is either $\top$ or
$\bot$. It is clear that for every formula $A$ we have an
isomorphism of $\mn_{\top,\bot}$ of the type ${A\str\nu(A)}$.

An arrow term of $\mn_{\top,\bot}$ is called \emph{directed} when
${}^\rts$ does not occur in it as a superscript. We can prove the
following.

\prop{Directedness Lemma}{If $f,g\!:A\str\nu(B)$ are directed
arrow terms of $\mn_{\top,\bot}$, then $f=g$ in
$\mn_{\top,\bot}$.}

\noindent The proof of this lemma is analogous to Mac Lane's proof
of a lemma that delivered monoidal coherence (see \cite{ML63},
\cite{ML71}, Section VII.2, or \cite{DP04}, Sections 4.3 and 4.6).
Whenever for ${i\in\{1,2\}}$ we have that the arrow terms
${f_i\!:A\str A_i}$ are factors of two directed developed arrow
terms, we establish that we have two arrow terms ${f_i'\!:A_i\str
C}$ such that ${f_1'\cirk f_1=f_2'\cirk f_2}$, and we have a
directed arrow term of the type ${C\str\nu(B)}$, because
${\nu(B)}$ is in normal form. For that we use bifunctorial and
naturality equations, except for the case where the heads of $f_1$
and $f_2$ are $\d{\xi}{\str}_\zeta$ and $\s{\xi}{\str}_\zeta$ for
${(\!\ks\!,\zeta)\in\{(\kon,\top),(\vee,\bot)\}}$. Then we use the
equation $\d{\xi}{\str}_\zeta\,=\;\s{\xi}{\str}_\zeta$ of biunital
categories (which we have also in monoidal categories). We can
then prove the following.

\prop{Normal Biunital Coherence}{The category $\mn_{\top,\bot}$ is
a preorder.}

\noindent This is obtained from the Directedness Lemma, as Mac
Lane obtained monoidal coherence (see the proof of Associative
Coherence in \cite{DP04}, Section 4.3).

Normal Biunital Coherence guarantees that there is a unique arrow
of $\mn_{\top,\bot}$ of the type ${A\str\nu(A)}$, and of the
converse type. These arrows are isomorphisms.

\section{\large\bf $\kappa$-Normal biunital categories}
A \emph{$\kappa$-normal biunital} category is a normal biunital
category ${\langle\aA,\kon,\vee,\top,\bot\rangle}$ (see the
preceding section) such that in \aA\ we have an arrow
${\kappa\!:\bot\str\top}$ that satisfies the equations
${(\kappa\delta\sigma)}$ (see Section~4). We call
$\mk^\emptyset_{\top,\bot}$ the $\kappa$-normal biunital category
freely generated by the empty set of objects. So there are no
objects of $\mk^\emptyset_{\top,\bot}$ in which letters occur. We
prove the following.

\prop{$\mk^\emptyset_{\top,\bot}$ Coherence}{The category
$\mk^\emptyset_{\top,\bot}$ is a preorder.}

\dkz For every arrow term $f$ of $\mk^\emptyset_{\top,\bot}$
either there is an arrow term $f'$ of $\mn_{\top,\bot}$ such that
${f=f'}$ in $\mk^\emptyset_{\top,\bot}$, or there are two arrow
terms $f'$ and $f''$ of $\mn_{\top,\bot}$ such that
${f=f''\cirk\kappa\cirk f'}$ in $\mk^\emptyset_{\top,\bot}$. This
is established by using the equations ${(\kappa\delta\sigma)}$ and
the following consequences of the naturality equations for the
$\delta$-$\sigma$-arrows:
\begin{tabbing}
\hspace{5em}\=$\kappa\kon\mj_\top$\=
$=\hat{\delta}^\rts_\top\cirk\kappa\cirk\hat{\delta}^\str_\bot$,
\hspace{5em}\=$\mj_\top\kon\kappa\,$\=
$=\hat{\sigma}^\rts_\top\cirk\kappa\cirk\hat{\sigma}^\str_\bot$,\\*[1ex]
\>$\kappa\vee\mj_\bot$\>
$=\check{\delta}^\rts_\top\cirk\kappa\cirk\check{\delta}^\str_\bot$,
\>$\mj_\bot\vee\kappa\;$\>
$=\check{\sigma}^\rts_\top\cirk\kappa\cirk\check{\sigma}^\str_\bot$.
\end{tabbing}
Note that there are no arrow terms of $\mk^\emptyset_{\top,\bot}$
of the form ${\kappa\cirk g\cirk\kappa}$.

Suppose now for ${i\in\{1,2\}}$ that ${f_i\!:A\str B}$ is an arrow
term of $\mk^\emptyset_{\top,\bot}$. We have either ${f_1=f_1'}$
and ${f_2=f_2'}$, or ${f_1=f_1''\cirk\kappa\cirk f_1'}$ and
${f_2=f_2''\cirk\kappa\cirk f_2'}$, for $f_i'$ and $f_i''$ arrow
terms of $\mn_{\top,\bot}$. It is impossible that ${f_i=f_i'}$ and
${f_{3-i}=f_{3-i}''\cirk\kappa\cirk f_{3-i}'}$ because
${f_i=f_i'}$ requires that $A$ be isomorphic to $B$, while
${f_{3-i}=f_{3-i}''\cirk\kappa\cirk f_{3-i}'}$ prevents that. Then
we just apply Normal Biunital Coherence. \qed

We suppose that if to $\kappa$-normal biunital categories we add
the equations
\begin{tabbing}
\hspace{9em}\=$\hat{\delta}^\str_{A\kon\bot}\cirk(\mj_{A\kon\bot}\kon\kappa)\,$\=
$=(\hat{\delta}^\str_A\cirk(\mj_A\kon\kappa))\kon\mj_\bot$,\\*[1ex]
\>$(\mj_{A\vee\top}\vee\kappa)\cirk\check{\delta}^\rts_{A\vee\top}$\>
$=((\mj_A\vee\kappa)\cirk\check{\delta}^\rts_A)\vee\mj_\top$,\\[2ex]
\>$\hat{\sigma}^\str_{\bot\kon A}\cirk(\kappa\kon\mj_{\bot\kon
A})$\>
$=\mj_\bot\kon(\hat{\sigma}^\str_A\cirk(\kappa\kon\mj_A))$,\\*[1ex]
\>$(\kappa\vee\mj_{\top\vee A})\cirk\check{\sigma}^\rts_{\top\vee
A}$\>
$=\mj_\top\vee((\kappa\vee\mj_A)\cirk\check{\sigma}^\rts_A)$,
\end{tabbing}
then we could prove coherence in the sense of preordering for the
resulting categories without assuming that the set of generators
for the freely generated category is empty. (For categories where
these equations hold see the end of the next section.)

\section{\large\bf Normal bimonoidal categories}
A \emph{normal bimonoidal} category is a bimonoidal category (see
Section~2) that is also normal biunital (see Section~5), in which,
moreover, the equations ${(wb)}$ (see Section~4) are satisfied. We
call $\mna_{\top,\bot}$ the normal bimonoidal category freely
generated by a set of objects. We can prove the following.

\prop{Normal Bimonoidal Coherence}{The category $\mna_{\top,\bot}$
is a preorder.}

\dkz By Biassociative Coherence (see Section~2) and the results of
\cite{DP04} (Chapter~3), we can replace the category
$\mna_{\top,\bot}$ by an equivalent strictified category
$\mna_{\top,\bot}^{st}$, where the $b$-arrows are identity arrows.
Then to show that $\mna_{\top,\bot}^{st}$ is a preorder, which
implies that $\mna_{\top,\bot}$ is a preorder, we proceed as for
the proof of Normal Biunital Coherence in Section~5. Here we apply
the equations ${(wb)}$. \qed

A \emph{$\kappa$-normal bimonoidal} category is a normal
bimonoidal category that is also $\kappa$-normal biunital (see the
preceding section). We call $\mka^\emptyset_{\top,\bot}$ the
$\kappa$-normal bimonoidal category freely generated by the empty
set of objects. We can prove the following.

\prop{$\mka^\emptyset_{\top,\bot}$ Coherence}{The category
$\mka^\emptyset_{\top,\bot}$ is a preorder.}

\dkz We enlarge the proof of $\mk^\emptyset_{\top,\bot}$ Coherence
of the preceding section by using the following equations of
monoidal categories:
\begin{tabbing}
\hspace{11em}\=$\b{\xi}{\str}_{A,\zeta,C}\,$\=$=\;\d{\xi}{\rts}_A\ks\s{\xi}{\str}_C$,\\[1ex]
\>$\b{\xi}{\str}_{A,B,\zeta}$\>$=\;\d{\xi}{\rts}_{A\kst
B}\!\cirk(\mj_A\ks\d{\xi}{\str}_B)$,\\[1ex]
\>$\b{\xi}{\str}_{\zeta,B,C}$\>$=(\s{\xi}{\rts}_B\ks\mj_C)\cirk\s{\xi}{\str}_{B\kst
C}$,
\end{tabbing}
for ${(\!\ks\!,\zeta)\in\{(\kon,\top),(\vee,\bot)\}}$. We use also
analogous equations derived by isomorphism, where the superscript
${}^\str$ of $b$ is replaced by ${}^\rts$. With these equations
and the equations ${(wb)}$ we can eliminate every occurrence of
$b$. To achieve that we use also equations derived from naturality
equations, like the following equation:
\[
\b{\xi}{\str}_{A,B,C}\,=((\mj_A\ks
i_B^{-1})\ks\mj_C)\cirk\b{\xi}{\str}_{A,\nu(B),C}\!\!\cirk(\mj_A\ks(i_B\ks\mj_C))
\]
where ${i_B\!:B\str\nu(B)}$ is an arrow term of $\mn_{\top,\bot}$
standing for an isomorphism, and $i^{-1}_B$ stands for its
inverse.\qed

The equations ${(wb)}$ can be justified in a manner somewhat
different from that of Section~4 by appealing to the isomorphisms
${i_B\!:B\str\nu(B)}$ mentioned at the end of the proof above. For
${\alpha_{\zeta,\ldots,\zeta}\!:M(\zeta,\ldots,\zeta)\str
N(\zeta,\ldots,\zeta)}$ a component of a natural transformation,
we can take that the equations ${(wb)}$ are obtained from the
following equation:
\[
(i\alpha)\quad\quad
i_{N(\zeta,\ldots,\zeta)}\cirk\alpha_{\zeta,\ldots,\zeta}=\mj_\zeta\cirk
i_{M(\zeta,\ldots,\zeta)}.
\]
We read this equation by saying that
${\alpha_{\zeta,\ldots,\zeta}}$ is \emph{shifted} by two
isomorphisms to $\mj_\zeta$ (see the end of Section~4 for shifting
to $\kappa$).

In $\kappa$-normal bimonoidal categories we have the equations
mentioned after the proof of $\mk^\emptyset_{\top,\bot}$ Coherence
at the end of the preceding section, and we can prove coherence in
the sense of preordering for these categories without assuming
that the set of generators for the freely generated
$\kappa$-normal bimonoidal category is empty. We will however not
dwell on this proof, for which we have no application in the rest
of the paper.

The categories treated in this and in the preceding two sections
make the following chart:
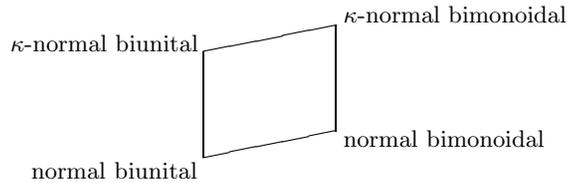
\begin{center}
\begin{picture}(200,90)

\put(70,30){\line(0,1){40}}

\put(120,40){\line(0,1){40}}

\put(70,30){\line(5,1){50}}

\put(70,70){\line(5,1){50}}

\put(5,22){\small normal biunital}

\put(-3,70){\small $\kappa$-normal biunital}

\put(123,34){\small normal bimonoidal}

\put(123,81){\small $\kappa$-normal bimonoidal}
\end{picture}
\end{center}

\section{\large\bf Rectangular notation}
To prove coherence results for categories with intermutation we
will find very helpful a planar notation for propositional
formulae involving $\kon$ and $\vee$, i.e.\ for the objects of our
freely generated categories. We first introduce a representation
of formula trees obtained by subdividing rectangles. To every
formula there will correspond a rectangle subdivided into further
rectangles corresponding to the subformulae; if the formula is a
letter, there is just one rectangle.

Out of this representation we obtain our planar notation, which we
call the \emph{rectangular notation} for formulae. We deal first
with binary trees involving two binary connectives, and next with
trees with finite branching, possibly bigger than binary,
involving two strictly associative binary connectives. Something
similar to our rectangular notation was used to explain the
Eckmann-Hilton argument in \cite{BD95} (Section VI) and
\cite{BN96} (Section 1.2). We do not know whether our rectangular
notation, and the use we make of it, may be connected
significantly with the little $n$-cubes studied in \cite{BFSV}
(Section~6).

To obtain our rectangular notation we go through the following
construction. At the beginning, we index every point in the plane
by the empty set. At the end of the construction the points that
are vertices of rectangles will receive as indices one of the
following sets:
\[
\emptyset,\{\downarrow\},\{\uparrow\},\{\downarrow,\uparrow\},\{\rightarrow\},\{\leftarrow\},
\{\rightarrow,\leftarrow\},
\]
while all the other points will be indexed by $\emptyset$.

The construction then proceeds as follows. For a formula $A$ we
draw a rectangle enclosing $A$. If during the construction we have
a rectangle enclosing a subformula $B\kon C$ of $A$:

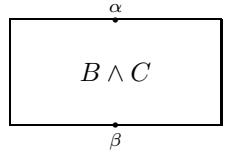
\begin{center}
\begin{picture}(80,60)

\put(0,10){\line(1,0){80}} \put(0,50){\line(1,0){80}}
\put(0,10){\line(0,1){40}} \put(80,10){\line(0,1){40}}

\put(40,10){\circle*{1.5}} \put(40,50){\circle*{1.5}}

\put(40,7){\makebox(0,0)[t]{\scriptsize$\beta$}}
\put(40,53){\makebox(0,0)[b]{\scriptsize$\alpha$}}
\put(40,30){\makebox(0,0){$B\kon C$}}

\end{picture}
\end{center}
then we delete $B\kon C$ and subdivide the rectangle in the middle
by a vertical line segment so as to obtain two rectangles
enclosing $B$ and $C$ respectively:

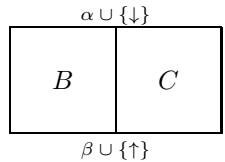
\begin{center}
\begin{picture}(80,60)

\put(0,10){\line(1,0){80}} \put(0,50){\line(1,0){80}}
\put(0,10){\line(0,1){40}} \put(80,10){\line(0,1){40}}
\put(40,10){\line(0,1){40}}

\put(40,7){\makebox(0,0)[t]{\scriptsize$\beta\cup\{\uparrow\}$}}
\put(40,53){\makebox(0,0)[b]{\scriptsize$\alpha\cup\{\downarrow\}$}}
\put(20,30){\makebox(0,0){$B$}} \put(60,30){\makebox(0,0){$C$}}

\end{picture}
\end{center}
Here $\alpha$ and $\beta$ are the index sets of the middle points
of the horizontal sides in the former figure, and these index sets
are changed after the subdivision as in the latter figure. All the
other points have index sets unchanged in passing from the former
to the latter figure. The vertical line segment introduced by the
subdivision \emph{corresponds} to the main $\kon$ of $B\kon C$,
and we say that the big rectangle that has previously enclosed
$B\kon C$ \emph{corresponds} to $B\kon C$.

If during the construction we have a rectangle enclosing a
subformula $B\vee C$ of $A$ as in the figure on the left, then by
subdividing it horizontally in the middle we pass to the figure on
the right:

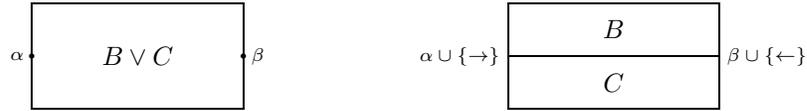
\begin{center}
\begin{picture}(260,60)

\put(0,10){\line(1,0){80}} \put(0,50){\line(1,0){80}}
\put(0,10){\line(0,1){40}} \put(80,10){\line(0,1){40}}

\put(180,10){\line(1,0){80}} \put(180,50){\line(1,0){80}}
\put(180,10){\line(0,1){40}} \put(260,10){\line(0,1){40}}
\put(180,30){\line(1,0){80}}

\put(0,30){\circle*{1.5}} \put(80,30){\circle*{1.5}}

\put(-3,30){\makebox(0,0)[r]{\scriptsize$\alpha$}}
\put(83,30){\makebox(0,0)[l]{\scriptsize$\beta$}}
\put(177,30){\makebox(0,0)[r]{\scriptsize$\alpha\cup\{\str\}$}}
\put(263,30){\makebox(0,0)[l]{\scriptsize$\beta\cup\{\rts\}$}}

\put(40,30){\makebox(0,0){$B\vee C$}}
\put(220,40){\makebox(0,0){$B$}} \put(220,20){\makebox(0,0){$C$}}

\end{picture}
\end{center}
with the new horizontal line segment corresponding to the main
$\vee$ of $B\vee C$, and the big rectangle on the right
corresponding to $B\vee C$. The construction is over when
rectangles enclose only letters, and the end result of the
construction is the \emph{binary rectangular grid} of $A$, which
we denote by $\rho(A)$. It is easy to see that $\rho$ is a one-one
map.

For example, for $A$ being $((p\kon q)\kon r)\vee(((s\kon
t)\vee(u\kon q))\kon((v\kon p)\vee w))$ we obtain the following
binary rectangular grid $\rho(A)$:

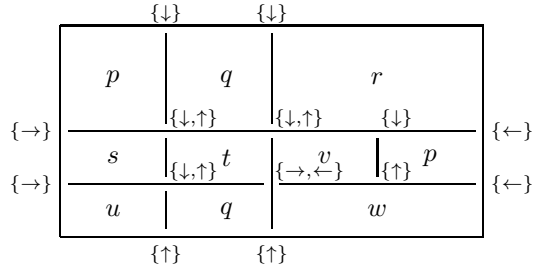
\begin{center}
\begin{picture}(160,100)

\put(0,10){\line(1,0){160}} \put(0,90){\line(1,0){160}}
\put(0,10){\line(0,1){80}} \put(160,10){\line(0,1){80}}
\put(3,50){\line(1,0){154}} \put(3,30){\line(1,0){74}}
\put(83,30){\line(1,0){74}} \put(40,13){\line(0,1){14}}
\put(40,33){\line(0,1){14}} \put(40,53){\line(0,1){34}}
\put(80,13){\line(0,1){34}} \put(80,53){\line(0,1){34}}
\put(120,33){\line(0,1){14}}

\put(-3,30){\makebox(0,0)[r]{\scriptsize$\{\str\}$}}
\put(-3,50){\makebox(0,0)[r]{\scriptsize$\{\str\}$}}
\put(163,30){\makebox(0,0)[l]{\scriptsize$\{\rts\}$}}
\put(163,50){\makebox(0,0)[l]{\scriptsize$\{\rts\}$}}
\put(40,7){\makebox(0,0)[t]{\scriptsize$\{\uparrow\}$}}
\put(80,7){\makebox(0,0)[t]{\scriptsize$\{\uparrow\}$}}
\put(121,33){\makebox(0,0)[lb]{\scriptsize$\{\uparrow\}$}}
\put(40,93){\makebox(0,0)[b]{\scriptsize$\{\downarrow\}$}}
\put(80,93){\makebox(0,0)[b]{\scriptsize$\{\downarrow\}$}}
\put(121,53){\makebox(0,0)[lb]{\scriptsize$\{\downarrow\}$}}
\put(41,33){\makebox(0,0)[lb]{\scriptsize$\{\downarrow,\!\uparrow\}$}}
\put(41,53){\makebox(0,0)[lb]{\scriptsize$\{\downarrow,\!\uparrow\}$}}
\put(81,53){\makebox(0,0)[lb]{\scriptsize$\{\downarrow,\!\uparrow\}$}}
\put(81,33){\makebox(0,0)[lb]{\scriptsize$\{\str,\!\rts\}$}}

\put(20,20){\makebox(0,0){$u$}} \put(63,20){\makebox(0,0){$q$}}
\put(120,20){\makebox(0,0){$w$}} \put(20,40){\makebox(0,0){$s$}}
\put(63,40){\makebox(0,0){$t$}} \put(100,40){\makebox(0,0){$v$}}
\put(140,40){\makebox(0,0){$p$}} \put(20,70){\makebox(0,0){$p$}}
\put(63,70){\makebox(0,0){$q$}} \put(120,70){\makebox(0,0){$r$}}

\end{picture}
\end{center}
All the points in the plane except those whose index set is
mentioned in this picture are indexed by $\emptyset$. By not
making the sides of the rectangles meet at crossing points we
indicate how these points are indexed, and we may omit mentioning
in the picture even the index sets we have mentioned.

If in the construction we have presented above we do not require
that the subdividing line segments, vertical and horizontal, be in
the middle, though they may be there, then we obtain a
\emph{rectangular grid} $\gamma(A)$ where $\gamma$ differs from
$\rho$ in not being one-one any more. To $\gamma(A)$ there will
correspond a set of formulae obtained from $A$ by associating
parentheses in different manners; more precisely, we will obtain
the set of all objects isomorphic to $A$ in the free biassociative
category \ma\ of Section~2.

\section{\large\bf Intermuting categories}

We call $\langle {\cal A},\kon,\vee\rangle$ an \emph{intermuting}
category when $\kon$ and $\vee$ are biendofunctors of \aA\ and
there is a natural transformation $c^k$ whose components are the
following arrows of $\cal A$:
\[
c^k_{A,B,C,D}\!:(A\kon B)\vee(C\kon D)\str(A\vee C)\kon(B\vee D).
\]
In the terminology of Section~3, we assume that $\vee$ intermutes
with $\kon$. Let \Ck\ be the free intermuting category generated
by a set of objects.

For every arrow term $f\!:A\str B$ of \Ck\ the binary rectangular
grids $\rho(A)$ and $\rho(B)$ differ only with respect to the
index sets; these grids are otherwise the same. This is because
$c^k_{A,B,C,D}$ corresponds to passing from the figure on the left
to the figure on the right:

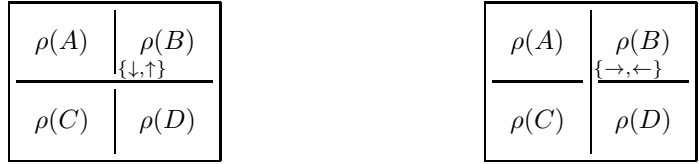
\begin{center}
\begin{picture}(260,80)

\put(0,10){\line(1,0){80}} \put(0,70){\line(1,0){80}}
\put(3,40){\line(1,0){74}} \put(0,10){\line(0,1){60}}
\put(80,10){\line(0,1){60}} \put(40,13){\line(0,1){24}}
\put(40,43){\line(0,1){24}}

\put(180,10){\line(1,0){80}} \put(180,70){\line(1,0){80}}
\put(180,10){\line(0,1){60}} \put(260,10){\line(0,1){60}}
\put(220,13){\line(0,1){54}} \put(183,40){\line(1,0){34}}
\put(223,40){\line(1,0){34}}

\put(41,43){\makebox(0,0)[lb]{\scriptsize$\{\downarrow,\!\uparrow\}$}}
\put(221,43){\makebox(0,0)[lb]{\scriptsize$\{\str,\!\rts\}$}}

\put(20,25){\makebox(0,0){$\rho(C)$}}
\put(60,25){\makebox(0,0){$\rho(D)$}}
\put(20,55){\makebox(0,0){$\rho(A)$}}
\put(60,55){\makebox(0,0){$\rho(B)$}}

\put(200,25){\makebox(0,0){$\rho(C)$}}
\put(240,25){\makebox(0,0){$\rho(D)$}}
\put(200,55){\makebox(0,0){$\rho(A)$}}
\put(240,55){\makebox(0,0){$\rho(B)$}}

\end{picture}
\end{center}
with all the index sets being unchanged except that indexing the
crossing in the middle of the figure. This crossing corresponds in
a unique way to $c^k_{A,B,C,D}$. So to every occurrence of $c^k$
in $f$ there corresponds a unique crossing in $\rho(A)$ and
$\rho(B)$, which is only indexed differently in $\rho(A)$ and
$\rho(B)$. Hence, for two arrow terms $f$ and $g$ of the same
type, there is a bijection from the set of occurrences of $c^k$ in
$f$ to the set of occurrences of $c^k$ in $g$ such that this
bijection maps an occurrence of $c^k$ in $f$ to the occurrence of
$c^k$ in $g$ with the same corresponding crossing.

We say that a $c^k$-term \emph{corresponds} to a crossing when its
head corresponds to this crossing. (The notions of $c^k$-term and
of its head are defined in Section~2.)

Note that not every crossing indexed with
$\{\downarrow,\uparrow\}$ is such that it can be reindexed with
$\{\rightarrow,\leftarrow\}$ through $c^k$, i.e.\ by intermuting.
In the example in the preceding section, the crossing involving
$p$, $q$, $s$ and  $t$ is not ready for intermuting. It will
become such after intermuting is performed at the two other
crossings indexed with $\{\downarrow,\uparrow\}$. But a crossing
indexed with $\{\downarrow,\uparrow\}$ may be such that it can
never become ready for intermuting, because it is not in the
centre of a rectangle.

Although in a composition $f_2\cirk f_1\!:A\str B$ of two
$c^k$-terms the two intermutings corresponding to $f_1$ and $f_2$
cannot always be ``permuted'', they can be ``permuted'' if in
$\rho(A)$ the crossings corresponding to $f_1$ and $f_2$ are ready
for intermuting. We can first easily establish the following
lemma.

\prop{Lemma~1}{If $f\!:A\str B$ and $g\!:A\str C$ are two
$c^k$-terms such that $B$ differs from $C$, then there are two
$c^k$-terms ${g'\!:B\str D}$ and ${f'\!:C\str D}$ such that
${g'\cirk f=f'\cirk g}$, and for ${h\in\{f,g\}}$ the $c^k$-terms
$h$ and $h'$ correspond to the same crossing in $\rho(A)$.}

\noindent Next we have the following lemma about ``permutation''.

\prop{Lemma~2}{Let ${f_n\cirk\ldots\cirk f_1\!:A\str B}$, for
${n\geq 2}$ be a composition of $c^k$-terms such that for
${i\in\{1,\ldots,n\}}$ the factor $f_i$ corresponds to the
crossing $x_i$ in $\rho(A)$, and let ${g_1\!:A\str C}$ be a
$c^k$-term that corresponds to $x_j$ for some
${j\in\{1,\ldots,n\}}$. Then there is a composition of $c^k$-terms
${g_n\cirk\ldots\cirk g_2\!:C\str B}$ such that
$f_n\cirk\ldots\cirk f_1=g_n\cirk\ldots\cirk g_1$ in \Ck.}

\dkz Let the type of $f_1$ be ${A\str A_1}$. If $A_1$ is $C$, then
$f_1=g_1$ and ${g_n\cirk\ldots\cirk g_2}$ is ${f_n\cirk\ldots\cirk
f_2}$. If $A_1$ is not $C$, then by Lemma~1, there are two
$c^k$-terms ${g_1'\!:A_1\str D}$ and ${f_1'\!:C\str D}$ such that
${g_1'\cirk f_1=f_1'\cirk g_1}$.

If $D$ is $B$, then $n=2$, and $g_2$ is $f_1'$. Since $f_2$ and
$g_1'$ must be the same $c^k$-term, we have ${f_2\cirk
f_1=g_2\cirk g_1}$.

If $D$ is not $B$, then $n>2$, and, because the crossing
corresponding to $g_1'$ is $x_j$ for $j\neq 1$, we may apply the
induction hypothesis to ${f_n\cirk\ldots\cirk f_2\!:A_1\str B}$
and $g_1'\!:A_1\str D$ to obtain $g_n\cirk\ldots\cirk g_3\!:D\str
B$ such that $f_n\cirk\ldots\cirk f_2=g_n\cirk\ldots\cirk g_3\cirk
g_1'$. We take that $g_2$ is $f_1'$, and we obtain

\begin{tabbing}
\mbox{\hspace{8em}}$f_n\cirk\ldots\cirk f_2\cirk f_1$ \= = \=
$g_n\cirk\ldots\cirk g_3\cirk g_1'\cirk f_1$
\\[.5ex]
\> = \> $g_n\cirk\ldots\cirk g_3\cirk g_2\cirk g_1$\`$\dashv$
\end{tabbing}

By applying Lemma~2 we can easily prove the following proposition
by induction on the length of a developed arrow term (see
Section~2).

\prop{Intermuting Coherence}{The category \Ck\ is a preorder.}

\section{\large\bf Biassociative intermuting categories}
A \emph{biassociative intermuting} category is a biassociative
category $\langle{\cal A},\kon,\vee\rangle$ (see Section~2) that
is an intermuting category (see the preceding section), and,
moreover, the equations $(c^kb)$ (see Section~4) are satisfied.
Let \ACk\ be the free biassociative intermuting category generated
by a set of objects.

One can show that every natural transformation defined by an arrow
term of \ma\ such that $\kon$ does not occur in it is upward
preserved by $\kon$ in \ACk, and analogously when $\kon$ is
replaced by $\vee$ and ``upward'' by ``downward''. For that we
rely essentially on the following. Suppose that $\beta'$ and
$\beta''$ are upward preserved by $\kon$. Then to show that
$\beta'\vee\beta''$ is also upward preserved by $\kon$ we rely on
bifunctorial equations and the naturality of $c^k$. To show that
the natural transformation $\beta'$ obtained from $\beta$ by
substituting a $\vee$-shape (see Section~4) in one of the indices
of $\beta$ is upward preserved by $\kon$ if $\beta$ is, we rely on
the naturality of $\beta$. In showing that, we have equations like
the following instance of the naturality of $\check{b}^\str$:

\begin{tabbing}
\mbox{\hspace{.5em}}$((\mj_{A_1\kon A_1'}\kon\mj_{A_2\kon
A_2'})\kon c^k_{A_3,A_3',A_4,A_4'})\cirk\check{b}^\str_{A_1\kon
A_1',A_2\kon A_2',(A_3\kon A_3')\vee(A_4\kon A_4')}=$
\\[.5ex]
\`$\check{b}^\str_{A_1\kon A_1',A_2\kon A_2',(A_3\vee
A_4)\kon(A_3'\vee A_4')}\cirk(\mj_{A_1\kon A_1'}\kon(\mj_{A_2\kon
A_2'}\kon c^k_{A_3,A_3',A_4,A_4'}))$.
\end{tabbing}
Something analogous can be shown for categories more complex than
\ACk, which we will encounter later in this paper (cf.\ Section
14).

As we did in the proof of normal Bimonoidal Coherence in
Section~7, we can apply Biassociative Coherence to obtain a
strictified category $\ACk^{st}$ equivalent to \ACk\ where the
$b$-arrows are identity arrows. Our ultimate goal is to show that
$\ACk^{st}$ is a preorder, which implies that \ACk\ is a preorder
too.

Before working towards that goal, we will consider in this section
the problem whether there is an arrow of a given type in
$\ACk^{st}$. This sort of problem (which is called the
\emph{theoremhood} problem in \cite{DP04}, Section 1.1) is
sometimes taken to be a part of a coherence result (cf.\
\cite{KML71}, Theorem 2.1, and \cite{BFSV}, Theorem 3.6.2). We do
not need to solve this problem to show that $\ACk^{st}$ is a
preorder, which is properly coherence, but the techniques used in
this section will be analogous to those used in the next section
to demonstrate coherence.

The objects of $\ACk^{st}$ are equivalence classes of formulae
$\lz A\dz$ such that $\lz A\dz$ is the set of all formulae
isomorphic to $A$ in the free biassociative category \ma\ of
Section~2. Such an equivalence class $\lz A\dz$ corresponds in a
one-to-one manner to the rectangular grid $\gamma(A)$ mentioned at
the end of Section~8. We may denote the equivalence class $\lz
A\dz$ by deleting from the formula $A$ parentheses tied to an
occurrence of $\!\ks\!$ within the immediate scope of another
occurrence of $\!\ks\!$, for $\!\ks\!\in\{\kon,\vee\}$. We call
the result of this deleting procedure a \emph{form sequence}
(which is short for \emph{form sequence in natural notation},
according to the terminology of \cite{DP04}, Section 6.2).

We call a form sequence \emph{diversified} when every letter
occurs in it at most once. To simplify matters, we speak from now
on only about diversified form sequences. It is easier if we speak
about letters, rather than their occurrences, and with diversified
form sequences we may do so. We denote by ${let(X)}$ the set of
letters occurring in the form sequence $X$.

If $B$ is in the equivalence class $\lz A\dz$, then the
rectangular grids $\gamma(A)$ and $\gamma(B)$ may be taken to be
the same, and instead of $\gamma(A)$ we write $\gamma(X)$, where
$X$ is the form sequence obtained from $A$ by the deleting
procedure above. Formally, we take $\gamma(X)$ to be an
equivalence class.

For a form sequence $X$ we define inductively four sequences of
letters taken from $X$ which we call $T(X)$, $B(X)$, $L(X)$ and
$R(X)$ ($T$ stands for \emph{top}, $B$ for \emph{bottom}, $L$ for
\emph{left} and $R$ for \emph{right}). The sequences $T(p)$,
$B(p)$, $L(p)$ and $R(p)$ are all the one-member sequence $p$. For
\[(S,\!\ks\!)\in\{(T,\kon),(B,\kon),(L,\vee),(R,\vee)\}
\]
we have that $S(X_1\ks X_2)$ is the sequence obtained by
concatenating the sequences $S(X_1)$ and $S(X_2)$. In the
remaining cases we have
\[
\begin{array}{ll}
T(X_1\vee X_2)=T(X_1), \quad & B(X_1\vee X_2)=B(X_2),
\\[.5ex]
L(X_1\kon X_2)=L(X_1), \quad & R(X_1\kon X_2)=R(X_2).
\end{array}
\]

If, for example, $X$ is the form sequence
\[(p_1\kon q_1\kon r)\vee(((s\kon t)\vee(u\kon q_2))\kon((v\kon
p_2)\vee w)),
\]
which is obtained from the formula in the example in Section~8,
then we have the following:
\[
\begin{array}{ll}
T(X)=p_1q_1r, \quad & B(X)=uq_2w,
\\[.5ex]
L(X)=p_1su, \quad & R(X)=rp_2w.
\end{array}
\]
Note that these sequences can easily be read from $\gamma(X)$,
which is obtained from the grid in Section~8 just by adding
subscripts to $p$ and $q$.

For $x$ an occurrence of $\kon$ in a form sequence $X$ such that
$X_1xX_2$ is a \emph{subformsequence} of $X$, i.e.\ a form
sequence that is subword of $X$, we define the sequences $L_x$ and
$R_x$ as $L(X_2)$ and $R(X_1)$ respectively. Intuitively, we may
read $L_x$ as ``the vertical sequence of letters that have $x$
immediately on the \emph{left} in the rectangular grid
$\gamma(X)$'', and analogously for $R_x$ (and also for $T_y$ and
$B_y$ below). Note that the same $x$ may occur in two different
subformsequences $X_1 x X_2$ and $X_1' x X_2'$ of $X$. For
example, $p\kon q$ and $p\kon q\kon r$ are both subformsequences
of $p\kon q\kon r$. It is however easy to check that the
definition above is correct, since if $X_1' x X_2'$ is also a
subformsequence of $X$, then $L(X_2)=L(X_2')$ and
$R(X_1)=R(X_1')$.

For $y$ an occurrence of $\vee$ in a form sequence $X$ such that
$X_1 y X_2$ is a subformsequence of $X$, we define the sequence
$T_y$ as $T(X_2)$ and $B_y$ as $B(X_1)$.

Note that if $x_1$ and $x_2$ are two different occurrences of
$\kon$ in $X$, then no letter is both in $L_{x_1}$ and $L_{x_2}$,
and the same holds for $R_{x_1}$ and $R_{x_2}$. We have an
analogous situation with $y$, $T$ and $B$.

For $X$ and $Y$ being form sequences, let ${f\!:X\str Y}$ be an
arrow of $\ACk^{st}$. (It is easy to see that $X$ is diversified
iff $Y$ is, and our assumption is that both are diversified.) By
considering $f$ in a developed form (see Section~2), and what is
intermuted by each of its factors, it can be checked easily that
the following conditions hold for $X$ and $Y$.

\prop{Condition $\kon$}{There is a function $c$ from the set of
occurrences of $\kon$ in $X$ onto the set of occurrences of $\kon$
in $Y$ such that for every occurrence $x$ of $\kon$ in $Y$ there
are occurrences $x_1,\ldots,x_k$ of $\kon$ in $X$ for which
\[L_x=L_{x_1}\ldots L_{x_k},\quad R_x=R_{x_1}\ldots R_{x_k}
\quad\mbox{\emph{and}}\quad c^{-1}(\{x\})=\{x_1,\ldots,x_k\}.
\]}
\vspace{-5ex}

\noindent We may say that $x_1,\ldots,x_k$ \emph{merge} into $x$.

\prop{Condition $\vee$}{There is a function $d$ from the set of
occurrences of $\vee$ in $Y$ onto the set of occurrences of $\vee$
in $X$ such that for every occurrence $y$ of $\vee$ in $X$ there
are occurrences $y_1,\ldots,y_l$ of $\vee$ in $Y$ for which
\[T_y=T_{y_1}\ldots T_{y_l},\quad B_y=B_{y_1}\ldots B_{y_l}
\quad\mbox{\emph{and}}\quad d^{-1}(\{y\})=\{y_1,\ldots,y_l\}.
\] }
\vspace{-5ex}

\noindent We may say that $y$ is \emph{split} into
$y_1,\ldots,y_l$.

For $X$ and $Y$ being form sequences, we say that $(X,Y)$ is a
\emph{legitimate pair} when ${let(X)}={let(Y)}$ and Conditions
$\kon$ and $\vee$ are satisfied.

So we know that if $X$ and $Y$ are respectively the source and the
target of an arrow of $\ACk^{st}$, then $(X,Y)$ is a legitimate
pair. Our purpose in this section is to show that the converse
holds too, which will give us a criterion for the existence of
arrows in $\ACk^{st}$. It is a decidable question whether $(X,Y)$
is a legitimate pair.

We introduce first the following notions. Let $X$ be a form
sequence, and let $p$ and $q$ be letters in $X$. Let \nad{q}{p} in
$X$ mean that there exists an occurrence $y$ of $\vee$ in $X$ such
that $q$ belongs to $B_y$ and $p$ belongs to $T_y$. If \nad{q}{p},
then in $\gamma(X)$ we have a horizontal dividing line segment to
which the top side of the rectangle enclosing $p$ and the bottom
side of the rectangle enclosing $q$ both belong. Let \nadd{}{} be
the transitive closure of the relation $\nad{}{}$ in $X$. Then we
have the following.

\prop{Lemma for \nadd{}{}}{Let $(X,Y)$ be a legitimate pair. If
\nadd{q}{p} in $Y$, then \nadd{q}{p} in~$X$.}

\noindent This is an easy consequence of Condition~$\vee$, and
becomes clear when we consider $\gamma(X)$ and $\gamma(Y)$. Note
that the converse of this lemma does not hold, but we will use
this lemma to establish a related equivalence in the
Position-Preservation Lemma below.

We define $q|p$ and $q|_*p$ analogously, and there is an analogous
lemma to the one above for $|_*$, which is a consequence of
Condition~$\kon$.

When there is no occurrence $x$ of $\kon$ in the form sequence $X$
such that $p$ belongs to $L_x$ we say that $p$ is a
\emph{left-border} letter of $X$. In the grid $\gamma(X)$, the
letter $p$ is at the left border of the rectangle corresponding to
$X$. We define analogously the \emph{right-border},
\emph{top-border} and \emph{bottom-border} letters of $X$; we just
replace $(L,\kon)$ by $(R,\kon)$, $(T,\vee)$ and $(B,\vee)$.

Let $y_1,\ldots,y_l$ for $l\geq 1$ be occurrences of $\vee$ in a
form sequence $X$, and let $T_{y_i}$, for $i\in\{1,\ldots,l\}$, be
the sequence $p_1^i\ldots p_{k_i}^i$ for $k_i\geq 1$. We say that
the sequence $y_1\ldots y_l$ is a \emph{horizontal transversal} of
$X$ when $p_1^1$ is a left-border letter of $X$ and $p_{k_l}^l$ is
a right-border letter of $X$, while for every $i\in\{1,\ldots,l\od
1\}$ we have $p_{k_i}^i|p_1^{i+1}$. (Horizontal transversals could
equivalently be defined by using $B_{y_i}$ instead of $T_{y_i}$.)
In $\gamma(X)$ a horizontal transversal appears as

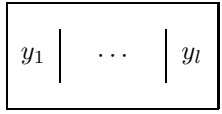
\begin{center}
\begin{picture}(80,60)

\put(0,10){\line(1,0){80}} \put(0,50){\line(1,0){80}}
\put(0,10){\line(0,1){40}} \put(80,10){\line(0,1){40}}
\put(20,20){\line(0,1){20}} \put(60,20){\line(0,1){20}}

\put(10,30){\makebox(0,0){$y_1$}}
\put(70,30){\makebox(0,0){$y_l$}}
\put(40,30){\makebox(0,0){$\cdots$}}

\end{picture}
\end{center}
where each $y_i$ is a horizontal line segment, and the rectangle
drawn corresponds to $X$. For example, if $y_1$, $y_2$ and $y_3$
are the three occurrences of $\vee$ (counting from the left) in
the example of Section~8, then the sequences $y_1$ and $y_2y_3$
are horizontal transversals. Contrary to what we have in this
example, horizontal transversals can be of length greater than 2,
and not every occurrence of $\vee$ in $X$ need be included in a
horizontal transversal. We can define analogously \emph{vertical
transversals}, and rely in the remainder of the exposition on this
other notion, rather than on the notion of horizontal transversal.

It is easy to ascertain the following lemma.

\prop{Transversal-Preservation Lemma}{Let $(X,Y)$ be a legitimate
pair. If $y_1\ldots y_l$ is a horizontal transversal of $X\!$,
then the members of
$\,\bigcup_{i\in\{1,\ldots,l\}}d^{-1}(\{y_i\})$ make a horizontal
transversal of $Y$.}

\noindent We denote by $d^{-1}(y_1\ldots y_l)$ the horizontal
transversal of $Y$ whose existence is claimed by this lemma.

We say that a letter $p$ is \emph{below} a horizontal transversal
$y_1\ldots y_l$ of a form sequence $X$ when there is a letter $q$
in $T_{y_i}$ for some $i\in\{1,\ldots,l\}$ such that either $p$ is
$q$ or \nadd{q}{p}. (We could say equivalently that there is a
letter $q$ in $B_{y_i}$ such that \nadd{q}{p}.) We say that an
occurrence $y$ of $\vee$ in a form sequence $X$ is \emph{below} a
horizontal transversal of $X$ when there is a letter $p$ below
this transversal such that $p$ is in $B_{y}$. These notions become
clear by considering $\gamma(X)$, where they really mean
``below''.

There are analogous definitions of being \emph{above} a horizontal
transversal for letters and occurrences of $\vee$. With vertical
transversals we would define being \emph{on the left-hand side}
and \emph{on the right-hand side} for letters and occurrences of
$\kon$. Any of these notions could be taken as central for the
exposition, as \emph{below} is in ours.

We can prove the following.

\prop{Position-Preservation Lemma}{Let $(X_1yX_2,Y)$ be a
legitimate pair for $y$ an occurrence of $\vee$. Then $p$ is below
the horizontal transversal $y$ of $X_1yX_2$ iff $p$ is below the
horizontal transversal $d^{-1}(y)$ of $Y$.}

\dkz From left to right we proceed by induction. If $p$ is in
$T_y$ in $X_1yX_2$, then we are done. If $p$, which is below $y$,
is not in $T_y$ in $X_1yX_2$, then there must be an occurrence
$y'$ of $\vee$ in $X_1yX_2$ below $y$, different from $y$, such
that $p$ is in $T_{y'}$, and by the induction hypothesis for every
$q$ in $B_{y'}$ we have that $q$ is below $d^{-1}(y)$. By
Condition~$\vee$ there is a sequence $y_1',\ldots,y_n'$ of
occurrences of $\vee$ in $Y$ such that for some
$i\in\{1,\ldots,n\}$ we have that $p$ is in $T_{y_i'}$ in $Y$, and
for every $q$ in $B_{y_i'}$ we have that $q$ is in $B_{y'}$. So,
by the induction hypothesis, for every $q$ in $B_{y_i'}$ we have
that $q$ is below $d^{-1}(y)$, and hence $p$ is below $d^{-1}(y)$.
From right to left we just apply the Lemma for \nadd{}{}. \qed

\noindent An analogous Position-Preservation Lemma holds for
\emph{above}, \emph{on the left-hand side} and \emph{on the
right-hand side}, defined as indicated above.

For $X$ a form sequence, let $p$ be a letter such that $X$ is not
$p$. We define inductively the form sequence $X^{-p}$:
\begin{itemize}
\item[] if $p$ is not in $X$, then $X^{-p}$ is $X$; \item[] for
${\!\ks\!\in\{\kon,\vee\}}$, if $X$ is of the form $Y\ks p$ or
$p\ks Y$, then $X^{-p}$ is $Y$, and if $X$ is of the form $Y_1\ks
Y_2$ for $p$ occurring in $Y_i$, for some $i\in\{1,2\}$, but
different from $Y_i$, then $X^{-p}$ is $Y_1^{-p}\ks Y_2^{-p}$.
\end{itemize}
Since the same $X$ can be of the form $Y_1 x Y_2$ and $Y_1'x'
Y_2'$ for $x$ and $x'$ different occurrences of $\!\ks\!$, and
$Y_i$ different from $Y_i'$, one can raise the question whether
the definition above is unambiguous. That it is indeed such can
easily be checked with the formal notation for form sequences of
\cite{DP04} (Section 6.2). For example, $(q\kon p\kon r)^{-p}$ is
$q\kon r$, irrespectively of whether we interpret $q\kon p\kon r$
as $(q\kon p)\kon r$ or as $q\kon(p\kon r)$.

Since it is easy to see that ${(X^{-p})^{-q}=(X^{-q})^{-p}}$, we
can define in the following manner $X^{-P}$ for $X$ a form
sequence and $P=\{p_1,\ldots,p_n\}$ a set of letters such that
${let(X)}-P\neq\emptyset$:
\[
X^{-P}=_{df}(\ldots(X^{-p_1})^{\ldots})^{-p_n}.
\]
For form sequences $X$ and $Y$, we abbreviate $X^{-let(Y)}$ by
$X^{-Y}$. We can then prove the following.

\prop{Interpolation Lemma $\vee$}{If $(X_1\vee X_2,Y)$ is a
legitimate pair, then ${(X_1,Y^{-X_2})}$, ${(X_2,Y^{-X_1})}$ and
$(Y^{-X_2}\vee Y^{-X_1},Y)$ are legitimate pairs.}

\dkz The proof of this lemma is pretty clear from $\gamma(X)$,
$\gamma(Y)$ and the Position-Preservation Lemma. We have the
following picture:

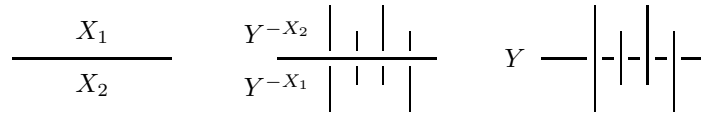
\begin{center}
\begin{picture}(260,60)

\put(0,30){\line(1,0){60}} \put(100,30){\line(1,0){60}}
\put(200,30){\line(1,0){17}} \put(223,30){\line(1,0){4}}
\put(233,30){\line(1,0){4}} \put(243,30){\line(1,0){4}}
\put(253,30){\line(1,0){7}}

\put(120,33){\line(0,1){17}} \put(120,27){\line(0,-1){17}}
\put(130,33){\line(0,1){7}} \put(130,27){\line(0,-1){7}}
\put(140,33){\line(0,1){17}} \put(140,27){\line(0,-1){7}}
\put(150,33){\line(0,1){7}} \put(150,27){\line(0,-1){17}}
\put(220,10){\line(0,1){40}} \put(230,20){\line(0,1){20}}
\put(240,20){\line(0,1){30}} \put(250,10){\line(0,1){30}}

\put(30,40){\makebox(0,0){$X_1$}}
\put(30,20){\makebox(0,0){$X_2$}}
\put(100,40){\makebox(0,0){$Y^{-X_2}$}}
\put(100,20){\makebox(0,0){$Y^{-X_1}$}}
\put(190,30){\makebox(0,0){$Y$}}

\end{picture}
\end{center}

Formally, for the legitimacy of ${(X_1,Y^{-X_2})}$, we demonstrate
what is the result $Z^{-X_2}$ of deleting the letters in
${let(X_2)}$ from a subformsequence $Z$ of $Y$. We distinguish
three cases depending on the form of $Z$. The first case is when
$Z$ is of the form ${Z_1y'Z_2}$ for $y'$ an occurrence of $\vee$
in $Y$ such that ${d(y')=y}$ where $X$ is $X_1yX_2$; then
$Z^{-X_2}$ is $Z_1$. In the second case, when $Z$ is of the form
$Z_1 y' Z_2$ such that ${d(y')\neq y}$, we have two subcases: if
$d(y')$ is below $y$, then $Z^{-X_2}$ is $Z_1^{-X_2}$, and if
$d(y')$ is above $y$, then $Z^{-X_2}$ is ${Z_1\vee Z_2^{-X_2}}$.
Finally, if $Z$ is of the form ${Z_1\kon Z_2}$, then $Z^{-X_2}$ is
${Z_1^{-X_2}\kon Z_2^{-X_2}}$.

We define the functions $c'$ and $d'$ that make legitimate
$(X_1,Y^{-X_2})$ in the following manner. First we define a
partial map $\kappa_{Z,X_2}$ from the set of occurrences of $\kon$
in a subformsequence $Z$ of $Y$ to the occurrences of $\kon$ in
$Z^{-X_2}$ by distinguishing the same three cases as above. When
from the domain of $\kappa_{Z,X_2}$ we omit the elements for which
$\kappa_{Z,X_2}$ is undefined, we obtain a bijection. The function
$c'$ is defined by composing the restriction of $c$ with
$\kappa_{Y,X_2}$. Next we define a map $\delta_{Z,X_2}$ from the
set of occurrences of $\vee$ in $Z^{-X_2}$ to the set of
occurrences of $\vee$ in $Z$ by distinguishing the same three
cases. There is also an obvious partial map $\varepsilon$ from the
set of occurrences of $\vee$ in ${X_1\vee X_2}$ to the set of
occurrences of $\vee$ in $X_1$. When from the domain of
$\varepsilon$ we omit the elements for which $\varepsilon$ is
undefined, we obtain a bijection. The function $d'$ is defined by
composing $\delta_{Y,X_2}$ with $d$ and $\varepsilon$.

We proceed analogously to show that ${(X_2,Y^{-X_1})}$ is
legitimate. To define the functions $c''$ and $d''$ that make
legitimate ${(Y^{-X_2}\vee Y^{-X_1},Y)}$ we use the inverses of
the maps $\kappa_{Y,X_2}$ and $\kappa_{Y,X_1}$ for $c''$, and the
maps $\delta_{Y,X_2}$ and $\delta_{Y,X_1}$ for $d''$. \qed

\noindent We can prove analogously the following dual lemma, where
$\vee$ is replaced by $\kon$.

\prop{Interpolation Lemma $\kon$}{If ${(X,Y_1\kon Y_2)}$ is a
legitimate pair, then ${(X^{-Y_2},Y_1)}$, ${(X^{-Y_1},Y_2)}$ and
${(X,X^{-Y_2}\kon X^{-Y_1})}$ are legitimate pairs.}

As a consequence of the foregoing results (cf.\ in particular the
first case in the proof of the Interpolation Lemma~$\vee$) we have
the following lemma, and its analogue immediately below.

\prop{Auxiliary Lemma $\vee$}{If ${(X,Y_1yY_2)}$ is a legitimate
pair for $y$ an occurrence of $\vee$, then $X$ is of the form
${X_1d(y)X_2}$ and ${(X_i,Y_i)}$ is a legitimate pair for every
${i\in\{1,2\}}$.} \vspace{-2ex}

\prop{Auxiliary Lemma $\kon$}{If ${(X_1xX_2,Y)}$ is a legitimate
pair for $x$ an occurrence of $\kon$, then $Y$ is of the form
${Y_1c(x)Y_2}$ and ${(X_i,Y_i)}$ is a legitimate pair for every
${i\in\{1,2\}}$.}

Then we can prove the following.

\prop{Theorem}{There is an arrow ${f\!:X\str Y}$ of $\ACk^{st}$
iff ${(X,Y)}$ is a legitimate pair.}

\dkz We have already established this theorem from left to right.
From right to left we proceed by induction on the sum $n$ of the
number of letters in $X$ and the number of occurrences of $\kon$
in $X$. If ${n=1}$, then $X$ and $Y$ are the same letter $p$.

If ${n>1}$ and $Y$ is ${Y_1\vee Y_2}$ or $X$ is ${X_1\kon X_2}$,
then we apply the Auxiliary Lemmata~$\vee$ and $\kon$ above.

If ${n>1}$ and $X$ is ${X_1\vee X_2}$, while $Y$ is ${Y'xY''}$,
for $x$ an occurrence of $\kon$, then we proceed as follows. By
applying the Interpolation Lemma~$\vee$ we obtain that
${(X_i,Y^{-X_{3-i}})}$, for every ${i\in\{1,2\}}$, and
${(Y^{-X_2}\vee Y^{-X_1},Y)}$ are legitimate pairs. By the
induction hypothesis, we obtain the arrows ${f_i\!:X_i\str
Y^{-X_{3-i}}}$, and hence the arrow ${f_1\vee f_2\!:X_1\vee
X_2\str Y^{-X_2}\vee Y^{-X_1}}$.

If either $f_1$ or $f_2$ is not an identity arrow, then
${Y^{-X_2}\vee Y^{-X_1}}$ has at least one occurrence of $\kon$
less than $X$, and we may apply the induction hypothesis to obtain
an arrow ${g\!:Y^{-X_2}\vee Y^{-X_1}\str Y}$. Hence we have
${g\cirk(f_1\vee f_2)\!:X\str Y}$.

If both $f_1$ and $f_2$ are identity arrows, then $Y^{-X_{3-i}}$
is of the form ${Y_i'x_iY_i''}$ for $x_i$ an occurrence of $\kon$
such that ${\kappa_{Y,X_{3-i}}(x)=x_i}$, where
$\kappa_{Y,X_{3-i}}$ is defined as in the proof of the
Interpolation Lemma~$\vee$. We have the arrow
\[
c^k_{Y_1',Y_1'',Y_2',Y_2''}\!:(Y_1'\kon Y_1'')\vee(Y_2'\kon
Y_2'')\str(Y_1'\vee Y_2')\kon(Y_1''\kon Y_2''),
\]
and we have to ascertain that ${(Y_1^j\vee Y_2^j,Y^j)}$ is a
legitimate pair for every ${j\in\{',''\}}$. This follows from the
Interpolation Lemma~$\kon$.

In terms of rectangular grids we have

\begin{center}
\begin{picture}(250,85)

\put(0,40){\line(1,0){60}} \put(100,40){\line(1,0){27}}
\put(133,40){\line(1,0){27}} \put(200,40){\line(1,0){7}}
\put(213,40){\line(1,0){14}} \put(233,40){\line(1,0){4}}
\put(243,40){\line(1,0){4}} \put(253,40){\line(1,0){7}}

\put(10,43){\line(0,1){7}}

\put(10,37){\line(0,-1){7}}

\put(30,43){\line(0,1){27}}

\put(30,37){\line(0,-1){27}}

\put(40,43){\line(0,1){7}}

\put(40,37){\line(0,-1){7}}

\put(50,43){\line(0,1){17}}

\put(50,37){\line(0,-1){7}}

\put(110,43){\line(0,1){7}}

\put(110,37){\line(0,-1){7}}

\put(130,10){\line(0,1){60}}

\put(140,43){\line(0,1){7}} \put(140,37){\line(0,-1){7}}
\put(150,43){\line(0,1){17}} \put(150,37){\line(0,-1){7}}
\put(210,30){\line(0,1){20}} \put(230,10){\line(0,1){60}}
\put(240,30){\line(0,1){20}} \put(250,30){\line(0,1){30}}

%%%%%%%%%%%%%%%%%%%%%%%%%%%%

\put(20,60){\makebox(0,0){$Y_1'$}}

\put(20,20){\makebox(0,0){$Y_2'$}}

\put(43,60){\makebox(0,0){$Y_1''$}}

\put(43,20){\makebox(0,0){$Y_2''$}}

%%%%%%%%%%%%%%%%%%%

\put(120,60){\makebox(0,0){$Y_1'$}}

\put(120,20){\makebox(0,0){$Y_2'$}}

\put(143,60){\makebox(0,0){$Y_1''$}}

\put(143,20){\makebox(0,0){$Y_2''$}}

%%%%%%%%%%%%%%%%%%%%

\put(-15,62){\makebox(0,0){$Y^{-X_2}$}}

\put(-15,22){\makebox(0,0){$Y^{-X_1}$}}

\put(220,77){\makebox(0,0){$Y'$}}

\put(243,77){\makebox(0,0){$Y''$}}

%%%%%%%%%%%%%%%%%%%%%%%%

\put(33,73){\makebox(0,0){$x_1$}}

\put(33,5){\makebox(0,0){$x_2$}}

\end{picture}
\end{center}

By applying the induction hypothesis we obtain the arrows
${f^j\!:Y_1^j\vee Y_2^j\str Y^j}$ for ${j\in\{',''\}}$. Hence we
have ${(f'\kon f'')\cirk c^k_{Y_1',Y_1'',Y_2',Y_2''}\!:X\str Y}$.
\qed

\vspace{-1ex}

\section{\large\bf Biassociative intermuting coherence}
Suppose ${(X,Y)}$ is a legitimate pair, $x$ is an occurrence of
$\kon$ in $Y$, and $y$ is an occurrence of $\vee$ in X. We say
that ${(x,y)}$ is a \emph{crossing} of ${(X,Y)}$ when the sets
${R_x\cap B_y}$, ${L_x\cap B_y}$, ${R_x\cap T_y}$ and ${L_x\cap
T_y}$ are all nonempty.

If these four sets are nonempty, then they are singletons. This is
because for every pair of letters $p$ and $q$ in a form sequence
there is a minimal subformsequence in which both $p$ and $q$
occur, and this subformsequence is either a conjunction or a
disjunction, but not both.

The arrow terms of \ACk\ where indices are replaced by form
sequences are the arrow terms of $\ACk^{st}$ (see \cite{DP04},
Section 3.2). In these arrow terms of $\ACk^{st}$ we may replace
all the $b$-arrows by identity arrows.

For ${f\!:X\str Y}$ an arrow term of $\ACk^{st}$ there is a
natural one-to-one correspondence between the occurrences of $c^k$
in $f$ and the crossings of ${(X,Y)}$. If the representatives for
the rectangular grids $\gamma(X)$ and $\gamma(Y)$ are well chosen,
then $c^k_{S,T,U,V}$ in $f$ corresponds to the crossing in the
following two fragments of these grids:

\begin{center}
\begin{picture}(160,60)

\put(0,30){\line(1,0){60}} \put(100,30){\line(1,0){27}}
\put(133,30){\line(1,0){27}} \put(30,10){\line(0,1){17}}
\put(30,33){\line(0,1){17}} \put(130,10){\line(0,1){40}}

\put(-1,29){\makebox(0,0)[r]{$y$}}
\put(130.5,53){\makebox(0,0)[b]{$x$}}
\put(15,45){\makebox(0,0){$S$}} \put(45,45){\makebox(0,0){$T$}}
\put(15,15){\makebox(0,0){$U$}} \put(45,15){\makebox(0,0){$V$}}

\put(115,45){\makebox(0,0){$S$}} \put(145,45){\makebox(0,0){$T$}}
\put(115,15){\makebox(0,0){$U$}} \put(145,15){\makebox(0,0){$V$}}

\end{picture}
\end{center}
These two pictures explain why we call ${(x,y)}$ a crossing. If
these representatives are not well chosen, we may have, for
example,
\begin{center}
\begin{picture}(160,60)

\put(0,30){\line(1,0){60}} \put(100,35){\line(1,0){27}}
\put(133,25){\line(1,0){27}} \put(20,10){\line(0,1){17}}
\put(40,33){\line(0,1){17}} \put(130,7){\line(0,1){43}}

\put(-1,29){\makebox(0,0)[r]{$y$}}
\put(130.5,53){\makebox(0,0)[b]{$x$}}
\put(20,45){\makebox(0,0){$S$}} \put(50,45){\makebox(0,0){$T$}}
\put(10,15){\makebox(0,0){$U$}} \put(40,15){\makebox(0,0){$V$}}

\put(115,48){\makebox(0,0){$S$}} \put(145,42){\makebox(0,0){$T$}}
\put(115,18){\makebox(0,0){$U$}} \put(145,12){\makebox(0,0){$V$}}

\end{picture}
\end{center}

We say that ${(x,y)}$ is the pair of \emph{coordinates} of
$c^k_{S,T,U,V}$ in ${(X,Y)}$. For a $c^k$-term $f$, we say that
the \emph{coordinates} of $f$ are the coordinates of the head
of~$f$. (The notions of $c^k$-term and of its head are defined in
Section~2.)

If ${f,g\!:X\str Y}$ are two arrow terms of $\ACk^{st}$, then
there is one-to-one correspondence between the occurrences of
$c^k$ in $f$ and the occurrences of $c^k$ in $g$, such that the
occurrences with the same coordinates correspond to each other.

Let ${f_2\cirk f_1}$ be a composition of two $c^k$-terms that is
not equal to an arrow term of the form ${f_1'\cirk f_2'}$ for
$f_i'$, where ${i\in\{1,2\}}$, a $c^k$-term with the same
coordinates as $f_i$. Then we say that $f_2$ is \emph{checked} by
$f_1$.

It is straightforward to verify the following lemma by appealing
to bifunctorial and naturality equations.

\prop{Checking Lemma}{Let ${f_2\cirk f_1}$ be a composition of two
$c^k$-terms such that $f_i$, where ${i\in\{1,2\}}$, has the
coordinates ${(x_i,y_i)}$ with ${x_1\neq x_2}$ and ${y_1\neq
y_2}$. Then $f_2$ is not checked by $f_1$.}

As a corollary of this lemma we have that if $f_2$ is checked by
$f_1$ and ${(x_i,y_i)}$ are the coordinates of $f_i$, then either
${x_1=x_2}$, and we say that $f_2$ is \emph{vertically checked} by
$f_1$, or ${y_1=y_2}$ and we say that $f_2$ is \emph{horizontally
checked} by $f_1$.

For ${f_n\cirk\ldots\cirk f_1}$, with ${n\geq 1}$, a composition
of $c^k$-terms, we say that $f_i$ is \emph{horizontally blocked}
in this composition when either ${i=1}$ or $f_i$ is horizontally
checked by $f_{i-1}$.

We say that a $\beta$-term $f$ is \emph{conjunctively headed} when
its head $\beta$ occurs in a subterm of $f$ of the form
${\mj_X\kon\beta}$ or ${\beta\kon\mj_X}$. We are interested in
this notion when $\beta$ is $c^k$. We can prove the following
series of lemmata, leading to the theorem below.

\prop{Horizontal Checking Lemma}{If ${f_2\cirk f_1}$ is a
composition of two $c^k$-terms such that $f_2$ is conjunctively
headed, then $f_2$ is not horizontally checked by $f_1$.}

\dkz Suppose $f_2$ is conjunctively headed, and the coordinates of
$f_i$, for ${i\in\{1,2\}}$, are ${(x_i,y_i)}$. If ${y_1\neq y_2}$,
then $f_2$ is not horizontally checked by $f_1$. If ${y_1=y_2=y}$,
then we have to consider many cases, for which the following three
fragments of rectangular grids corresponding to the source of
$f_1$ are typical:
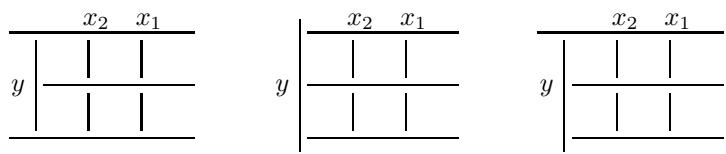
\begin{center}
\begin{picture}(295,55)

\put(10,5){\line(1,0){70}} \put(10,45){\line(1,0){70}}
\put(23,25){\line(1,0){57}} \put(20,8){\line(0,1){34}}
\put(40,8){\line(0,1){14}} \put(40,28){\line(0,1){14}}
\put(60,8){\line(0,1){14}} \put(60,28){\line(0,1){14}}

\put(16,24){\makebox(0,0)[r]{$y$}}
\put(43,48){\makebox(0,0)[b]{$x_2$}}
\put(63,48){\makebox(0,0)[b]{$x_1$}}

\put(123,5){\line(1,0){57}} \put(123,45){\line(1,0){57}}
\put(123,25){\line(1,0){57}} \put(120,0){\line(0,1){50}}
\put(140,8){\line(0,1){14}} \put(140,28){\line(0,1){14}}
\put(160,8){\line(0,1){14}} \put(160,28){\line(0,1){14}}

\put(116,24){\makebox(0,0)[r]{$y$}}
\put(143,48){\makebox(0,0)[b]{$x_2$}}
\put(163,48){\makebox(0,0)[b]{$x_1$}}

\put(223,5){\line(1,0){57}} \put(210,45){\line(1,0){70}}
\put(223,25){\line(1,0){57}} \put(220,0){\line(0,1){42}}
\put(240,8){\line(0,1){14}} \put(240,28){\line(0,1){14}}
\put(260,8){\line(0,1){14}} \put(260,28){\line(0,1){14}}

\put(216,24){\makebox(0,0)[r]{$y$}}
\put(243,48){\makebox(0,0)[b]{$x_2$}}
\put(263,48){\makebox(0,0)[b]{$x_1$}}

\end{picture}
\end{center}
We can always apply the equation $(\overline{\psi}\hat{b})$ of
Section~4. Note that in the second and third case we could not
have
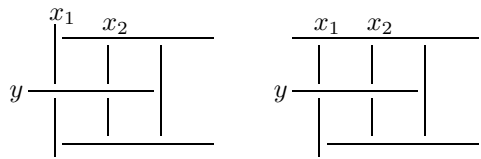
\begin{center}
\begin{picture}(295,55)

\put(123,5){\line(1,0){57}} \put(123,45){\line(1,0){57}}
\put(110,25){\line(1,0){47}} \put(120,0){\line(0,1){22}}
\put(120,28){\line(0,1){22}} \put(140,8){\line(0,1){14}}
\put(140,28){\line(0,1){14}} \put(160,8){\line(0,1){34}}

\put(108,24){\makebox(0,0)[r]{$y$}}
\put(143,48){\makebox(0,0)[b]{$x_2$}}
\put(123,52){\makebox(0,0)[b]{$x_1$}}

\put(223,5){\line(1,0){57}} \put(210,45){\line(1,0){70}}
\put(210,25){\line(1,0){47}} \put(220,0){\line(0,1){22}}
\put(220,28){\line(0,1){14}} \put(240,8){\line(0,1){14}}
\put(240,28){\line(0,1){14}} \put(260,8){\line(0,1){34}}

\put(208,24){\makebox(0,0)[r]{$y$}}
\put(243,48){\makebox(0,0)[b]{$x_2$}}
\put(223,48){\makebox(0,0)[b]{$x_1$}}

\end{picture}
\end{center}
because then the source of $f_1$ would not be a form sequence.
\qed

\vspace{-2ex}

\prop{Horizontal Blocking Corollary}{Let ${f_n\cirk\ldots\cirk
f_1\!:X_1yX_2\str Y}$, for ${n\geq 1}$ and $y$ an occurrence of
$\vee$, be a composition of $c^k$-terms, and let
${i\in\{1,\ldots,n\}}$ be such that the second coordinate of $f_i$
is $y$. If $f_i$ is conjunctively headed, then $f_i$ is not
horizontally blocked in ${f_n\cirk\ldots\cirk f_1}$.}

\vspace{-2ex}

\prop{Vertical Checking Lemma}{If ${f_2\cirk f_1}$ is a
composition of two $c^k$-terms such that $f_1$ is not
conjunctively headed, then $f_2$ is not vertically checked by
$f_1$.}

\dkz Suppose $f_1$ is not conjunctively headed, and the
coordinates of $f_i$, for ${i\in\{1,2\}}$, are ${(x_i,y_i)}$. If
${x_1\neq x_2}$, then $f_2$ is not vertically checked by $f_1$. If
${x_1=x_2=x}$, then we have a situation as in the following
fragment of a rectangular grid corresponding to the source of
$f_1$:
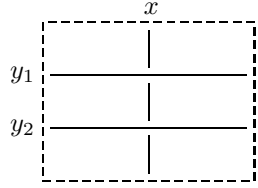
\begin{center}
\begin{picture}(80,75)

\multiput(0,5)(5,0){16}{\line(1,0){3}}
\multiput(80,5)(0,5){12}{\line(0,1){3}}
\multiput(80,65)(-5,0){16}{\line(-1,0){3}}
\multiput(0,65)(0,-5){12}{\line(0,-1){3}}

\put(3,25){\line(1,0){74}} \put(3,45){\line(1,0){74}}
\put(40,8){\line(0,1){14}} \put(40,28){\line(0,1){14}}
\put(40,48){\line(0,1){14}}

\put(-3,25){\makebox(0,0)[r]{$y_2$}}
\put(-3,45){\makebox(0,0)[r]{$y_1$}}
\put(41,68){\makebox(0,0)[b]{$x$}}

\end{picture}
\end{center}
We can always apply the equation $(\psi\check{b})$ of Section~4.
\qed

\vspace{-2ex}

\prop{Permuting Lemma}{Let ${g\cirk f\!:X_1yX_2\str Y}$, for $y$
an occurrence of $\vee$, be such that $g$ is a $c^k$-term whose
second coordinate is different from $y$, and every occurrence of
$c^k$ in $f$ has the second coordinate $y$. Then ${g\cirk
f=f'\cirk g'}$ where $g'$ has the same coordinates as $g$.}

\dkz Let ${f_n\cirk\ldots\cirk f_1\cirk\mj_{X_1yX_2}}$ be a
developed arrow term equal to $f$ (see Section~2 for the notion of
developed arrow term). We proceed by induction on $n$. If ${n=0}$,
then ${g\cirk\mj_{X_1yX_2}=\mj_{X_1yX_2}\cirk g}$ and ${g'=g}$.

Let ${n>0}$. Suppose the coordinates of $g$ are ${(x,y_1)}$. If
there is an $i\in\{1,\ldots,n\}$ such that $f_i$ has the
coordinates ${(x,y)}$, then we can assume that
${f_n\cirk\ldots\cirk f_1}$ is such that $f_i$ is horizontally
blocked in this composition. Suppose $g$ is checked by $f_n$.
Then, since ${y_1\neq y}$, the $c^k$-term $g$ must be vertically
checked by $f_n$. Hence the coordinates of $f_n$ are ${(x,y)}$. By
the Vertical Checking Lemma, $f_n$ is conjunctively headed. By the
Horizontal Blocking Corollary, $f_n$ is not horizontally blocked
in the composition ${f_n\cirk\ldots\cirk f_1}$, which contradicts
the assumption. So, ${g\cirk f_n=f_n'\cirk g'}$ where $g'$ has the
same coordinates as $g$, and we may apply the induction
hypothesis. \qed

\vspace{-2ex}

\prop{Theorem}{If ${f',f\!:X\str Y}$ are arrows of $\ACk^{st}$,
then ${f=f'}$.}

\dkz As in the proof of the Theorem of the preceding section, we
proceed by induction on the sum $n$ of the number of letters in
$X$ and the number of occurrences of $\kon$ in $X$. If ${n=1}$,
then ${X=Y=p}$, and ${\mj_p\!:p\str p}$ is the unique arrow from
$p$ to $p$.

If $X$ is ${X_1\kon X_2}$, then $Y$ must be of the form ${Y_1\kon
Y_2}$, and it is easy to prove that every ${f\!:X\str Y}$ is equal
to an arrow ${f_1\kon f_2}$ for ${f_i\!:X_i\str Y_i}$, for
${i\in\{1,2\}}$. Then we may apply the induction hypothesis to
$f_i$. We proceed analogously when $Y$ is ${Y_1\vee Y_2}$.

Suppose now that $X$ is ${X_1yX_2}$ for $y$ an occurrence of
$\vee$, and $Y$ is ${Y_1xY_2}$ for $x$ an occurrence of $\kon$. By
the Permuting Lemma, we may conclude that ${f=h\cirk g}$ where the
second coordinate of every occurrence of $c^k$ in $h$ is $y$, and
there is no occurrence of $c^k$ in $g$ with the second coordinate
$y$. Then it is easy to see that $g$ must be equal to an arrow
${g_1\vee g_2\!:X_1yX_2\str Y^{-X_2}\vee Y^{-X_1}}$ (cf.\ the
Interpolation Lemma in the preceding section). If ${g_1\vee g_2}$
is not an identity arrow, then we may apply the induction
hypothesis to $g_i$, for ${i\in\{1,2\}}$, and $h$. This is because
the ``interpolated'' ${Y^{-X_2}\vee Y^{-X_1}}$ is uniquely
determined by $X$ and $Y$, and is hence the same for $f$ and $f'$,
and because in the source of $h$ there is at least one occurrence
of $\kon$ less than in $X$. If, on the other hand, ${g_1\vee g_2}$
is the identity arrow $\mj_{X_1yX_2}$, then
${h=h_m\cirk\ldots\cirk h_1}$, for ${m\geq 1}$, where
${h_m\cirk\ldots\cirk h_1}$ is a composition of $c^k$-terms, and
for some ${j\in\{1,\ldots,m\}}$ we have that $h_j$ has the
coordinates ${(x,y)}$. Then, since $x$ is the main connective of
$Y$, we have that $h_j$ is either $h_1$ or it is conjunctively
headed. In the second case, by the Horizontal Checking Lemma we
may again assume that $h_j$ is $h_1$. Since $h_1$ is of the type
${X_1yX_2\str Y_1'\kon Y_2'}$ (it is a $c^k$-term), we may apply
the induction hypothesis to ${h_m\cirk\ldots\cirk h_2\!:Y_1'\kon
Y_2'\str Y_1xY_2}$, in whose source there is one occurrence of
$\kon$ less than in $X$. \qed

\noindent The proof of the foregoing theorem suggests how to
construct a unique normal form for arrow terms.

The Theorem we have proved above does not amount to showing that
$\ACk^{st}$ is a preorder because it is formulated with the
assumption that $X$ and $Y$ are diversified. We can however lift
this assumption. To achieve that we need the following.

We say that the sequence of letters ${p_1\ldots p_n}$ for ${n\geq
1}$ is the \emph{left border} of a form sequence $X$ when
${\{p_1,\ldots,p_n\}}$ is the set of all the left-border letters
of $X$ and for every ${i\in\{1,\ldots,n\od 1\}}$ we have
\naad{p_1}{p_{i+1}}$\;$ (see the preceding section for the
definition of left-border letter and of the relation \nad{}{}).
The left border of $X$ appears indeed as the left border of the
rectangular grid ${\gamma(X)}$. We could define analogously the
\emph{right border}, the \emph{top border} and the \emph{bottom
border}. We can then prove the following.

\prop{Border Lemma}{If ${(X,Y)}$ is a legitimate pair, then
${p_1\ldots p_n}$ is the left border of $X$ iff ${p_1\ldots p_n}$
is the left border of $Y$.}

\dkz We establish first that for a legitimate pair ${(X,Y)}$ we
have that $p$ is a left-borer letter of $X$ iff $p$ is a
left-border letter of $Y$. From left to right we appeal to the
fact that the function $c$ of Condition~$\kon$ is onto. For the
other direction we appeal to the fact that $c$ is a function. The
matter becomes clear by considering the rectangular grids
${\gamma(X)}$ and ${\gamma(Y)}$. So ${\{p_1,\ldots,p_n\}}$ is the
set of left-border letters of $X$ iff ${\{p_1,\ldots,p_n\}}$ is
the set of left-border letters of $Y$.

Suppose ${p_1\ldots p_n}$ is the left border of $Y$, and suppose
${n>1}$. Then by Condition~$\vee$ from \naad{p_i}{p_{i+1}}$\;$ in
$Y$ we may infer \naad{p_i}{p_{i+1}}$\;$ in $X$. So ${p_1\ldots
p_n}$ is the left border of $X$. \qed

The same lemma holds when we replace ``left'' by ``right'',
``top'' and ``bottom''. For ``top'' and ``bottom'' we appeal to
the function $d$ of Condition~$\vee$ in the first part of the
proof, and to Condition~$\kon$ in the second part. We can prove
the following.

\prop{Uniqueness Lemma}{If ${f_1\!:X\str Y_1}$ and ${f_2\!:X\str
Y_2}$ are arrow terms of $\ACk^{st}$ such that $Y_1$ and $Y_2$ are
substitution instances of each other, then $Y_1$ and $Y_2$ are the
same form sequence.}

\dkz The form sequences $Y_1$ and $Y_2$ are substitution instances
of each other iff their rectangular grids ${\gamma(Y_1)}$ and
${\gamma(Y_2)}$ are the same when letters are omitted. We will
show that $X$ imposes the letters to be put in this grid. So $Y_1$
and $Y_2$ will be the same form sequence.

Since ${(X,Y_1)}$ and ${(X,Y_2)}$ are legitimate pairs, by the
Border Lemma, the left border ${p_1\ldots p_n}$ of $X$ is also the
left border of $Y_1$ and $Y_2$. If this left border is also the
right border, then we are done. Otherwise, there must be an
occurrence $x$ of $\kon$ in $Y_1$ and $Y_2$ such that the sequence
${p_1\ldots p_n}$ is
\[
p_1\ldots p_{j-1}R_xp_{j+m}\ldots p_n
\]
for some ${m\geq 1}$. Then, by Condition~$\kon$, for some ${k\geq
1}$ there are occurrences ${x_1,\ldots,x_k}$ of $\kon$ in $X$ such
that ${R_x=R_{x_1}\ldots R_{x_k}}$. So $Y_1$ and $Y_2$ must be
such that ${L_x=L_{x_1}\ldots L_{x_k}}$.

If the sequence ${p_1\ldots p_{j-1}L_xp_{j+m}\ldots p_n}$ is the
right border of $Y_1$ and $Y_2$, then we are done. Otherwise,
there is in $Y_1$ and $Y_2$ an occurrence $x'$ of $\kon$ such that
this sequence is of the form
\[
p_1'\ldots p_{j'-1}'R_{x'}p_{j'+m'}'\ldots p_{n'}',
\]
and we proceed as before until we are done. \qed

We are now ready to prove the following.

\prop{$\ACk^{st}$ Coherence}{The category $\ACk^{st}$ is a
preorder.}

\dkz For ${f_1,f_2\!:X\str Y}$ arrows of $\ACk^{st}$ with $X$ and
$Y$ not diversified, we find ${f_1'\!:X'\str Y_1}$ and
${f_2'\!:X'\str Y_2}$ such that $X'$, and hence also $Y_1$ and
$Y_2$, are diversified, while $f_1$ and $f_2$ are substitution
instances of $f_1'$ and $f_2'$ respectively. Then by the
Uniqueness Lemma we obtain that $Y_1$ is $Y_2$, and by the Theorem
established in this section we have that ${f_1'=f_2'}$. So
${f_1=f_2}$. \qed

As a consequence of this and of the equivalence of the categories
\ACk\ and $\ACk^{st}$ we obtain the following.

\prop{Biassociative Intermuting Coherence}{The category \ACk\ is a
preorder.}

\section{\large\bf Bimonoidal intermuting categories}
In this section we add unit objects to the biassociative
intermuting categories of the preceding two sections, and prove a
restricted coherence result for the ensuing notion (more general
than the notion of two-fold monoidal category of \cite{BFSV}; see
the next section). A \emph{bimonoidal intermuting} category is a
$\kappa$-normal bimonoidal category ${\langle\aA,
\kon,\vee,\top,\bot\rangle}$ (see Section~6) such that
${\langle\aA,\kon,\vee\rangle}$ is a biassociative intermuting
category (see Section 10), and, moreover, the equations
$(c^k\delta\sigma)$ and $(c^k\kappa)$ (see Section~4) are
satisfied. This means that besides the assumptions for bimonoidal
categories (see Section~2) we have the natural transformation
$c^k$, the isomorphisms
${\hat{w}^{\rts}_\bot\!:\bot\str\bot\kon\bot}$ and
${\check{w}^{\str}_\top\!:\top\vee\top\str\top}$, and the arrow
${\kappa\!:\bot\str\top}$, for which we assume the equations
$(c^kb)$, $(c^k\delta\sigma)$, $(wb)$, $(\kappa\delta\sigma)$ and
$(c^k\kappa)$ of Section~4.

Let $\mbox{\bf AC}^{\bf k}_{\top,\bot}$ and $\mbox{\bf AC}^{{\bf
k}\emptyset}_{\top,\bot}$ be the free bimonoidal intermuting
categories generated respectively by a nonempty set of objects and
the empty set of objects. In $\mbox{\bf AC}^{\bf k}_{\top,\bot}$
there are no arrows of the type ${A\str B}$ such that in one of
$A$ and $B$ there are no letters and in the other there are. The
coherence result we will prove below for $\mbox{\bf AC}^{{\bf
k}\emptyset}_{\top,\bot}$ will enable us to ascertain that
$\mbox{\bf AC}^{{\bf k}\emptyset}_{\top,\bot}$ is a full
subcategory of $\mbox{\bf AC}^{\bf k}_{\top,\bot}$.

According to what we had for the category $\mn_{\top,\bot}$ in
Section~5, there is an isomorphism of the type ${A\str\nu(A)}$ in
$\mbox{\bf AC}^{\bf k}_{\top,\bot}$. In $\mbox{\bf AC}^{{\bf
k}\emptyset}_{\top,\bot}$ the formula $\nu(A)$ is either $\top$ or
$\bot$. We can prove the following.

\prop{$\mbox{\bf AC}^{{\bf k}\emptyset}_{\top,\bot}$ Coherence}
{The category $\mbox{\bf AC}^{{\bf k}\emptyset}_{\top,\bot}$ is a
preorder.}

\dkz We enlarge the proof of $\mka_{\top,\bot}^\emptyset$
Coherence of Section~7. We now have additional equations obtained
from the equations $(c^k\delta\sigma)$ and $(c^k\kappa)$, which
enable us to eliminate every occurrence of $c^k$, together with
equations derived from the naturality equations for $c^k$. \qed

We say that an object $A$ of $\mbox{\bf AC}^{\bf k}_{\top,\bot}$
is $\zeta$-\emph{pure} for ${\zeta\in\{\top,\bot\}}$ when there is
no occurrence of $\zeta$ in $\nu(A)$. For
${(\zeta,\!\ks\!)\in\{(\bot,\kon),(\top,\vee)\}}$, it is easy to
see that $A$ is not $\zeta$-pure iff either $\nu(A)$ is $\zeta$ or
there is subformula of $A$ of the form ${B_1\ks B_2}$ such that
for some ${i\in\{1,2\}}$ we have that $\nu(B_i)$ is $\zeta$ and a
letter occurs in $B_{3-i}$. A formula is called \emph{pure} when
it is both $\bot$-pure and $\top$-pure. We can prove the
following.

\prop{Purity Lemma}{Let ${f\!:A\str B}$ be an arrow of $\mbox{\bf
AC}^{\bf k}_{\top,\bot}$. If $A$ is $\bot$-pure, then $B$ is
$\bot$-pure, and if $B$ is $\top$-pure, then $A$ is $\top$-pure.}

\dkz By the Development Lemma (see Section~2), it is enough to
verify the lemma for ${f\!:A\str B}$ a $\beta$-term of $\mbox{\bf
AC}^{\bf k}_{\top,\bot}$, where $\beta$ is $\b{\xi}{\str}$,
$\b{\xi}{\rts}$, $\d{\xi}{\str}$, $\d{\xi}{\rts}$,
$\s{\xi}{\str}$, $\s{\xi}{\rts}$, $\w{\xi}{\str}_\zeta$,
$\w{\xi}{\rts}_\zeta$, $\kappa$ or $c^k$. The only interesting
case arises when ${f\!:A\str B}$ is a $c^k$-term.

Suppose $B$ is not $\bot$-pure. If $\nu(B)$ is $\bot$, then we can
easily conclude that $\nu(A)$ is $\bot$ too. Suppose, on the other
hand, that there is a subformula of $B$ of the form ${KxD}$ or
${DxK}$ for $x$ an occurrence of $\kon$ such that $\nu(K)$ is
$\bot$ and a letter occurs in $D$. Suppose the head of $f$ is
\[
c^k_{E,F,G,H}\!:(E\kon F)\vee(G\kon H)\str(E\vee G)\kon(F\vee H).
\]
The only interesting case is when $x$ is the main $\kon$ of
${(E\vee G)\kon(F\vee H)}$. Suppose $K$ is ${E\vee G}$. Then
${\nu(E)=\nu(G)=\bot}$, and there is a letter in either $F$ or
$H$. So $A$ is not $\bot$-pure. We reason analogously in other
cases. \qed

As an immediate corollary of this lemma we have the following.

\prop{Purity Corollary}{If ${f\!:A\str B}$ and ${g\!:B\str C}$ are
two arrows of $\mbox{\bf AC}^{\bf k}_{\top,\bot}$ such that both
$A$ and $C$ are pure, then $B$ is pure.}

We say that a formula is \emph{constant-free} when there is no
occurrence of $\top$ or $\bot$ in it. For
${\!\ks\!\in\{\kon,\vee\}}$, the arrow terms
$\b{\xi}{\str}_{A,B,C}$, $\b{\xi}{\rts}_{A,B,C}$ and
$c^k_{A,B,C,D}$ are \emph{constant-free} when their indices $A$,
$B$, $C$ and $D$ are constant-free. We can then prove the
following.

\prop{Constant-Free Indices Lemma}{If ${f\!:A\str B}$ is an arrow
of $\mbox{\bf AC}^{\bf k}_{\top,\bot}$ such that $A$ and $B$ are
pure, then there is an arrow term ${f'\!:A\str B}$ of $\mbox{\bf
AC}^{\bf k}_{\top,\bot}$ in which every occurrence of
$\b{\xi}{\str}$, $\b{\xi}{\rts}$ and $c^k$ is in a constant-free
arrow term, $\kappa$ does not occur in $f'$ and ${f=f'}$ in
$\mbox{\bf AC}^{\bf k}_{\top,\bot}$.}

\dkz We call the assumption that $A$ and $B$ are pure \emph{the
purity assumption}. By the Development Lemma (see Section~2) and
the Purity Corollary, it is enough to verify the lemma for
${f\!:A\str B}$ a $\beta$-term of $\mbox{\bf AC}^{\bf
k}_{\top,\bot}$. By reasoning as in the proof of
$\mka_{\top,\bot}^\emptyset$ Coherence in Section~7, we may assume
that every index of every $\beta$-term, for $\beta$ being
$\b{\xi}{\str}$, $\b{\xi}{\rts}$ or $c^k$, is $\nu(D)$ for some
$D$; so each of these indices is either $\top$ or $\bot$, or it is
constant-free.

Suppose $\beta$ is $\hat{b}^\str$. If any of its indices is
$\top$, then we eliminate $\hat{b}^\str$ by relying on the
equations mentioned in the proof of $\mka_{\top,\bot}^\emptyset$
Coherence. If all its indices are $\bot$, then we eliminate
$\hat{b}^\str$ by relying on the ${(wb)}$ equation
$(\pp{\bot}\hat{b})$ of Section~4. If any of these indices were
$\bot$ without all of them being such, then this would contradict
the purity assumption. We proceed in such a manner for all
$\b{\xi}{\str}$-terms and $\b{\xi}{\rts}$-terms.

Suppose $\beta$ is $c^k$. Then we have the following cases. If
only one of its indices is $\top$ or $\bot$, while the remaining
three indices are constant-free, then this contradicts the purity
assumption. If two of its indices are $\top$ or $\bot$, while the
remaining two indices are constant-free, then we may either apply
the equations $(c^k\delta\sigma)$ to eliminate $c^k$, or we would
contradict again the purity assumption. We have an analogous
situation when three indices of $c^k$ are $\top$ or $\bot$, while
the remaining one is constant-free. The last case is when all the
indices of $c^k$ are $\top$ or $\bot$. Then we apply either the
equations $(c^k\delta\sigma)$ to eliminate $c^k$, or the equations
$(c^k\kappa)$ to reduce $c^k$ to $\kappa$. If we apply
$(c^k\kappa)$, then we obtain a new occurrence of $\kappa$, with
which we deal as in the last part of the proof, into which we go
now.

Suppose we have a $\kappa$-term. By the purity assumption, this
term cannot be just $\kappa$. So it has a subterm of the form
${\mj_E\ks \kappa}$ or ${\kappa\ks \mj_E}$. We may assume as above
that $E$ is $\nu(D)$ for some $D$. No letter can occur in $E$,
because this would contradict the purity assumption. So $E$ is
$\top$ or $\bot$. If
${(\!\ks\!,E)\in\{(\kon,\bot),(\vee,\top)\}}$, then we apply the
equations $(\kappa\delta\sigma)$ to eliminate $\kappa$. If
${(\!\ks\!,E)\in\{(\kon,\top),(\vee,\bot)\}}$, then we apply the
equations mentioned in the proof of $\mk_{\top,\bot}^\emptyset$
Coherence in Section~6, which are consequences of naturality
equations. These equations do not eliminate $\kappa$, but they
replace a $\kappa$-term of greater complexity by a $\kappa$-term
of lesser complexity. This enables us to proceed by induction.
\qed

We prove now the main result of this section.

\prop{Restricted Bimonoidal Intermuting Coherence}{If
${f,g\!:A\str B}$ are arrows of $\mbox{\bf AC}^{\bf
k}_{\top,\bot}$ such that both $A$ and $B$ are pure, or no letter
occurs in them, then ${f=g}$ in $\mbox{\bf AC}^{\bf
k}_{\top,\bot}$.}

\dkz When no letter occurs in $A$ and $B$, we apply $\mbox{\bf
AC}^{{\bf k}\emptyset}_{\top,\bot}$ Coherence. When $A$ and $B$
are pure, we proceed as follows. We replace first $f$ and $g$ by
$f'$ and $g'$ respectively, which satisfy the conditions specified
in the Constant-Free Indices Lemma above; we have ${f=f'}$ and
${g=g'}$ in $\mbox{\bf AC}^{\bf k}_{\top,\bot}$.

As we applied Biassociative Coherence to obtain $\ACk^{st}$ in
Section 10 (and as we did previously in the proof of Normal
Bimonoidal Coherence in Section~7), we can apply Normal Biunital
Coherence of Section~5 to obtain a strictified category $\mbox{\bf
AC}^{{\bf k}^{\mbox{\scriptsize \emph{st}}}}_{\top,\bot}$
equivalent to $\mbox{\bf AC}^{\bf k}_{\top,\bot}$, where the
$\delta$-$\sigma$-arrows (see Section~2) and the $w$-arrows (see
Section~5) are identity arrows. So the arrows $f'$ and $g'$ are
mapped by the functor underlying this equivalence of categories
into the arrows ${f'',g''\!:\nu(A)\str\nu(B)}$ of $\mbox{\bf
AC}^{{\bf k}^{\mbox{\scriptsize \emph{st}}}}_{\top,\bot}$.

The arrows $f''$ and $g''$ may be represented by two arrow terms
of \ACk, and, by Biassociative Intermuting Coherence of the
preceding section, we may conclude that ${f''=g''}$ in \ACk, and
hence also in $\mbox{\bf AC}^{{\bf k}^{\mbox{\scriptsize
\emph{st}}}}_{\top,\bot}$. So, by equivalence of categories,
${f'=g'}$, and hence also ${f=g}$ in $\mbox{\bf AC}^{\bf
k}_{\top,\bot}$. \qed

The restriction imposed by this restricted coherence result is of
the same kind as the restriction Kelly and Mac Lane had for their
restricted coherence result for symmetric monoidal closed
categories (see \cite{KML71} and \cite{DP05}, end of Section 3.1,
or \cite{DP06}, Section~8).

\section{\large\bf Bimonoidal intermuting categories and two-fold loop spaces}
In this section we consider the relationship of our results to the
paper \cite{BFSV}, mentioned in the introduction. We will just
summarize matters and will not go into details, known either from
\cite{BFSV} or other references.

Restricted Bimonoidal Intermuting Coherence, which we have proved
in the preceding section, enables us to strengthen up to a point
Theorem 2.1 of \cite{BFSV}, one of the main results of that paper.
This theorem is about the notion of $n$-fold monoidal category,
which for ${n=2}$ is a bimonoidal intermuting category in which
the $b$-arrows, the $\delta$-$\sigma$-arrows, the $w$-arrows and
the arrow $\kappa$ are identity arrows. In two-fold monoidal
categories $\top$ and $\bot$ coincide. The theorem also relies on
the notion of lax functor from a category to the category
\emph{Cat} of categories (see \cite{S72}; the notion originates in
\cite{B67}) and in the topologists' notion of simplicial category,
called $\Delta$ (which in \cite{ML71}, Section VII.5, is called
$\Delta^+$). This theorem is formulated as follows.

\prop{Theorem 2.1 \cite{BFSV}}{An $n$-fold monoidal category $\cal
C$ determines a lax functor ${{\cal
C}_{\ast\ast\ldots\ast}\!:\Delta^{op}\times
\Delta^{op}\times\ldots \Delta^{op}\str Cat}$ such that ${{\cal
C}_{p_1,p_2,\ldots,p_n}={\cal C}^{p_1p_2\ldots p_n}}$.}

From this theorem one obtains immediately the main result of
\cite{BFSV}, which says that the group completion of the nerve of
an $n$-fold monoidal category is an $n$-fold loop space. We can
however prove the following theorem for \emph{strict} bimonoidal
intermuting categories, in which the $b$-arrows and the
$\delta$-$\sigma$-arrows are identity arrows.

\prop{Theorem}{A strict bimonoidal intermuting category $\cal C$
determines a lax functor ${{\cal C}_{\ast\ast}\!:\Delta^{op}\times
\Delta^{op}\str Cat}$ such that ${{\cal C}_{p_1,p_2}={\cal
C}^{p_1p_2}}$.}

\dkz We proceed as in the proof of Theorem 2.1 in \cite{BFSV}, and
rely on the notions introduced in that paper. For
${\gamma=(\alpha,\beta)}$ a pair of arrows of $\Delta^{op}$ such
that ${\alpha\!:m\str m'}$ and ${\beta\!:n\str n'}$, we denote by
$\gamma^\ast$ the functor ${{\cal C}_{\alpha,\beta}\!:{\cal
C}^{m\cdot n}\str{\cal C}^{m'\cdot n'}}$ defined as in
\cite{BFSV}. For ${}^\ast$ one can easily prove the following.

For ${k\geq 1}$ and ${i\in\{1,\ldots,k\}}$, let
${\gamma_i=(\alpha_i,\beta_i)}$ be a pair of arrows of
$\Delta^{op}$ such that ${\alpha_i\!:m_{i-1}\str m_i}$ and
${\beta_i\!:n_{i-1}\str n_i}$. Then
${\gamma_k^\ast\cirk\ldots\cirk\gamma_1^\ast}$ is a functor from
${\cal C}^{m_0\cdot n_0}$ to ${\cal C}^{m_k\cdot n_k}$.

For $\vec{p}$ being $p_1^1,\ldots,p_{n_0}^1,
\ldots,p_1^{m_0},\ldots,p_{n_0}^{m_0}$, and
\[
\gamma_k^\ast\cirk\ldots\cirk\gamma_1^\ast(\vec{p\,})=(A_1^1,\ldots,A_{n_k}^1,\ldots,
A_1^{m_k},\ldots,A_{n_k}^{m_k}),
\]
where we treat every object $p_j^i$ of $\cal C$ as a letter, we
have that if ${\nu(A_j^i)=\top}$, then for every ${l\in\{1,\ldots,
m_k\}}$ no letter occurs in $A_j^l$, and for every
${l\in\{1,\ldots,m_k\}}$ the elements of
${\{\nu(A_1^l),\ldots,\nu(A_{n_k}^l)\}}$ different from $\top$ are
all $\bot$ or all constant-free. As a corollary we have that every
$A_j^i$ is either pure or no letter occurs in it.

Hence, by our Restricted Bimonoidal Intermuting Coherence of the
preceding section, we conclude that any two ${m_3\cdot
n_3}$-tuples of ``canonical'' arrows of $\cal C$ from
${\gamma_3^\ast\cirk\gamma_2^\ast\cirk\gamma_1^\ast (\vec{p\,})}$
to ${(\gamma_3\cirk\gamma_2\cirk\gamma_1)^\ast(\vec{p\,})}$ are
equal, which is sufficient to show that ${\cal C}_{\ast\ast}$ is a
lax functor. \qed

\noindent We leave open the question whether this theorem could be
generalized to cover a notion of $n$-fold intermuting category for
$n\geq 2$, which generalizes the notion of $n$-fold monoidal
category as our notion of bimonoidal intermuting category
generalizes the notion of two-fold monoidal category. As we said
in the introduction, we suppose that this can be achieved by
relying on the technique of the next section.

Our strengthening of Theorem 2.1 of \cite{BFSV} should be
considered in the light of the statement of the main result of
\cite{BFSV}, made in the first sentence of the abstract and in the
last paragraph of the Introduction. This is the statement that
\emph{for all} $n$ the notion of $n$-fold monoidal category
corresponds \emph{precisely} to the notion of $n$-fold loop space.
There are interesting categories that are not two-fold monoidal,
but which, according to our strengthening, give rise to a two-fold
loop space. The dicartesian categories of \cite{DP04} (Section
9.6), which are categories with all finite products and coproducts
(including the empty ones) and with $w$-isomorphisms added, are
such. A concrete dicartesian category obtained by extending the
category of pointed sets with the empty set is described in
\cite{DP04} (Section 9.7).

Bimonoidal intermuting categories do not make like two-fold
monoidal categories assumptions that preclude a logical
interpretation. Such is the assumption that the two unit objects,
$\top$ and $\bot$, are isomorphic, which delivers the arrows of
the type ${A\vee B\str A\kon B}$. The assumptions of dicartesian
categories can be naturally interpreted in conjunctive-disjunctive
logic, including the constants $\top$ and $\bot$.

Dicartesian categories are lattice categories, in the sense of
\cite{DP04} (Section 9.4). We will examine lattice categories in
Section 15 below, and we will show that every such category has as
its fragment a symmetric biassociative intermuting structure,
which we define in the next section.

Our approach differs from \cite{BFSV} also in the following
respect. Our Biassociative Intermuting Coherence of Section 11
states that \ACk\ is a preorder. For the category of \cite{BFSV}
closest to \ACk\ one cannot obtain this result, but only that it
is what we call in the next section a \emph{diversified} preorder,
which is a result like the Theorem of Section 11, where a
diversification proviso is involved (see Section~10).

\section{\large\bf Symmetric biassociative intermuting categories}
In this section we add symmetry, i.e.\ natural commutativity
isomorphisms for $\kon$ and $\vee$, to the biassociative
intermuting categories of Section 10, and prove an appropriate
coherence result for the ensuing notion.

A \emph{symmetric biassociative} category is a biassociative
category ${\langle{\cal A},\kon,\vee\rangle}$ (see Section 2) in
which there are two natural self-inverse isomorphisms $\c{\xi}$,
for $\!\ks\!\in\{\kon,\vee\}$, with the following components in
\aA:
\[
\c{\xi}_{A,B}:A\ks B\str B\ks A,
\]
which satisfy Mac Lane's hexagonal equations:
\[
\c{\xi}_{A,B\kst
C}=\:\b{\xi}{\str}_{B,C,A}\!\!\cirk(\mj_B\ks\!\c{\xi}_{A,C})
\cirk\b{\xi}{\rts}_{B,A,C}\!\!\cirk(\c{\xi}_{A,B}\!\ks\mj_C)\cirk
\b{\xi}{\str}_{A,B,C}.
\]
We call the arrows $\hat{c}$ and $\check{c}$ collectively
\emph{$c$-arrows}. Let \ms\ be the free symmetric biassociative
category generated by a set of objects.

We say that a formula is \emph{diversified} when every letter
occurs in it at most once. (We already had diversified form
sequences in Section 10). We say that a category whose objects are
formulae is a \emph{diversified preorder} when for every pair of
arrows ${f,g\!:A\str B}$ such that $A$ and $B$ are diversified we
have ${f=g}$ in this category. In \cite{DP04} (Section 6.3) one
can find a proof, based on Mac Lane's symmetric monoidal
coherence, that \ms\ is a diversified preorder. We call this fact
Symmetric Biassociative Coherence. In \cite{DP04} we call by the
same name an equivalent result about the existence of a faithful
functor from the category \ms\ to the category whose objects are
finite ordinals and whose arrows are bijections.

A \emph{symmetric biassociative intermuting} category is a
symmetric biassociative and biassociative intermuting category
(see Section 10) that satisfies moreover the equations
\begin{tabbing}
\hspace{.5em}\=$(\psi\check{c})$\quad\=
$c^k_{A_2,A'_2,A_1,A'_1}\cirk\check{c}_{A_1\kon A'_1,A_2\kon
A'_2}$ \= = \=
$(\check{c}_{A_1,A_2}\kon\check{c}_{A'_1,A'_2})\cirk
c^k_{A_1,A'_1,A_2,A'_2}$,
\\[2ex]
\>$(\overline{\psi}\hat{c})$\> $\hat{c}_{A_1\vee A'_1,A_2\vee
A'_2}\cirk c^k_{A_1,A_2,A'_1,A'_2}$ \> = \>
$c^k_{A_2,A_1,A'_2,A'_1}\cirk(\hat{c}_{A_1,A_2}\vee\hat{c}_{A'_1,A'_2})$,
\end{tabbing}
which we call collectively $(c^kc)$ (these equations may be found
in \cite{L06}, Section 2.3). The equation $(\psi\check{c})$ says
that $\check{c}$ is upward preserved by $\kon$, while
$(\overline{\psi}\hat{c})$ says that $\hat{c}$ is downward
preserved by $\vee$ (see Section~4).

In every symmetric biassociative intermuting category the
following equations hold:
\begin{tabbing}
\hspace{.5em}\= for ${\!\ks\!\in\{\kon,\vee\}}$ and\\[.5ex]
\>${\c{\xi}{\!^m}{\!\!_{A,B,C,D}}}=_{df}\:\b{\xi}{\str}_{A,C,B\kst
D}\!\!\cirk(\mj_A\ks(\b{\xi}{\rts}_{C,B,D}\!\!\cirk
(\c{\xi}_{B,C}\ks\mj_D)\cirk\b{\xi}{\str}_{B,C,D}))\cirk
\b{\xi}{\rts}_{A,B,C\kst D}:$
\\*[.5ex]
\`$(A\ks B)\ks(C\ks D)\str(A\ks C)\ks(B\ks D)$,
\\[2ex]
\> $(\psi\check{c}^m)$\quad\= $c^k_{A_1\vee A_3,A'_1\vee
A'_3,A_2\vee A_4,A'_2\vee A'_4}\cirk(c^k_{A_1,A'_1,A_3,A'_3}\vee
c^k_{A_2,A'_2,A_4,A'_4})\cirk$
\\*[.5ex]
\`$\cirk\check{c}^m_{A_1\kon A'_1,A_2\kon A'_2,A_3\kon
A'_3,A_4\kon A'_4}=$
\\*[.5ex]
\>
\>$(\check{c}^m_{A_1,A_2,A_3,A_4}\kon\check{c}^m_{A'_1,A'_2,A'_3,A'_4})\cirk
c^k_{A_1\vee A_2,A'_1\vee A'_2,A_3\vee A_4,A'_3\vee A'_4}\cirk$
\\*[.5ex]
\`$\cirk(c^k_{A_1,A'_1,A_2,A'_2}\vee c^k_{A_3,A'_3,A_4,A'_4})$,
\\[2ex]
\>$(\overline{\psi}\hat{c}^m)$\> $\hat{c}^m_{A_1\vee A_3,A'_1\vee
A'_3,A_2\vee A_4,A'_2\vee A'_4}\cirk(c^k_{A_1,A'_1,A_3,A'_3}\kon
c^k_{A_2,A'_2,A_4,A'_4})\cirk$
\\*[.5ex]
\`$\cirk c^k_{A_1\kon A'_1,A_2\kon A'_2,A_3\kon A'_3,A_4\kon
A'_4}=$
\\*[.5ex]
\> \>$(c^k_{A_1,A_2,A_3,A_4}\kon c^k_{A'_1,A'_2,A'_3,A'_4})\cirk
c^k_{A_1\kon A_2,A'_1\kon A'_2,A_3\kon A_4,A'_3\kon A'_4}\cirk$
\\*[.5ex]
\`$\cirk(\hat{c}^m_{A_1,A'_1,A_2,A'_2}\vee\hat{c}^m_{A_3,A'_3,A_4,A'_4})$.
\end{tabbing}

These two equations, which we call collectively ${(c^kc^m)}$, are
both obtained from the equation ${(\psi c^m)}$ of \cite{DP06a}
(Section~2) by replacing some occurrences of $\kon$ by $\vee$.
(The equation ${(\psi c^m)}$ is the equation
${(\overline{\psi}\hat{c}^m)}$ with $c^k$ and $\vee$ replaced
everywhere by $\hat{c}^m$ and $\kon$ respectively.) The equation
${(\psi \check{c}^m)}$ says that $\check{c}^m$ is upward preserved
by $\kon$, and ${(\overline{\psi}\hat{c}^m)}$ says that $\hat{c}$
is downward preserved by $\vee$. The equations $(c^kc^m)$ are
derived from the equations assumed for symmetric biassociative
intermuting categories by relying on the fact that the natural
transformations used to define $\check{c}^m$ are upward preserved
by $\kon$, and those used to define $\hat{c}^m$ are downward
preserved by $\vee$ (cf.\ the beginning of Section 10).

The equations ${(c^kc^m)}$ are related to the equation
corresponding to the hexagonal interchange diagram used for
defining $n$-fold monoidal categories for $n>2$ in \cite{BFSV}
(end of Definition 1.7). Let us call this equation \emph{HI}.
Taking that $\check{c}^m$ and $\hat{c}^m$ are $\eta^{11}$ and
$\eta^{22}$ respectively, the equation $(\psi\check{c}^m)$ is
obtained by putting ${i=j=1}$ and ${k=2}$ in \emph{HI}, and the
equation $(\overline{\psi}\hat{c}^m)$ is obtained by putting
${i=1}$ and ${j=k=2}$ in \emph{HI}. The equation ${(\psi c^m)}$
mentioned above is obtained by putting ${i=j=k}$ in \emph{HI},
which amounts to omitting the superscripts of $\eta$.

Let \SCk\ be the free symmetric biassociative intermuting category
generated by a set of objects. We can proceed as in \cite{DP04}
(Section 7.6) to obtain a strictified category $\SCk^{st}$ where
the $b$-arrows and the $c$-arrows are identity arrows. The
category $\SCk^{st}$ is not equivalent to \SCk, but we have that
if $\SCk^{st}$ is a diversified preorder in a sense to be defined
below, then \SCk\ is a diversified preorder. The remainder of this
section is devoted to proving that $\SCk^{st}$ is indeed such a
diversified preorder. This proof could be adapted to prove the
Theorem of Section 11, which says that $\ACk^{st}$ is a
diversified preorder. (We have however preferred to rely on
rectangular grids in Section 11, because they give a clearer
picture.)

The objects of $\SCk^{st}$ are equivalence classes of formulae
$\lz A\dz$ such that $\lz A\dz$ is the set of all formulae
isomorphic to $A$ in the free symmetric biassociative category
$\mathbf{S}$. We may represent unequivocally the equivalence class
$\lz A\dz$ by an equivalence class of form sequences, such that
form sequences in the class differ from each other in the order of
conjuncts and disjuncts. We call such equivalence classes of form
sequences \emph{form multisets} (see \cite{DP04}, Section 7.7).
For example, ${p\kon q\kon(p\vee r\vee p)}$ and ${q\kon(r\vee
p\vee p)\kon p}$ stand for the same form multiset.

As we did in Section 10 we assume that we deal with
\emph{diversified} form multisets, i.e.\ with form multisets in
which every letter occurs at most once. Such form multisets are
called \emph{form sets}. That $\SCk^{st}$ is a \emph{diversified
preorder} means that if ${f',f''\!:X\str Y}$ are arrows of
$\SCk^{st}$ where $X$ and $Y$ are form sets, i.e.\ diversified
form multisets, then ${f'=f''}$.

We denote by ${let(X)}$ the set of letters occurring in the form
set $X$. Let $p$ be a letter such that $X$ is not $p$. We define
inductively the form set $X^{-p}$ as we did in Section 10 when $X$
was a form sequence, and we use for $X$ and $Y$ form sets and $P$
a set of letters the notation $X^{-P}$ and $X^{-Y}$ defined
analogously, with the same provisos, as in Section 10.

For $X$ a form set, let $X'$ be a \emph{subformset} of $X$ when
$X'$, which stands for a form set, is a subword of a
representative of $X$. We can then state the following.

\prop{Lemma 1}{If ${f\!: X\str Y}$ is an arrow of $\SCk^{st}$, and
$P$ is a set of letters such that for every subformset ${U\kon V}$
of $X$
\[ {let(U)}\subseteq
P\quad {\mbox {\it iff}} \quad{let(V)}\subseteq P,
\]
then this equivalence holds for every subformset ${U\kon V}$ of
$Y$.}

\dkz It is sufficient to prove this lemma for $f$ being a
$c^k$-term. We proceed by induction on the complexity of this
$c^k$-term. For the basis, suppose that $f$ is
\[
c^k_{X_1,X_2,X_3,X_4}\!:(X_1\kon X_2)\vee(X_3\kon X_4)\str(X_1\vee
X_3)\kon(X_2\vee X_4).
\]
\begin{tabbing}
We have\quad${let(X_1\vee X_3)\subseteq P}$\quad\= iff\quad\=
${let (X_1)\subseteq P}$ and ${let(X_3)\subseteq P}$,
\\[1ex]
\> iff\> ${let (X_2)\subseteq P}$ and ${let(X_4)\subseteq P}$, by
the
\\*[.5ex]
\` assumption concerning the source $X$ of $f$,
\\[1ex]
\> iff\> ${let(X_2\vee X_4)\subseteq P}$,
\end{tabbing}
and the main $\kon$ of ${(X_1\vee X_3)\kon(X_2\vee X_4)}$ is the
only ``new'' $\kon$ in the target of~$f$.

For the induction step, suppose $f$ is of the form
${f'\kon\mj_{X_2}}$ for ${f'\!:X_1\str Z}$. If $Z$ is ${Z_1\vee
Z_2}$, then we just apply the induction hypothesis to $f'$.

Suppose $Z$ is ${Z_1\kon Z_2}$. If ${let(X_2)\subseteq P}$, then
${let(X_1)\subseteq P}$, which implies that ${let(Z_1\kon
Z_2)}\subseteq P$, and hence ${let(Z_i)\subseteq P}$ for every
${i\in\{1,2\}}$. If ${let(Z_i)\subseteq P}$ for some
${i\in\{1,2\}}$, then by the induction hypothesis
${let(Z_{3-1})\subseteq P}$, and so ${let(X_1)\subseteq P}$. Hence
${let(X_2)\subseteq P}$.

If $f$ is of the form ${f'\vee\mj_{X_2}}$, then we just apply the
induction hypothesis to~$f'$. \qed

As a corollary we have the following lemma.

\prop{Lemma 2}{If ${f\!:X_1\vee X_2\str Y}$ is an arrow of
$\SCk^{st}$, then for every subformset ${Y_1\kon Y_2}$ of $Y$ and
for every ${i\in\{1,2\}}$
\[
let(Y_1)\subseteq let(X_i)\quad{\mbox{\it iff}}\quad
let(Y_2)\subseteq let(X_i).
\]}
\vspace{-3ex}

The arrow terms of \SCk\ where indices are replaced by form sets
are the arrow terms of $\SCk^{st}$. In these arrow terms of
$\SCk^{st}$ we may replace all the $b$-arrows and $c$-arrows by
identity arrows. Since in developed $\SCk^{st}$ arrow terms we
need to consider only $c^k$-terms, we easily establish the
following.

\prop{Lemma 3}{Every arrow term ${f\!:X_1\kon X_2\str Y}$ of
$\SCk^{st}$ is equal to ${f_1\kon f_2}$ for some ${f_i\!:X_i\str
Y_i}$ where ${i\in\{1,2\}}$.}

Since the equations assumed for \SCk\ are such that the number of
occurrences of $c^k$ is equal on the two sides of the equations,
we have in general the following.

\prop{Lemma 4}{If ${f=g}$ in \SCk\ or $\SCk^{st}$, then the number
of occurrences of $c^k$ is equal in the arrow terms $f$ and $g$.}

Let ${f\!:X\str Y}$ be an arrow term of $\SCk^{st}$, and let $P$
be a set of letters such that ${let(X)-P\neq\emptyset}$ (which
implies that ${let(Y)-P\neq\emptyset}$, since ${let(X)=let(Y)}$),
and such that, as in Lemma~1, for every subformset ${U\kon V}$ of
$X$ we have ${let(U)\subseteq P}$ iff ${let(V)\subseteq P}$. Then
we define inductively the arrow $f^{-P}\!:X^{-P}\str Y^{-P}$ of
$\SCk^{st}$ in the following manner:
\begin{itemize}
\item[] if ${let(X)\cap P=\emptyset}$, then ${f^{-P}=f}$; \item[]
if $f$ is $\mj_X$, then ${f^{-P}=\mj_{X^{-P}}}$; \item[] if $f$ is
$c^k_{X_1,X_2,X_3,X_4}$, then ${f^{-P}=\mj_{X^{-P}}}$ in case
either ${let(X_1\kon X_2)\subseteq P}$ or ${let(X_3\kon
X_4)\subseteq P}$, and
${f^{-P}=c^k_{X_1^{-P},X_2^{-P},X_3^{-P},X_4^{-P}}}$ otherwise;
\item[] if $f$ is ${f_1\ks f_2}$, for ${\!\ks\!\in\{\kon,\vee\}}$
and ${f_i\!:X_i\str Y_i}$, where ${i\in\{1,2\}}$, then
${f^{-P}=f^{-P}_{3-i}}$ in case ${let(X_i)\subseteq P}$, and
${f^{-P}=f^{-P}_1\ks f^{-P}_2}$ in case for no ${i\in\{1,2\}}$ the
inclusion ${let(X_i)\subseteq P}$ holds; \item[] if $f$ is
${f_2\cirk f_1}$, then by Lemma~1, both $f_1^{-P}$ and $f_2^{-P}$
are defined, and ${f^{-P}=f_2^{-P}\cirk f_1^{-P}}$.
\end{itemize}

For $X_1$ and $X_2$ form sets, we say that $c^k_{S,T,U,V}$ is
${(X_1,X_2)}$-\emph{splitting} when ${let(S\kon T)}\subseteq {let
(X_i)}$ and ${let(U\kon V)}\subseteq {let(X_{3-i})}$ for some
${i\in\{1,2\}}$. We say that an arrow term of $\SCk^{st}$ is
${(X_1,X_2)}$-\emph{splitting} when every occurrence of $c^k$ in
it is ${(X_1,X_2)}$-splitting, and we say that it is
${(X_1,X_2)}$-\emph{nonsplitting} when every occurrence of $c^k$
in it is not ${(X_1,X_2)}$-splitting. One can easily check that if
${f=g}$ and $f$ is an ${(X_1,X_2)}$-splitting arrow term, then $g$
is an ${(X_1,X_2)}$-splitting arrow term too. (This is not the
case when we replace ``splitting'' by ``nonsplitting''.) It is
clear that every ${(X_1,X_2)}$-splitting arrow term is equal to a
developed ${(X_1,X_2)}$-splitting arrow term, and analogously with
``splitting'' replaced by ``nonsplitting''. We can then prove the
following.

\prop{Lemma 5}{If ${f\!:X_1\vee X_2\str Y}$ is
${(X_1,X_2)}$-nonsplitting, then $f$ is equal to ${f_1\vee f_2}$
for ${f_i\!:X_i\str Y_i}$ and ${i\in\{1,2\}}$.}

\dkz It is enough to establish this lemma for $f$ being a
$c^k$-term, which is trivial. We then proceed by induction on the
number of factors in a developed arrow term equal to $f$. \qed

Lemmata 6-14 below will yield a normal form for
${(X_1,X_2)}$-splitting arrow terms of $\SCk^{st}$ with source
${X_1\vee X_2}$. After these lemmata we establish with the help of
this normal form and of Lemma~5 that every arrow term
${f\!:X_1\vee X_2\str Y}$ of $\SCk^{st}$ is equal to
${h\cirk(g_1\vee g_2)}$ for $h$ being ${(X_1,X_2)}$-splitting and
for ${g_i\!:X_i\str Y^{-X_{3-i}}}$. After that we will be able to
proceed by induction to show that $\SCk^{st}$ is a diversified
preorder. First we prove the following.

\prop{Lemma 6}{If ${f\!:X_1\vee X_2\str Y}$ is
${(X_1,X_2)}$-splitting, then $Y^{-X_i}$ is $X_{3-i}$ for every
${i\in\{1,2\}}$.}

\dkz Since ${f^{-X_i}\!:(X_1\vee X_2)^{-X_i}\str Y^{-X_i}}$ is
$\mj_{Y^{-X_i}}$, and since ${(X_1\vee X_2)^{-X_i}}$ is $X_{3-i}$,
the lemma follows. \qed

\prop{Lemma 7}{If ${f\!:X_1\vee X_2\str Y_1\kon Y_2}$ is
${(X_1,X_2)}$-splitting, then $X_i$ is of the form ${X_i'\kon
X_i''}$ for every ${i\in\{1,2\}}$.}

\dkz We have that
\[
f^{-X_{3-i}}=\mj_{X_i}\!:X_i\str(Y_1\kon Y_2)^{-X_{3-i}}.
\]
By Lemma~2, we have that ${(Y_1\kon Y_2)^{-X_{3-i}}}$ is
${Y_1^{-X_{3-i}}\kon Y_2^{-X_{3-i}}}$, which is $X_i$. \qed

For $X_1$ being ${X_1'\kon X_1''}$ and $X_2$ being ${X_2'\kon
X_2''}$, consider an arrow term
\[
f\cirk c^k_{X_1',X_1'',X_2',X_2''}\!: X_1\vee X_2\str Y_1\kon Y_2
\]
of $\SCk^{st}$. For ${i,j\in\{1,2\}}$ and ${l\in\{',''\}}$ we say
that $X_j^l$ is \emph{exactly} in $Y_i$ when $X_j^l$ is
$Y_i^{-X_{3-j}}$. We say that $X_j^l$ is \emph{properly} in $Y_i$
when $X_j^l$ is $P_j$ for ${P_1\kon P_2}$ being $Y_i^{-X_{3-j}}$.
We say that $X_j^l$ \emph{extends} $Y_i$ when $X_j^l$ is
${Y_i^{-X_{3-j}}\kon Q_j}$ for ${Q_1\kon Q_2}$ being
$Y_{3-i}^{-X_{3-j}}$. Finally, we say that $X_j^l$ \emph{partakes}
in $Y_1$ and $Y_2$ when $X_j^l$ is ${P_j\kon Q_j}$ for ${P_1\kon
P_2}$ being $Y_1^{-X_{3-j}}$ and ${Q_1\kon Q_2}$ being
$Y_2^{-X_{3-j}}$. We say that $X_1^l$ and $X_2^l$ are
\emph{related in the same way} to $Y_i$ when either they are both
exactly in $Y_i$, or they are both properly in $Y_i$, or they both
extend $Y_i$, or they both partake in $Y_1$ and~$Y_2$.

For Lemmata 8-10 below we assume that ${f\cirk
c^k_{X_1',X_1'',X_2',X_2''}}$, of the type displayed above, is
${(X_1,X_2)}$-splitting.

\prop{Lemma 8}{For every $i$, $j$ and $l$, we have that $X_j^l$ is
either exactly in $Y_1$ or $Y_2$, or $X_j^l$ is properly in $Y_i$,
or $X_j^l$ extends $Y_i$, or $X_j^l$ partakes in $Y_1$ and $Y_2$.}

\dkz By Lemma~6 and Lemma~2 we have that
\[
X_j'\kon X_j''\quad{\mbox{\rm is}}\quad Y_1^{-X_{3-j}}\kon
Y_2^{-X_{3-j}}.
\]
From this the lemma follows. \qed

\vspace{-2ex}

\prop{Lemma 9}{For every $i$, $j$ and $l$, we have that $X_j^l$ is
properly in $Y_i$ iff $X_{3-j}^l$ extends $Y_{3-i}$, and $X_j^l$
partakes in $Y_1$ and $Y_2$ iff $X_{3-j}^l$ partakes in $Y_1$ and
$Y_2$.}

\noindent The proof of Lemma~8 yields this lemma too.

\prop{Lemma 10}{For every $i$, $j$, and $l$, we have that $X_1^l$
and $X_2^l$ are related in the same way to $Y_i$.}

\dkz Suppose $X_1^l$ and $X_2^l$ are not related in the same way
to $Y_i$. Then, by Lemma~8, we have to consider a number of cases,
and show that we have a contradiction in each of these cases. We
will consider just one, typical, case. In all the other cases we
proceed more or less in the same way.

Suppose $X_1'$ is exactly in $Y_1$ and $X_2'$ is properly in
$Y_1$. Then, by Lemma~9, we have that $X_1''$ is exactly in $Y_2$
and $X_2''$ extends $Y_2$, which means that $X_1''$ is
$Y_2^{-X_2}$ and $X_2''$ is ${Y_2^{-X_1}\kon Q_2}$ for ${Q_1\kon
Q_2}$ being $Y_1^{-X_1}$. By Lemma~3, we have that
\[
f=f'\kon f''\!:(X_1'\vee X_2')\kon(X_1''\vee X_2'')\str Z'\kon Z''
\]
for ${f^l\!:X_1^l\vee X_2^l\str Z^l}$ where ${l\in\{',''\}}$, and
${Z'\kon Z''}$ the same form set as ${Y_1\kon Y_2}$. So $f''$ is
of the type ${Y_2^{-X_2}\vee(Y_2^{-X_1}\kon Q_2)\str Z''}$ and
${let(Z'')}={let(Y_2)}\cup{let(Q_2)}$.

Since ${let(Q_2)}\subseteq {let(Y_1)}$ and ${Z'\kon Z''}$  is the
same form set as ${Y_1\kon Y_2}$, the form set $Z''$ must be
${Z_1''\kon Z_2''}$ for ${let(Z_2'')}={let(Q_2)}$. Hence ${f\cirk
c^k_{X_1',X_1'',X_2',X_2''}}$ is of the type ${X_1\vee
X_2\str(Z'\kon Z_1'')\kon Z_2''}$. We have
${let(Z_2'')}\subseteq{let(X_2)}$, because ${let(Q_2)}\subseteq
{let(Y_1^{-X_1})}$. So, by Lemma~2, ${let(Z'\kon Z_1'')}\subseteq
{let(X_2)}$. But this would mean that ${let(X_1)=\emptyset}$,
which is a contradiction. \qed

We can also prove the following.

\prop{Lemma 11}{If ${f\!:X_1\vee X_2\str Y_1\kon Y_2}$ is
${(X_1,X_2)}$-splitting, then for some ${(X_1,X_2)}$-splitting
arrow term $g$ we have ${f=g\cirk
c^k_{Y_1^{-X_2},Y_2^{-X_2},Y_1^{-X_1},Y_2^{-X_1}}}$.}

\dkz Suppose ${f\!:X_1\vee X_2\str Y_1\kon Y_2}$ is
${(X_1,X_2)}$-splitting. By the remark made after the definition
of ${(X_1,X_2)}$-splitting arrow terms, and by Lemma~7, we have
that ${f=g\cirk c^k_{X_1',X_1'',X_2',X_2''}}$ for $g$ an
${(X_1,X_2)}$-splitting developed arrow term, where $X_j$ is
${X_j'\kon X_j''}$ for ${j\in\{1,2\}}$. We proceed by induction on
the number $n$ of $c^k$-terms in $g$. If ${n=0}$, then
${g=\mj_{(X_1'\vee X_2')\kon(X_1''\vee X_2'')}}$ for $Y_1$ being
${X_1'\vee X_2'}$ and $Y_2$ being ${X_1''\vee X_2''}$. So
${f=c^k_{Y_1^{-X_2},Y_2^{-X_2},Y_1^{-X_1},Y_2^{-X_1}}}$.

If ${n>0}$, and ${c^k_{X_1',X_1'',X_2',X_2''}=
c^k_{Y_1^{-X_2},Y_2^{-X_2},Y_1^{-X_1},Y_2^{-X_1}}}$, then we are
done.

If this equation does not hold, then, by Lemmata 10 and~9, we have
two cases to consider:

\begin{tabbing}
\hspace{2em}\= (I)\quad\= $X_1^l$ and $X_2^l$ both extend $Y_i$
for some ${i\in\{1,2\}}$,
\\[1ex]
\>(II)\> $X_1^l$ and $X_2^l$ both partake in $Y_1$ and $Y_2$.
\end{tabbing}
We first consider case (I), and assume that $l$ is $'$ and $i$ is
1. So $X_1'$ is ${Y_1^{-X_2}\kon Q}$ for ${Q\kon X_1''}$ being
$Y_2^{-X_2}$, and $X_2'$ is ${Y_1^{-X_1}\kon R}$ for ${R\kon
X_2''}$ being $Y_2^{-X_1}$. We first conclude by Lemma~3 that
${g=g'\kon g''}$ for ${g'\!:X_1'\vee X_2'\str Y_1\kon Y_2'}$ and
${g''\!:X_1''\vee X_2''\str Y_2''}$, where ${Y_2'\kon Y_2''}$ is
$Y_2$, because ${let(Y_1)}$ is a proper subset of ${let(X_1'\vee
X_2')}$. Since $g$ is ${(X_1,X_2)}$-splitting, $g'$ is
${(X_1,X_2)}$-splitting, and hence also ${(X_1',X_2')}$-splitting.
So, by the induction hypothesis, for an ${(X_1',X_2')}$-splitting
arrow term $h$ we have that
\[
g'=h\cirk
c^k_{Y_1^{-X_2'},{Y_2'}^{-X_2'},Y_1^{-X_1'},{Y_2'}^{-X_1'}}.
\]
By Lemma~6 applied to $g'$, we have that ${(Y_1\kon
Y_2')^{-X_i'}}$ is $X_{3-i}'$. So ${Y_1^{-X_2'}\kon Y_2'^{-X_2'}}$
is ${Y_1^{-X_2}\kon Q}$ and ${Y_1^{-X_1'}\kon Y_2'^{-X_1'}}$ is
${Y_1^{-X_1}\kon R}$. Since ${let(Y_1^{-X_i'})}\subseteq
{let(Y_1)}$ and ${let(Q_1)}\cup {let(R_1)}\subseteq {let(Y_2)}$,
we conclude that
\[
Y_1^{-X_2'}\;{\mbox{\rm is }}Y_1^{-X_2},\quad
Y_2'^{-X_2'}\;{\mbox{\rm is }}Q,\quad Y_1^{-X_1'}\;{\mbox{\rm is
}}Y_1^{-X_1}\quad{\mbox{\rm and }}\quad Y_2'^{-X_1'}\;{\mbox{\rm
is }}R.
\]
So
\begin{tabbing}
\hspace{2.5em}$f$ \= $=(h\kon
g'')\cirk(c^k_{Y_1^{-X_2},Q,Y_1^{-X_1},R}\kon\mj_{X_1''\vee
X_2''})\cirk c^k_{Y_1^{-X_2}\kon Q,X_1'',Y_1^{-X_1}\kon R,X_2''}$
\\*[1ex]
\> $=(h\kon g'')\cirk(\mj_{Y_1^{-X_2}\vee Y_1^{-X_1}}\kon
c^k_{Q,X_1'',R,X_2''})\cirk c^k_{Y_1^{-X_2},Q\kon
X_1'',Y_1^{-X_1},R\kon X_2''}$
\end{tabbing}
by the equation $(\overline{\psi}\hat{b})$ of Section~4. Here
${Q\kon X_1''}$ is $Y_2^{-X_2}$ and ${R\kon X_2''}$ is
$Y_2^{-X_1}$, as desired. We proceed analogously when $l$ is $'$
and $i$ is 2, and when $l$ is $''$ and $i$ is 1 or 2.

Next we consider case (II). So $X_1'$ is ${P'\kon Q'}$ for
${P'\kon P''}$ being $Y_1^{-X_2}$ and ${Q'\kon Q''}$ being
$Y_2^{-X_2}$, while $X_2'$ is ${O'\kon R'}$ for ${O'\kon O''}$
being $Y_1^{-X_1}$ and ${R'\kon R''}$ being $Y_2^{-X_1}$. We have
that $X_1''$ is ${P''\kon Q''}$ and $X_2''$ is ${O''\kon R''}$.

We conclude by Lemma~3 that ${g=g'\kon g''}$ for ${g^l\!:X_1^l\vee
X_2^l\str Y_1^l\kon Y_2^l}$, where ${l\in\{',''\}}$ and ${Y_i'\kon
Y_i''}$ is $Y_i$ for ${i\in\{1,2\}}$. Since $g$ is
${(X_1,X_2)}$-splitting, $g^l$ for ${l\in\{',''\}}$ is
${(X_1,X_2)}$-splitting, and hence also
${(X_1^l,X_2^l)}$-splitting. So, by the induction hypothesis, for
an ${(X_1^l,X_2^l)}$-splitting arrow terms $h^l$ we have that
\[
g^l=h^l\cirk
c^k_{{Y_1^l}^{-X_2^l},{Y_2^l}^{-X_2^l},{Y_1^l}^{-X_1^l},{Y_2^l}^{-X_1^l}}.
\]
By Lemma~6 applied to $g^l$ we have that ${(Y_1^l\kon
Y_2^l)^{-X_i^l}}$ is $X_{3-i}^l$. So ${{Y_1'}^{-X_2'}\kon
{Y_2'}^{-X_2'}}$ is ${P'\kon Q'}$. Since ${let(P')}\subseteq
{let(Y_1)}$ and ${let(Q')}\subseteq {let(Y_2)}$, we conclude that
\begin{tabbing}
\hspace{10em}\= $Y_1''^{-X_2''}$ \= is \= $P''$\quad\= and\quad\=
$Y_2''^{-X_2''}$ \= is \= $Q''$\kill

\> $Y_1'^{-X_2'}$ \> is \> $P'$ \> and \> $Y_2'^{-X_2'}$ \> is \>
$Q'$.
\\[1ex]
We conclude similarly that
\\[1ex]
\> $Y_1'^{-X_1'}$ \> is \> $O'$ \> and \> $Y_2'^{-X_1'}$ \> is \>
$R'$,
\\[1ex]
\> $Y_1''^{-X_2''}$ \> is \> $P''$ \> and \> $Y_2''^{-X_2''}$ \>
is \> $Q''$,
\\[1ex]
\> $Y_1''^{-X_1''}$ \> is \> $O''$ \> and \> $Y_2''^{-X_1''}$ \>
is \> $R''$.
\end{tabbing}
So
\begin{tabbing}
\hspace{1em}$f$ \= $= (h'\kon h'')\cirk(c^k_{P',Q',O',R'}\kon
c^k_{P'',Q'',O'',R''})\cirk c^k_{P'\kon Q',P''\kon Q'',O'\kon
R',O''\kon R''}$
\\*[1ex]
\> $= (h'\kon h'')\cirk(c^k_{P',P'',O',O''}\kon
c^k_{Q',Q'',R',R''})\cirk c^k_{P'\kon P'',Q'\kon Q'',O'\kon
O'',R'\kon R''}$,
\end{tabbing}
by the equation $(\overline{\psi}\hat{c}^m)$ from the beginning of
this section. Here ${P'\kon P''}$ is $Y_1^{-X_2}$, ${Q'\kon Q''}$
is $Y_2^{-X_2}$, ${O'\kon O''}$ is $Y_1^{-X_1}$ and ${R'\kon R''}$
is $Y_2^{-X_1}$, as desired. \qed

Next we prove the following.

\prop{Lemma 12}{If $f_1,f_2\!:(p_1\kon\ldots\kon
p_n)\vee(q_1\kon\ldots\kon q_n)\str(p_1\vee
q_1)\kon\ldots\kon(p_n\vee q_n)$,  for ${n\geq 1}$ and
$p_1,\ldots, p_n,q_1,\ldots, q_n$ distinct letters, are arrow
terms of $\SCk^{st}$, then ${f_1=f_2}$.}

\dkz If $n=1$, then it is clear that ${f_1=f_2=\mj_{p_1\vee
q_1}}$. If $n>1$, then it is easy to see that $f_1$ and $f_2$ are
${(p_1\kon\ldots\kon p_n,q_1\kon\ldots\kon q_n)}$-splitting. Then
by Lemmata 11 and~3, it follows that for every ${i\in\{1,2\}}$
\[
f_i=(\mj_{p_1\vee q_1}\kon f_i')\cirk c^k_{p_1,p_2\kon\ldots\kon
p_n,q_1,q_2\kon\ldots\kon q_n}
\]
for ${f_i'\!:(p_2\kon\ldots\kon p_n)\vee(q_2\kon\ldots\kon
q_n)\str(p_2\vee q_2)\kon\ldots\kon(p_n\vee q_n)}$. By the
induction hypothesis, ${f_1'=f_2'}$. \qed

Let $c^k_{p_1,\ldots,p_n,q_1,\ldots,q_n}$ stand for any arrow term
of $\SCk^{st}$ of the type
\[
(p_1\kon\ldots\kon p_n)\vee(q_1\kon\ldots\kon q_n)\str(p_1\vee
q_1)\kon\ldots\kon(p_n\vee q_n).
\]
In every such arrow term there are ${n\od 1}$ occurrences of
$c^k$. If ${n=1}$, then $c^k_{p_1,q_1}$ stands for $\mj_{p_1\vee
q_1}$. We write $c^k_{P_1,\ldots,P_n,Q_1,\ldots,Q_n}$ for the
arrow term obtained from $c^k_{p_1,\ldots,p_n,q_1,\ldots,q_n}$ by
substituting the form sets $P_i$ and $Q_i$ for $p_i$ and $q_i$
respectively.

We can then easily prove the following.

\prop{Lemma 13}{If ${f\!:X_1\vee X_2\str Y_1\kon\ldots\kon Y_n}$,
for ${n\geq 1}$, is ${(X_1,X_2)}$-splitting, then
\[
f=(f_1\kon\ldots\kon f_n)\cirk
c^k_{Y_1^{-X_2},\ldots,Y_n^{-X_2},Y_1^{-X_1},\ldots,Y_n^{-X_1}}
\]
for ${f_i\!:Y_i^{-X_2}\vee Y_i^{-X_1}\str Y_i}$ an
${(Y_i^{-X_2},Y_i^{-X_1})}$-splitting arrow term, where $1\leq
i\leq n$.}

This lemma is of particular interest when $Y_1,\ldots,Y_n$ are
\emph{all} the conjuncts of the target of $f$; i.e., $Y_i$ is not
of the form ${Y_i'\kon Y_i''}$. We say that in this case $Y_i$ is
a \emph{prime conjunct}. We define analogously a \emph{prime
disjunct} of a form set, just by replacing $\kon$ by $\vee$.

By the dual of Lemma~3 and by Lemma~7 we have the following lemma.

\prop{Lemma 14}{If ${f\!:X_1\vee X_2\str Y'\vee Y''}$ is an
${(X_1,X_2)}$-splitting arrow term such that $Y'$ is a prime
disjunct of ${Y'\vee Y''}$, then ${f=f'\vee f''}$ for
${f'\!:X'\str Y'}$ where either \begin{itemize}\item[] $X'$ is a
prime disjunct of $X_i$ for some ${i\in\{1,2\}}$ and
${f'=\mj_{X'}}$, or \item[] ${X'=X_1'\vee X_2'}$ for $X_i'$ being
a prime disjunct of $X_i$ for ${i\in\{1,2\}}$ and $f'$ being
${(X_1',X_2')}$-splitting.
\end{itemize}}

We define as follows a class $\cal N$ of arrow terms, for which we
say that they are in ${(X_1,X_2)}$-\emph{splitting normal form}.
For every $X$ we have that $\mj_X$ is in $\cal N$.

If for ${n\geq 2}$ we have that
$c^k_{P_1,\ldots,P_n,Q_1,\ldots,Q_n}$ is ${(X_1,X_2)}$-splitting
and that $f_1,\ldots,f_n$, which are not all identity arrows, are
in $\cal N$, then
\begin{itemize} \item[] $c^k_{P_1,\ldots,P_n,Q_1,\ldots,Q_n}$ is
in $\cal N$,\vspace{-1ex} \item[] ${f_1\vee\ldots\vee f_n}$ is in
$\cal N$,\vspace{-1ex} \item[] ${(f_1\kon\ldots\kon f_n)\cirk
c^k_{P_1,\ldots,P_n,Q_1,\ldots,Q_n}}$ is in $\cal N$ (provided the
composition is defined).
\end{itemize}
It is easy to see that an arrow term in ${(X_1,X_2)}$-splitting
normal form is ${(X_1,X_2)}$-splitting. Lemmata 13 and 14
guarantee that every ${(X_1,X_2)}$-splitting arrow term with the
source ${X_1\vee X_2}$ is equal to an arrow term in
${(X_1,X_2)}$-splitting normal form.

For the proof of next lemma we define the notion of
${(\kon,\vee)}$-\emph{arrow-shape}, analogous up to a point to the
notion of $\!\ks\!$-shape of Section~4. We have:
\begin{itemize}
\item[]$\koc$ is a ${(\kon,\vee)}$-arrow-shape;\vspace{-1ex}
\item[]if $M$ is a ${(\kon,\vee)}$-arrow-shape, then
${M\kon\mj_X}$ and ${M\vee\mj_X}$ are
${(\kon,\vee)}$-arrow-shapes.
\end{itemize}
For $M$ a ${(\kon,\vee)}$-arrow-shape and $f$ an arrow term, the
arrow term $M(f)$ is obtained from $M$ by replacing $\koc$ by $f$.

We will also use in the proof of the next lemma the following
abbreviations for ${m\geq 1}$:
\begin{tabbing}
\mbox{\hspace{11.5em}}\=$\bigwedge X_m$ \=$=_{df}X_1\kon\ldots\kon
X_m$.\kill

\>\>$\vec{X}_m\!$\'$=_{df}X_1,\ldots,
X_m$,\\[.5ex]
\> $\bigwedge X_m$ \> $=_{df}X_1\kon\ldots\kon X_m$.
\end{tabbing}

\prop{Lemma 15}{Let ${g\cirk h}$ be such that $g$ is a $c^k$-term
that is not ${(X_1,X_2)}$-splitting, and $h$ is an
${(X_1,X_2)}$-splitting arrow term whose source is ${X_1\vee
X_2}$. Then there exists a $c^k$-term $g'$ that is not
${(X_1,X_2)}$-splitting such that ${g\cirk h=f\cirk g'}$ for some
arrow term $f$.}

\dkz As remarked above, by Lemmata 13 and 14, we may assume that
$h$ is in ${(X_1,X_2)}$-splitting normal form. We proceed by
induction on the number $n$ of occurrences of $c^k$ in $h$. (We
find ${m\od 1}$ occurrences of $c^k$ in each
$c^k_{P_1,\ldots,P_m,Q_1,\ldots,Q_m}$.)

By the assumption that $g$ is not ${(X_1,X_2)}$-splitting and that
$h$ is in ${(X_1,X_2)}$-splitting normal form, we have only two
interesting cases, for which the following two cases are typical:
\begin{tabbing}
\mbox{\hspace{1em}}\= (1)\hspace{.5em}\=$g\cirk h=g\cirk h''\cirk
h'=M(c^k_{\Kon(P_m^1\vee P_m^2),\Kon(Q_n^1\vee
Q_n^2),R^2,S^2}\cirk$
\\*[.5ex]
\` $\cirk(c^k_{\vec{P}_m^1,\vec{Q}_n^1,\vec{P}_m^2,\vec{Q}_n^2}
\vee\mj_{R^2\kon S^2}))\cirk h'$,
\\[2ex]
\> (2) \> $g\cirk h=g\cirk h''\cirk h'=M(c^k_{\Kon(P_m^1\vee
P_m^2),\Kon(Q_n^1\vee Q_n^2), \Kon(R_l^1\vee R_l^2),\Kon(S_k^1\vee
S_k^2)} \cirk$
\\*[.5ex]
\` $\cirk(c^k_{\vec{P}_m^1,\vec{Q}_n^1,\vec{P}_m^2,\vec{Q}_n^2}
\vee c^k_{\vec{R}_l^1,\vec{S}_k^1,\vec{R}_l^2,\vec{S}_k^2}))\cirk
h'$
\end{tabbing}
where $M$ is a ${(\kon,\vee)}$-arrow-shape, and for every form set
$Y^i$ in whose name there is a superscript ${i\in\{1,2\}}$, we
have that ${let(Y^i)}\subseteq {let(X_i)}$.

In case (1) we have that
\begin{tabbing}
$c^k_{\Kon(P_m^1\vee P_m^2),\Kon(Q_n^1\vee Q_n^2),R^2,S^2}\cirk
(c^k_{\vec{P}_m^1,\vec{Q}_n^1,\vec{P}_m^2,\vec{Q}_n^2}
\vee\mj_{R^2\kon S^2})=$
\\[1ex]
\` $((c^k_{\vec{P}_m^1,\vec{P}_m^2}\vee\mj_{R^2}) \kon
(c^k_{\vec{Q}_n^1,\vec{Q}_n^2}\vee\mj_{S^2})\cirk c^k_{\Kon
P_m^1,\Kon Q_n^1,\Kon P_m^2\vee R^2,\Kon Q_n^2\vee S^2}\cirk$
\\*[.5ex]
\` $\cirk (\mj_{\Kon P_m^1 \kon \Kon Q_n^1}\vee c^k_{\Kon
P_m^2,\Kon Q_n^2,R^2,S^2})$
\end{tabbing}
by taking that
$c^k_{\vec{P}_m^1,\vec{Q}_n^1,\vec{P}_m^2,\vec{Q}_n^2}$ is
$(c^k_{\vec{P}_m^1,\vec{P}_m^2}\vee
c^k_{\vec{Q}_n^1,\vec{Q}_n^2})\cirk c^k_{\Kon P_m^1,\Kon
Q_n^1,\Kon P_m^2,\Kon Q_n^2}, $ and by applying the naturality
equation for $c^k$ and the equation $(\psi\check{b})$ of
Section~4.

Then we apply the induction hypothesis to the composition of
\[
M(\mj_{\Kon P_m^1\kon\Kon Q_n^1}\vee c^k_{\Kon P_m^2,\Kon
Q_n^2,R^2,S^2}),
\]
which is not ${(X_1,X_2)}$-splitting, with $h'$, which is in
${(X_1,X_2)}$-splitting normal form, and with at least one
occurrence of $c^k$ less than in $h$.

In case (2) we have that
\begin{tabbing}
$c^k_{\Kon(P_m^1\vee P_m^2),\Kon(Q_n^1\vee Q_n^2),\Kon(R_l^1\vee
R_l^2),\Kon(S_k^1\vee S_k^2)}\cirk
(c^k_{\vec{P}_m^1,\vec{Q}_n^1,\vec{P}_m^2,\vec{Q}_n^2} \vee
c^k_{\vec{R}_l^1,\vec{S}_k^1,\vec{R}_l^2,\vec{S}_k^2})=$
\\[1ex]
\` $((c^k_{\vec{P}_m^1,\vec{P}_m^2}\!\!\vee\!
c^k_{\vec{R}_l^1,\vec{R}_l^2}) \kon
(c^k_{\vec{Q}_n^1,\vec{Q}_n^2}\!\!\vee\!
c^k_{\vec{S}_k^1,\vec{S}_k^2}))\cirk c^k_{\Kon P_m^1 \vee \Kon
R_l^1,\Kon Q_n^1 \vee \Kon S_k^1,\Kon P_m^2 \vee \Kon R_l^2,\Kon
Q_n^2 \vee \Kon S_k^2}\cirk$
\\*[.5ex]
\` $\cirk (c^k_{\Kon P_m^1,\Kon Q_n^1,\Kon R_l^1,\Kon S_k^1}\vee
c^k_{\Kon P_m^2,\Kon Q_n^2,\Kon R_l^2,\Kon S_k^2})$
\end{tabbing}
by taking that
$c^k_{\vec{P}_m^1,\vec{Q}_n^1,\vec{P}_m^2,\vec{Q}_n^2}$ is as
above, and analogously for
$c^k_{\vec{R}_l^1,\vec{S}_k^1,\vec{R}_l^2,\vec{S}_k^2}$, and by
applying the naturality equation for $c^k$ and the equation
$(\psi\check{c}^m)$ from the beginning of this section.

Then we apply the induction hypothesis to the composition of
\[
M(\mj_{(\Kon P_m^1\kon\Kon Q_n^1)\vee(\Kon R_l^1\kon \Kon
S_k^1)}\vee c^k_{\Kon P_m^2,\Kon Q_n^2,\Kon R_l^2,\Kon S_k^2}),
\]
which is not ${(X_1,X_2)}$-splitting, with $h'$, which is in
${(X_1,X_2)}$-splitting normal form, and with at least two
occurrence of $c^k$ less than in $h$. \qed

Then we can prove the following key lemma.

\prop{Lemma 16}{For every arrow ${f\!:X_1\vee X_2\str Y}$ of
$\SCk^{st}$ there is an ${(X_1,X_2)}$-nonsplitting arrow term $f'$
and an ${(X_1,X_2)}$-splitting arrow term $f''$ such that
${f=f''\cirk f'}$.}

\dkz By the Development Lemma $f$ is equal to a developed arrow
term. Every developed arrow term of $\SCk^{st}$ such that some of
its factors are ${(X_1,X_2)}$-splitting $c^k$-terms is of the form
\[
f_{n+k+l}''\cirk\ldots\cirk f_{n+k+1}''\cirk
g_{n+k}\cirk\ldots\cirk g_{n+1}\cirk f_n''\cirk\ldots\cirk
f_1''\cirk f_m'\cirk\ldots\cirk f_1',
\]
for ${m,n\geq 1}$ and ${k,l\geq 0}$, where ${f_1',\ldots,f_m'}$,
$g_{n+1}$ and $g_{n+k}$ are ${(X_1,X_2)}$-non\-splitting, while
$f_1'',\ldots,f_n'',f_{n+k+1}'',\ldots,f_{n+k+l}''$ are
${(X_1,X_2)}$-splitting. (Note that identity arrow terms are
${(X_1,X_2)}$-nonsplitting, as well as ${(X_1,X_2)}$-splitting.)

By Lemma~5, the target of ${f_m'\cirk\ldots\cirk f_1'}$ is of the
form ${Y_1\vee Y_2}$ such that ${let(X_i)}\subseteq {let(Y_i)}$
for every ${i\in\{1,2\}}$. Hence being ${(X_1,X_2)}$-splitting is
the same as being ${(Y_1,Y_2)}$-splitting.

If ${k=0}$, then we are done. If ${k>0}$, then we call ${n+k}$ the
\emph{tail length}. By applying Lemma 15, the Development Lemma
and Lemma~4, we obtain that
\[
g_{n+1}\cirk f_n''\cirk\ldots\cirk f_1''=h_n\cirk\ldots\cirk
h_1\cirk g'
\]
for $g'$ a $c^k$-term that is not ${(X_1,X_2)}$-splitting. Then
after replacing above the left-hand side by the right-hand side we
have either obtained ${f=f''\cirk f'}$ as desired, or we have
obtained an arrow term with a strictly smaller tail length.
(Formally, we make an induction on a multiset ordering; see
\cite{DM79}.) \qed

We can now establish the following.

\prop{$\SCk^{st}$ Coherence}{If ${f',f''\!:X\str Y}$ are arrows of
$\SCk^{st}$, then ${f'=f''}$; i.e., the category $\SCk^{st}$ is a
diversified preorder.}

\dkz As in the proofs of the Theorems of Sections 10 and 11, we
proceed by induction on the sum $n$ of the number of letters in
the form set $X$ and the number of occurrences of $\kon$ in $X$.
If ${n=1}$, then ${X=Y=p}$, and ${\mj_p\!:p\str p}$ is the unique
arrow from $p$ to $p$.

If ${X=X_1\kon X_2}$, then, by Lemma~3, we have that $Y$ must be
of the form ${Y_1\kon Y_2}$ and ${f^l=f_1^l\kon f_2^l}$ for every
${l\in\{',''\}}$, where both $f_i'$ and $f_i''$ are of the type
${X_i\str Y_i}$ for ${i\in\{1,2\}}$. Then we may apply the
induction hypothesis to $f_i'$ and $f_i''$. We proceed analogously
when $Y$ is ${Y_1\vee Y_2}$ just by relying on the dual of
Lemma~3.

Suppose now that $X$ is ${X_1\vee X_2}$ and $Y$ is ${Y_1\kon
Y_2}$. Then, by Lemma 16 and Lemmata 5-6, we have that
${f^l=h^l\cirk(g_1^l\vee g_2^l)}$  for every ${l\in\{',''\}}$,
where $h^l$ is ${(X_1,X_2)}$-splitting (which is the same as being
${(Y^{-X_2},Y^{-X_1})}$-splitting) and $g_i^l$ is of the type
${X_i\str Y^{-X_{3-i}}}$ for every ${i\in\{1,2\}}$. If there is at
least one occurrence of $c^k$ in $g_i^l$ for some ${i\in\{1,2\}}$
and some ${l\in\{',''\}}$, then we may apply the induction
hypothesis to $g_i'$ and $g_i''$, since at least the number of
letters has decreased in their source, and we may apply the
induction hypothesis to $h'$ and $h''$, since in their source
${Y^{-X_2}\vee Y^{-X_1}}$ there is at least one occurrence of
$\kon$ less than in ${X_1\vee X_2}$.

If ${g_i'=g_i''=\mj_{X_i}}$ for ${i\in\{1,2\}}$, then, by Lemma
11, for every ${l\in\{',''\}}$
\[
f^l=u^l\cirk c^k_{Y_1^{-X_2},Y_2^{-X_2},Y_1^{-X_1},Y_2^{-X_1}}
\]
and we may apply the induction hypothesis to $u'$ and $u''$. \qed

From this result we obtain as a corollary the main result of this
section.

\prop{Symmetric Biassociative Intermuting Coherence}{The category
\SCk\ is a diversified preorder.}

\section{\large\bf Lattice categories and symmetric biassociative intermuting categories}
In this section we consider the relationship between the symmetric
biassociative intermuting categories of the preceding section and
an important type of categories for which coherence was previously
established. We will just summarize matters, and will not go into
all the details, known either from \cite{DP04}, or other, earlier,
references.

A \emph{lattice} category is a category with all finite nonempty
products and coproducts. We call \ml\ the free lattice category
generated by a set of objects. A detailed equational presentation
of \ml, in several possible languages, may be found in \cite{DP04}
(Chapter~9). The language on which we rely in this section has as
primitive arrow terms those of \SCk\ extended with the arrow terms
corresponding to the diagonal maps, the projections, and their
duals:
\begin{tabbing}
\hspace{5em}\=$\hat{w}_A\!:A\str A\kon
A$,\hspace{5em}\=$\hat{k}^i_{A_1,A_2}\!\!:A_1\kon A_2\str
A_i$,\\*[1ex] \>$\check{w}_A\!:A\vee A\str
A$,\>$\check{k}^i_{A_1,A_2}\!\!:A_i\str A_1\vee A_2$,
\end{tabbing}
for ${i\in\{1,2\}}$; the arrow terms are closed under composition
$\cirk$ and the operations $\kon$ and $\vee$. The equations
assumed for \ml\ are first categorial, bifunctorial and naturality
equations. Next we have equations that guarantee that $\kon$ is
binary product, with $\hat{w}$ and ${(k^1,k^2)}$ being
respectively the unit and counit of the underlying adjunction (see
the equations ${(\hat{w}\hat{k})}$ and ${(\hat{w}\hat{k}\hat{k})}$
in \cite{DP04}, Section 9.1). We also have dual equations that
guarantee that $\vee$ is binary coproduct. Next we have
definitional equations for $\hat{b}$, $\hat{c}$ and $c^k$:
\begin{tabbing}
\hspace{2em}\=$c^k_{A,B,C,D}\,$\=$=((\hat{k}^1_{A,B}\vee
\hat{k}^1_{C,D})\kon(\hat{k}^2_{A,B}\vee\hat{k}^2_{C,D}))\cirk\hat{w}_{(A\kon
B)\kon(C\kon D)}$\kill

\>\hspace{.9em}$\hat{b}^\str_{A,B,C}$\>$=((\mj_A\kon
\hat{k}^1_{B,C})\kon(\hat{k}^2_{B,C}\cirk\hat{k}^2_{A,B\kon
C}))\cirk\hat{w}_{A\kon(B\kon C)}$,\\*[1ex]
\>\hspace{.9em}$\hat{b}^\rts_{C,B,A}$\>$=((\hat{k}^1_{C,B}\cirk\hat{k}^1_{C\kon
B,A})\kon(\hat{k}^2_{C,B}\kon\mj_A))\cirk\hat{w}_{(C\kon B)\kon A}$,\\[1ex]
\>\hspace{1.7em}$\hat{c}_{A,B}$\>$=(\hat{k}^2_{A,B}\kon\hat{k}^1_{A,B})\cirk\hat{w}_{A\kon
B}$,\\[1ex]
\>$c^k_{A,B,C,D}$\>$=((\hat{k}^1_{A,B}\vee
\hat{k}^1_{C,D})\kon(\hat{k}^2_{A,B}\vee\hat{k}^2_{C,D}))\cirk\hat{w}_{(A\kon
B)\kon(C\kon D)}$,\\[1.5ex]
or, alternatively,\\[1.5ex]
\>$c^k_{A,B,C,D}$\>$=\check{w}_{(A\vee C)\kon(B\vee
D)}\cirk((\check{k}^1_{A,C}\kon\check{k}^1_{B,D})\vee
(\check{k}^2_{A,C}\kon\check{k}^2_{B,D}))$,
\end{tabbing}
and dual equations for $\check{b}$ and $\check{c}$.

There exists a faithful functor $G$ from the category \ml\ to the
category \emph{Rel} whose arrows are relations between finite
ordinals. The existence of this faithful functor is called Lattice
Coherence (see \cite{DP04}, Section 9.4). Symmetric Biassociative
Intermuting Coherence could also be expressed by stating that a
functor $G$ from \SCk\ to \emph{Rel}, which amounts to a
restriction of $G$ from \ml\ to \emph{Rel}, is a faithful functor.
(Coherence in the sense of preordering that we had previously in
this paper can also be expressed as the existence of a faithful
functor into \emph{Rel}, the image under this functor being a
discrete subcategory of \emph{Rel}; see \cite{DP04}, Section 2.9.)
From Symmetric Biassociative Intermuting Coherence it follows that
\SCk\ is isomorphic to a subcategory of \ml.

We can then prove that \SCk\ catches an interesting fragment of
\ml. (This result may be understood as extending the result of
\cite{DP99}.)

\prop{Proposition}{If for the arrow ${f\!:A\str B}$ of \ml\ we
have that ${G(f)}$ is a bijection, then there is an arrow term
${f'\!:A\str B}$ of \SCk\ such that ${f=f'}$ in \ml.}

\dkz By relying on the following equations of \ml:
\begin{tabbing}
\hspace{2.5em}\=$\hat{w}_{A\kon
B}\,$\=$=\hat{c}^m_{A,A,B,B}\cirk(\hat{w}_A\kon\hat{w}_B)$,\hspace{3em}\=
$\check{w}_{A\vee
B}\,$\=$=(\check{w}_A\vee\check{w}_B)\cirk\check{c}^m_{A,B,A,B}$,\\*[1ex]
\>$\hat{w}_{A\vee
B}\,$\>$=c^k_{A,A,B,B}\cirk(\hat{w}_A\vee\hat{w}_B)$,\>
$\check{w}_{A\kon B}\,$\>$=(\check{w}_A\kon\check{w}_B)\cirk
c^k_{A,B,A,B}$,
\end{tabbing}
we may assume that every $\w{\xi}{}$, for
${\!\ks\!\in\{\kon,\vee\}}$, has as its index a letter; i.e.\ we
have only $\w{\xi}{}_p$ for $p$ a letter (see the preceding
section for the definition of ${\c{\xi}{\!^m}}$). By relying on
the following equation of \ml:
\[
\hat{k}^1_{C,A\kon
B}=\hat{k}^1_{C,A}\cirk(\mj_C\kon\hat{k}^1_{A,B}),
\]
which follows from the naturality of $\hat{k}^1$, we may assume
that for every occurrence of $\hat{k}^1_{C,D}$ the second index
$D$ is either a letter or of the form ${D_1\vee D_2}$. With the
analogous equations for $\hat{k}^2$, $\check{k}^1$ and
$\check{k}^2$, we are allowed to make analogous assumptions
concerning $\hat{k}^2_{D,C}$, $\check{k}^1_{C,D}$ and
$\check{k}^2_{D,C}$.

Next we rely on bifunctorial and naturality equations, and the
following equations of \ml:
\[
\begin{array}{rl}
(\hat{k}^1_{A,D}\kon\mj_C)\cirk\hat{b}^\str_{A,D,C} \!\!\!\!\!& =
\mj_A\kon\hat{k}^2_{D,C},\\[1ex]
\hat{k}^1_{A\kon B,D}\cirk\hat{b}^\str_{A,B,D} \!\!\!\!\!& =
\mj_A\kon\hat{k}^1_{B,D},\\[1ex]
(\hat{k}^2_{D,B}\kon\mj_C)\cirk\hat{b}^\str_{D,B,C} \!\!\!\!\!& =
\hat{k}^2_{D,B\kon C},
\end{array}
\]
\[
\hat{k}^1_{A,D}\cirk\hat{c}_{D,A}\!=\hat{k}^2_{D,A},\hspace{1em}
\hat{k}^2_{D,A}\cirk\hat{c}_{A,D}\!=\hat{k}^1_{A,D},
\]
\[
\begin{array}{l}
\hat{k}^1_{A\vee B,D_1\vee D_2}\cirk
c^k_{A,D_1,B,D_2}\!=\hat{k}^1_{A,D_1}\vee\hat{k}^1_{B,D_2},\\[1ex]
\hat{k}^2_{D_1\vee D_2,A\vee B}\cirk
c^k_{D_1,A,D_2,B}\!=\hat{k}^2_{D_1,A}\vee\hat{k}^2_{D_2,B},
\end{array}
\]
\[
\hspace{1.8em}\hat{k}^i_{p,p}\cirk\hat{w}_p\!=\mj_p,
\]
together with the analogous equations involving $\hat{b}^\rts$ and
the dual equations involving $\check{k}^i$, $\check{b}^\str$,
$\check{b}^\rts$, $\check{c}$ and $\check{w}$, in order to
eliminate every occurrence of $\k{\xi}{i}$. (The first three
equations displayed above are related to the three equations
displayed in the proof of $\mka^\emptyset_{\top,\bot}$ Coherence
in Section~7.)

Since ${G(f)}$ is a bijection, no occurrence of $\k{\xi}{i}$ can
remain. That no occurrence of $\w{\xi}{}$ can remain after having
eliminated all the occurrences of $\k{\xi}{i}$ follows from the
Composition Elimination result of \cite{DP04} (Section 9.4)
combined with the existence of a functor $G'$ from \ml\ to the
category \emph{Mat}, which is isomorphic to the skeleton of the
category of finite-dimensional vector spaces over a fixed number
field with linear transformations as arrows (see \cite{DP04},
Section 12.5). For example, $c^k_{p,q,r,s}$ is mapped by $G$ to
the relation whose diagram is
\begin{center}
\begin{picture}(97,60)

\put(11,15){\line(0,1){30}}

\put(83,15){\line(0,1){30}}

\put(35,15){\line(2,3){20}}

\put(57,15){\line(-2,3){20}}

\put(5,5){$(p\;\vee\;r)\;\kon\;(q\;\vee\;s)$}

\put(5,50){$(p\;\kon\;q)\;\vee\;(r\;\kon\;s)$}

\end{picture}
\end{center}
while by $G'$ it is mapped to the matrix
\[
\left[%
\begin{array}{cccc}
  1 & 0 & 0 & 0 \\
  0 & 0 & 1 & 0 \\
  0 & 1 & 0 & 0 \\
  0 & 0 & 0 & 1 \\
\end{array}%
\right]
\]
The category \emph{Mat} takes into account whether in the diagram
corresponding to ${G(f)}$ two occurrences of the same letter may
be joined by more than one line.\qed

\section{\large\bf Symmetric bimonoidal intermuting categories}
In this section, which is parallel to Section 12, we establish as
our final result a symmetric variant of Restricted Bimonoidal
Intermuting Coherence. This result is based essentially on
Symmetric Biassociative Intermuting Coherence of Section 14.
Before defining symmetric bimonoidal intermuting categories, we
introduce some preliminary notions of bimonoidal categories with
natural commutativity isomorphisms, and prove auxiliary coherence
results for them.

A \emph{symmetric bimonoidal} category is a bimonoidal category
${\langle\aA,\kon,\vee,\top,\bot\rangle}$ (see Section~2) such
that ${\langle\aA,\kon,\vee\rangle}$ is a symmetric biassociative
category (see Section 14). Symmetric bimonoidal categories are
coherent in the sense that the symmetric bimonoidal category
freely generated by a set of objects is a diversified preorder
(see \cite{DP04}, Section 6.4).

In \cite{DP06a} one can find a justification in the spirit of
Sections 3 and 4 of all the assumptions made for symmetric
bimonoidal categories. In these categories $\!\ks\!$, for
${\!\ks\!\in\{\kon,\vee\}}$, intermutes with itself, and
$\b{\xi}{\str}$, $\c{\xi}$, $\d{\xi}{\str}$ and $\s{\xi}{\str}$
are upward and downward preserved by $\!\ks\!$, where this
preservation is understood quite analogously to what we had in
Section~4. Mac Lane's pentagonal and hexagonal equations follow
from this preservation. The role of $c^k$ in that is played by the
natural isomorphisms ${\c{\xi}{\!^m}}$ of Section 14. In terms of
these isomorphisms and of the unit objects one defines the
$b$-arrows and the $c$-arrows.

A \emph{symmetric normal bimonoidal} category is a symmetric
bimonoidal category that is also normal bimonoidal (see
Section~7), and that satisfies moreover the two equations
\begin{tabbing}
\hspace{11em}\=${(\p{\top}\!\check{c})}$\quad\quad\=$\check{w}^\str_\top\!\cirk
\check{c}_{\top,\top}\,$\=$=\check{w}^\str_\top$,\\*[1ex]
\>${(\pp{\bot}\!\hat{c})}$\quad\quad\=$\hat{c}_{\bot,\bot}\!\cirk\hat{w}^\rts_\bot
\,$\=$=\hat{w}^\rts_\bot$,
\end{tabbing}
called collectively ${(wc)}$. The equation
${(\p{\top}\!\check{c})}$ says that $\check{c}$ is upward
preserved by $\top$, and ${(\pp{\bot}\!\hat{c})}$ says that
$\hat{c}$ is downward preserved by $\bot$ (see Section~4). From
these two equations we obtain immediately the equations
\begin{tabbing}
\hspace{7em}\=${(c\mj)}$\quad\quad$\check{c}_{\top,\top}=\mj_{\top\vee\top}$,\quad\quad\quad
$\hat{c}_{\bot,\bot}=\mj_{\bot\kon\bot}$.
\end{tabbing}

We call $\mns_{\top,\bot}$ the symmetric normal bimonoidal
category freely generated by a set of objects, and we can prove
the following.

\prop{Symmetric Normal Bimonoidal Coherence}{The category
$\mns_{\top,\bot}$ is a diversified preorder.}

\dkz By Symmetric Biassociative Coherence (see the beginning of
Section 14) and the results of \cite{DP04} (Chapter~3, in
particular Section 3.3), we can replace the category
$\mns_{\top,\bot}$ by a strictified category
$\mns_{\top,\bot}^{st}$ where the $b$-arrows and the $c$-arrows
are identity arrows. We can show that $\mns_{\top,\bot}^{st}$ is a
diversified preorder, in the sense in which $\SCk^{st}$ is a
diversified preorder (see the beginning of Section 14). For that
we proceed as for the proof of Normal Biunital Coherence in
Section~6. This will imply that $\mns_{\top,\bot}$ is a
diversified preorder. \qed

A \emph{symmetric $\kappa$-normal bimonoidal} category is a
symmetric normal bimonoid\-al category that is also
$\kappa$-normal bimonoidal (see Section~7). We call
$\mks^\emptyset_{\top,\bot}$ the symmetric $\kappa$-normal
bimonoidal category freely generated by the empty set of objects.
We can prove the following.

\prop{$\mks^\emptyset_{\top,\bot}$ Coherence}{The category
$\mks^\emptyset_{\top,\bot}$ is a preorder.}

\dkz We enlarge the proofs of $\mk^\emptyset_{\top,\bot}$
Coherence and $\mka^\emptyset_{\top,\bot}$ Coherence of Sections 6
and 7 by using the following equations of symmetric monoidal
categories:
\[
\c{\xi}_{\zeta,A}\,=\;\d{\xi}{\rts}_A\!\cirk\s{\xi}{\str}_A,\quad\quad\quad
\c{\xi}_{A,\zeta}\,=\;\s{\xi}{\rts}_A\!\cirk\d{\xi}{\str}_A,
\]

\vspace{1ex}

\noindent for ${(\!\ks\!,\zeta)\in\{(\kon,\top),(\vee,\bot)\}}$.
With these equations and the equations ${(c\mj)}$ we can eliminate
every occurrence of $c$. \qed

A \emph{symmetric bimonoidal intermuting} category is a symmetric
$\kappa$-normal bimonoidal category
${\langle\aA,\kon,\vee,\top,\bot\rangle}$ that is also a
bimonoidal intermuting category (see Section 12) such that
${\langle\aA,\kon,\vee\rangle}$ is a symmetric biassociative
intermuting category (see Section~14). This means that we assume
all the equations we have assumed for various notions of
categories considered up to now, excluding the lattice categories
of the preceding section. In addition to the equations listed for
bimonoidal intermuting categories at the beginning of Section 12
we assume the equations ${(c^kc)}$ of Section 14 and the equations
${(wc)}$ from the beginning of this section.

The $b$-arrows and the $c$-arrows are definable in terms of the
isomorphisms ${\c{\xi}{\!^m}}$ and of the unit objects, as we
mentioned above. We could moreover replace the equations
${(c^kb)}$ and ${(c^kc)}$ by the equations ${(c^kc^m)}$ of Section
14. According to \cite{DP06a}, the equation ${(\psi c^m)}$ (see
Section 14), and an analogous equation with $\kon$ replaced by
$\vee$, deliver Mac Lane's pentagonal and hexagonal equations. We
have mentioned in Section 14 that the equations ${(c^kc^m)}$ and
${(\psi c^m)}$ may be viewed as instances of the hexagonal
interchange equation of \cite{BFSV} (end of Definition 1.7), which
we have called \emph{HI}, and which serves to define $n$-fold
monoidal categories for $n>2$. In an $n$-fold version of our
notion of symmetric bimonoidal intermuting category, for $n\geq
2$, the scheme of the equation \emph{HI} would deliver all the
coherence conditions, save those involving the unit objects.

Let $\mbox{\bf SC}^{{\bf k}}_{\top,\bot}$ and $\mbox{\bf SC}^{{\bf
k}\emptyset}_{\top,\bot}$ be the free symmetric bimonoidal
intermuting categories generated respectively by a nonempty set of
objects and the empty set of objects. By combining the proof of
$\mbox{\bf AC}^{{\bf k}\emptyset}_{\top,\bot}$ Coherence of
Section 12 and the proof of $\mks^\emptyset_{\top,\bot}$ Coherence
above we obtain a proof of $\mbox{\bf SC}^{{\bf
k}\emptyset}_{\top,\bot}$ \emph{Coherence}, which says that
$\mbox{\bf SC}^{{\bf k}\emptyset}_{\top,\bot}$ is a preorder.

By proceeding as in Section 12, and by relying on Symmetric
Biassociative Intermuting Coherence of Section 14, we can then
establish the following.

\vspace{2ex}

\noindent{\sc Restricted Symmetric Bimonoidal Intermuting
Coherence.}

\noindent {\it If $f,g\!:A\str B$ are arrows of $\mbox{\bf
SC}^{{\bf k}}_{\top,\bot}$ such that both $A$ and $B$ are either
pure and diversified, or no letter occurs in them, then ${f=g}$ in
$\mbox{\bf SC}^{{\bf k}}_{\top,\bot}$.}

\vspace{3ex}

\noindent {\small {\it Acknowledgement$\,$}. Work on this paper
was supported by the Ministry of Science of Serbia (Grants 144013
and 144029).}

\end{document}